\documentclass{amsart}
\usepackage{amssymb}
\usepackage{amsmath}
\usepackage{amsthm}
\usepackage{mathrsfs}
\usepackage{colortbl}

\newtheorem{theorem}{Theorem}
\newtheorem{lemma}[theorem]{Lemma}
\newtheorem{assumption}[theorem]{Assumption}
\newtheorem{remark}[theorem]{Remark}

\theoremstyle{definition}
\newtheorem{definition}[theorem]{Definition}

\newtheorem{proposition}[theorem]{Proposition}

\makeatletter

\@addtoreset{theorem}{section}
\makeatother

\makeatletter

\@addtoreset{equation}{section}
\makeatother

\newcommand{\dual}[1]{\langle #1 \rangle }
\begin{document}                                                 
\title[On $L^{3,\infty}$-stability of the Navier-Stokes system]{On $L^{3,\infty}$-stability of the Navier-Stokes system in exterior domains}                                 
\author[Hajime Koba]{Hajime Koba}                                
\address{Graduate School of Engineering Science, Osaka University\\
1-3 Machikaneyamacho, Toyonaka, Osaka, 560-8531, Japan}                                  
\email{iti@sigmath.es.osaka-u.ac.jp}                                      
\date{}                                       
\begin{abstract}
This paper studies the stability of a stationary solution of the Navier-Stokes system with a constant velocity at infinity in an exterior domain. More precisely, this paper considers the stability of the Navier-Stokes system governing the stationary solution which belongs to the weak $L^3$-space $L^{3,\infty}$. Under the condition that the initial datum belongs to a solenoidal $L^{3 , \infty}$-space, we prove that if both the $L^{3,\infty}$-norm of the initial datum and the $L^{3,\infty}$-norm of the stationary solution are sufficiently small then the system admits a unique global-in-time strong $L^{3,\infty}$-solution satisfying both $L^{3,\infty}$-asymptotic stability and $L^\infty$-asymptotic stability. Moreover, we investigate $L^{3,r}$-asymptotic stability of the global-in-time solution. Using $L^p$-$L^q$ type estimates for the Oseen semigroup and applying an equivalent norm on the Lorentz space are key ideas to establish both the existence of a unique global-in-time strong (or mild) solution of our system and the stability of our solution.
\end{abstract}       
\maketitle

\tableofcontents

\section{Introduction}\label{sect1}
\subsection{Purposes} 
This paper has two purposes. The first one is to construct a unique global-in-time strong $L^{3,\infty}$-solution of the Navier-Stokes system with a constant velocity at infinity. The second one is to investigate the stability of global-in-time $L^{3,\infty}$-solutions of the system when the initial datum is in an intermediate space between the $L^{3.1}$-space and the $L^{3,\infty}$-space. Here $L^{q,r}$ denotes the Lorentz space.

Let $\Omega$ be an exterior domain with smooth boundary in $\mathbb{R}^3$. Throughout this paper we fix $\Omega$. We are concerned with the stability of the Navier-Stokes system with a constant velocity at infinity:
\begin{equation}\label{eq11}
\begin{cases}
\partial_t u - \Delta u +(u , \nabla ) u + \nabla \Pi = \nabla \cdot F & \text{ in } \Omega \times (0, \infty ),\\
\nabla \cdot u=0 & \text{ in } \Omega \times (0, \infty ),\\
u|_{\partial \Omega} =0, { \ }\displaystyle{\lim_{|x| \to \infty}u=u_\infty}, { \ }u|_{t = 0} = u_0,
\end{cases}
\end{equation}
where the unknown function $u = u(x,t) = (u^1,u^2,u^3)$ is the velocity of the fluid, the unknown function $\Pi = \Pi (x,t)$ is the pressure of the fluid, while the given function $F = F (x) = (F_{jk}(x))_{j,k=1,2,3}$ is the external force, the given constant $u_\infty = (u_\infty^1, u_\infty^2, u_\infty^3 ) \in \mathbb{R}^3$ is the velocity of the fluid at infinity, and $u_0=u_0(x)$ is the given initial datum. Here we use the convention: $\Delta :=\partial_1^2 + \partial_2^2 + \partial_3^2$ and $\nabla := (\partial_1 , \partial_2 , \partial_3 )$. Note that $\nabla \cdot F = (\sum_{k=1}^3 \partial_k F_{jk})_{j=1,2,3}$. The model \eqref{eq11} illustrates the motion of an incompressible viscous fluid past an obstacle.

Putting $\tilde{u} = \tilde{u}(x,t) = ( \tilde{u}^1, \tilde{u}^2 , \tilde{u}^3 ) : =u (x,t) - u_\infty$, we have
\begin{equation}\label{eq12}
\begin{cases}
\partial_t \tilde{u} - \Delta \tilde{u} + (u_\infty , \nabla ) \tilde{u} + (\tilde{u} , \nabla ) \tilde{u} + \nabla \Pi = \nabla \cdot F & \text{ in } \Omega \times (0,\infty ),\\
\nabla \cdot \tilde{u}=0 & \text{ in } \Omega \times (0,\infty ),\\
\tilde{u}|_{\partial \Omega} = - u_\infty, { \ }\displaystyle{\lim_{|x| \to \infty} \tilde{u}=0}, { \ }\tilde{u}|_{t = 0} = u_0 - u_\infty.
\end{cases}
\end{equation}
To study the system \eqref{eq12}, we consider the stationary Navier-Stokes equations:
\begin{equation}\label{eq13}
\begin{cases}
- \Delta w + (u_\infty , \nabla ) w + (w , \nabla ) w + \nabla \pi = \nabla \cdot F &\text{ in } \Omega,\\
\nabla \cdot w=0 &\text{ in } \Omega ,\\
w|_{\partial \Omega} = - u_\infty, { \ }\displaystyle{\lim_{|x| \to \infty }w (x)} = 0.
\end{cases}
\end{equation}
Here $w=w(x) = (w^1,w^2,w^3)$ and $\pi = \pi (x)$. This paper studies the stability of solutions to the stationary Navier-Stokes equations under the three conditions:
\begin{assumption}\label{assA}
The function $w$ satisfies $\nabla \cdot w =0$ in $\Omega$ and
\begin{equation*}
\| w \|_{L^{3,\infty}} + \| w \|_{L^{3 \beta /(3 - \beta), \infty}} + \| \nabla w \|_{L^{\beta, \infty}} < + \infty \text{ for some }3/2 < \beta < 3.
\end{equation*}
\end{assumption}
\begin{assumption}\label{assB}
The function $w$ satisfies Assumption \ref{assA} and
\begin{equation*}
\| \nabla w \|_{L^2} < + \infty .
\end{equation*}
\end{assumption}
\begin{assumption}\label{assC}
The function $w$ satisfies Assumption \ref{assA} and
\begin{equation*}
\| w \|_{L^\infty} + \| \nabla w \|_{L^{3,\infty}}  + \| \nabla w \|_{L^2}< + \infty .
\end{equation*}
\end{assumption}
\noindent Here $\| \cdot \|_{L^{q,r}}$ denotes the Lorentz norm (see Subsection \ref{subsec21}). Remark that there exists a solution $(w, \pi)$ of the system \eqref{eq13} satisfying Assumption \ref{assA}, Assumption \ref{assB}, or Assumption \ref{assC} if $|u_\infty|$ and $F$ are sufficiently small in a suitable function space. See Finn \cite{Fin59}, Kozono-Yamazaki \cite[Main Theorem]{KY98A}, Shibata-Yamazaki \cite[Theorems 2.2 and 4.1]{SY05}, and Kim-Kozono \cite[Theorem 1.1]{KK12}. Many researchers have been studying both the existence and the uniqueness of solutions to the system \eqref{eq13} since Leray \cite{Ler33} and Finn \cite{Fin59}. See \cite{KY98A}, \cite{SY05}, \cite{Gal11}, \cite{KK12}, \cite{HKK13}, and the references given there.

To study the stability of a solution $(w , \pi )$ of the system \eqref{eq13}, we set $v = v(x,t) = (v^1,v^2,v^3) := \tilde{u}(x,t)-w(x)$ and $\mathfrak{p}= \mathfrak{p}(x,t) := \Pi (x,t) - \pi (x)$. It is easy to check that $(v , \mathfrak{p})$ formally satisfies the system: 
\begin{equation}\label{eq14}
\begin{cases}
\partial_t v - \Delta v + (u_\infty, \nabla ) v + (v, \nabla ) v + (v, \nabla ) w +(v , \nabla ) w + \nabla \mathfrak{p} = 0,\\
\nabla \cdot v =0 { \ \ \ \ \ \ \ \ \ \ \ \ \ \ \ \ \ \ \ \ \ \ \ \ \ \ \ \ \ \ \ \ \ \ \ \ \ \ \ \ \ \ \ \ \ \ \ \ \ }\text{ in } \Omega \times (0,\infty),\\
v |_{\partial \Omega} =0, { \ }\displaystyle{\lim_{|x| \to \infty}v=0}, { \ } v|_{t=0} = v_0.
\end{cases}
\end{equation}
Here $v_0 := u_0 - u_\infty - w$. This paper studies the system \eqref{eq14} instead of the system \eqref{eq11}. Applying the Helmholtz projection $P$ into the system \eqref{eq14}, we have
\begin{equation}\label{eq15}
\begin{cases}
v_t + L v + P \{ (v , \nabla ) v +  ( v , \nabla) w + (w , \nabla ) v \} = 0, { \ \ \ }t>0,\\
v|_{t=0} = P v_0.
\end{cases}
\end{equation}
Here $L v = P \{ -\Delta v + (u_\infty, \nabla ) v \}$. We call the operator $L$ the \emph{Oseen operator}. In particular, we call $L$ the \emph{Stokes operator} if $u_\infty =0$. See Subsection \ref{subsec22} for the Helmholtz projection $P$ and Subsections \ref{subsec23}-\ref{subsec26} for the Stokes and the Oseen operators.

Heywood \cite{Hey70}, \cite{Hey72} first studied the stability of $L^2$-solutions to the system \eqref{eq15} under the conditions that $w \in W^{2,2}$ and $\sup \{|x| |w(x)| \}$ is sufficiently small. He applied the Galerkin method and an energy inequality to show the existence of a unique global-in-time strong $L^2$-solution of \eqref{eq15} with the property that for each $\Omega' \subset \Omega$, $\| v (t)\|_{L^2 ( \Omega' )} \to 0$ as $t \to \infty$. {\color{red}{Masuda \cite{Mas75} considered the stability of a weak solution of \eqref{eq15} when $\nabla w \in L^3$ and $\sup \{|x| |w(x)| \}$ is sufficiently small. Masuda \cite{Mas75} made use of fractional powers of the Stokes operator and an energy inequality to show the existence of a weak solution of \eqref{eq15} satisfying $L^\infty$-asymptotic stability. Maremonti \cite{Mar84} constructed a global-in-time $L^2$-solution, satisfying $\| v (t ) \|_{L^\infty} \leq C t^{-1/2}$ as $t \to \infty$, of \eqref{eq15} when $u_\infty =0$, $w \in L^6 \cap W^{1,3} \cap \dot{W}^{1,p}(p>3)$, and the Reynold number is sufficiently small. Maremonti \cite{Mar85} considered $L^2$-asymptotic stability of a $L^2$-solution of \eqref{eq15}, and investigated $L^2$-decay rate with respect to time of the solution when the Reynold number is sufficiently smal and the initial datum belongs to $L^q(1 \leq q < 2)$. Miyakawa and Sohr \cite{MS88} constructed a weak solution of \eqref{eq15} satisfying the strong energy inequality to derive $L^2$-asymptotic stability of the solution when $\nabla w \in L^3$ and $\sup \{|x| |w(x)| \}$ is sufficiently small. Note that $w \in L^{3, \infty}$ if $\sup \{|x| |w (x)| \}$ is finite. Koba \cite{Kob13} applied maximal $L^p$-regularity for Hilbert space-valued functions to investigate $L^2$-asymptotic stability of energy solutions to the generalized Navier-Stokes-Boussinesq system including \eqref{eq15}. See also Koba \cite{Kob14}.}}

This paper studies $L^{3,\infty}$-asymptotic stability of $L^{3,\infty}$-solutions to the system \eqref{eq15}. Let us now introduce $L^n$-asymptotic stability of the Navier-Stokes system (i.e. \eqref{eq15} when $u_\infty = 0$ and $w \equiv 0$). Kato \cite{Kat84} used $L^p$-$L^q$ estimates for the heat kernel to establish the existence of a unique global-in-time $L^n$-solution of the Navier-Stokes system in the whole space $\mathbb{R}^n$ satisfying $L^n$-asymptotic stability when the $L^n$-norm of the initial datum is sufficiently small. Kozono \cite{Koz89} applied the Stokes semigroup and an implicit function theorem for Banach spaces to show the existence of a unique global-in-time $L^n$-solution of the Navier-Stokes system in the halfspace $\mathbb{R}^n_+$ satisfying $L^n$-asymptotic stability when the initial datum is sufficiently small in the $L^n$-norm. Iwashita \cite{Iwa89} obtained stability results similar to those in \cite{Kat84} on the Navier-Stokes system in an exterior domain by studying $L^p$-$L^q$ estimates for the Stokes semigroup. Remark that $L^n \subset L^{n,\infty}$ (see Lemma \ref{lem22}).

In this paper we show the existence of a unique global-in-time strong $L^{3,\infty}$-solution, satisfying the asymptotic stability, of the system \eqref{eq15}. To this end, we construct a unique global-in-time (generalized) mild $L^{3,\infty}$-solution of \eqref{eq15} in a Banach space $X_3$, and investigate the asymptotic stability of the solution. Before stating main results, we introduce some notation.
\begin{definition}[Solenoidal spaces]\label{def13}
For $1 < q < \infty$ and $1 \leq r \leq \infty$,
\begin{align*}
&C^\infty_{0, \sigma} = C^\infty_{0 , \sigma} (\Omega) := \{ f = (f^1,f^2,f^3 ) \in [C_0^\infty (\Omega)]^3 ; { \ } \nabla \cdot f = 0 \},\\
& L^q_\sigma =L^q_\sigma (\Omega) := \overline{C_{0,\sigma}^\infty (\Omega)}^{\| \cdot \|_{L^q}},{ \ \ }L^{q,\infty}_{0,\sigma} = L^{q,\infty}_{0,\sigma}(\Omega) := \overline{C_{0,\sigma}^\infty (\Omega)}^{\| \cdot \|_{L^{q,\infty}}},\\
& L^{q,r}_\sigma = L^{q,r}_\sigma (\Omega ) := (L^{q_1}_\sigma ( \Omega ) , L^{q_2}_\sigma (\Omega ) )_{\theta , r},
\end{align*}
where $1 < q_1 < q < q_2 < \infty$ and $0 < \theta < 1$ such that $1/q = (1 - \theta )/{q_1} + \theta / {q_2}$. Here $(\cdot , \cdot)_{\theta ,r}$ is the real interpolation couple. See \cite{BL76} and \cite{Lun09} for the interpolation theory. See Subsection \ref{subsec22} for the solenoidal spaces $L^{q,\infty}_{0,\sigma} ( \Omega )$ and $L^{q,r}_\sigma  ( \Omega )$. Note that $L^{q, \infty}_{0, \sigma} ( \Omega ) \subset L^{q,\infty}_\sigma ( \Omega )$.
\end{definition}
\begin{definition}[Function space $X_p$]
Let $3 \leq p < \infty$. Define
\begin{equation*}
X_p:= \{ f \in BC((0, \infty ) ; L^{3,\infty}_{\sigma } (\Omega ) ) \cap C ((0, \infty ) ; L^{p, \infty}_\sigma ( \Omega )); { \ }\| f \|_{X_p} < \infty  \}.
\end{equation*}
Here $\displaystyle{\| f \|_{X_p} := \sup_{t >0} \{ \| f (t) \|_{L^{3,\infty}} \} + \sup_{t>0} \{ t^{\frac{p - 3}{2 p}} \| f (t) \|_{L^{p,\infty}} \} }$.\\
Note that $X_3 = BC((0,\infty) ; L^{3,\infty}_\sigma ( \Omega ) )$.
\end{definition}
\begin{definition}[Duality pairing]
Let $1 \leq q , q' ,r,r' \leq \infty$ such that $1/q+1/{q'}=1$ and $1/r + 1/{r'}=1$, where $1/\infty:= 0$. By $\dual{\cdot , \cdot}$, we define the paring between $L^{q,r}$ and $L^{q',r'}$ as follows: for all $f \in [L^{q,r}( \Omega )]^3$ and $g \in [L^{q',r'}( \Omega )]^3$
\begin{equation*}
\dual{f, g} := \int_\Omega f(x) \cdot g( x) d x.
\end{equation*}
Here $L^{1, \ell} ( \Omega ) : = L^1 ( \Omega)$ and $L^{\infty, \ell} (\Omega ) := L^\infty (\Omega )$ for $1 \leq \ell \leq \infty$.
\end{definition}
\begin{definition}[Generalized mild solutions]{ \ }\\
Assume that $v_0 \in L_{\sigma}^{3 , \infty} ( \Omega )$ and that $v \in X_3$. We call the function $v$ a \emph{global-in-time generalized mild solution} of the system \eqref{eq15} with the initial datum $v_0$ if the following two properties hold:
\begin{enumerate}
\renewcommand{\labelenumi}{\rm{(\roman{enumi})}}
\item for every $t >0$ and $\psi \in L^{3/2 , 1}_\sigma (\Omega)$
\begin{equation*}
\dual{v (t) , \psi} = \dual{ v_0 , \mathrm{e}^{- t L^*} \psi} + \int_0^t \dual{v \otimes v(s) + v(s) \otimes w + w \otimes v(s) , \nabla \mathrm{e}^{-(t-s) L^*} \psi} d s ,
\end{equation*}
\item for each $\psi \in L^{3/2,1}_\sigma ( \Omega )$
\begin{equation*}
\lim_{t \to 0 + 0} \dual{v (t) , \psi } = {\color{red}{\dual{v_0 , \psi}}}.
\end{equation*}
\end{enumerate}
Here $L^* f = P \{ - \Delta f - (u_\infty , \nabla ) f \}$ and $\mathrm{e}^{- t L^*}$ is the semigroup whose generator is the dual operator $-L^*$ of the Oseen operator $- L$.
\end{definition}
Kozono and Ogawa \cite{KOK94} introduced the notion of the generalized mild solution.
\begin{definition}[Mild solutions]
Assume that $v_0 \in L_{\sigma}^{3 , \infty} ( \Omega )$ and that $v \in X_3$. We call the function $v$ a \emph{global-in-time mild solution} of the system \eqref{eq15} with the initial datum $v_0$ if the following two properties hold:
\begin{enumerate}
\renewcommand{\labelenumi}{\rm{(\roman{enumi})}}
\item for every $t >0$
\begin{equation*}
v (t) = \mathrm{e}^{-t L} v_0- \int_0^t \mathrm{e}^{- (t - s)L} P\{ (v(s) , \nabla ) v(s) + (v(s) , \nabla )w + (w , \nabla )v(s) \} d s,
\end{equation*}
\item $v ( t ) \to v_0$ in $L^{3 ,\infty} (\Omega)$ as $t \to 0 + 0$.
\end{enumerate}
\end{definition}
\begin{definition}[Strong solutions]
Assume that $v_0 \in L_{\sigma}^{3 , \infty} ( \Omega )$ and that $v \in X_3$. We call the function $v$ a \emph{global-in-time strong solution} of the system \eqref{eq15} with the initial datum $v_0$ if the following two properties hold:
\begin{enumerate}
\renewcommand{\labelenumi}{\rm{(\roman{enumi})}}
\item $v$ is a global-in-time mild solution of the system \eqref{eq15} with the initial datum $v_0$,
\item for each $0 < t < \infty$, $v (t)$ satisfies the system \eqref{eq15} in $L^{3, \infty}$.
\end{enumerate}
\end{definition}
\subsection{Main Results and Key Ideas} 
We now state four main results.
\begin{theorem}[Existence of a unique generalized mild solution]\label{thm19}{ \ }\\
Assume that $w$ is as in Assumption \ref{assA} and that $|u_\infty| \leq \gamma$ for some $\gamma > 0$. Let $3 < p  <\infty$. Then there are two positive constants $\delta_0 = \delta_0 ( \gamma , p )$ and $K_0 = K_0 ( \gamma , p ) >0$ such that if
\begin{equation*}
v_0 \in L_{\sigma}^{3,\infty} (\Omega) , { \ }\| w \|_{L^{3,\infty}} < \delta_0, { \ }\text{ and }{ \ }\| v_0 \|_{L^{3,\infty}} < \delta_0 ,
\end{equation*}
then there exists a unique global-in-time generalized mild solution $v$ in $X_3$ of the system \eqref{eq15} with the initial datum $v_0$, satisfying $v \in X_p$ and
\begin{equation}\label{eq16}
\| v \|_{X_p} \leq K_0 \| v_0 \|_{L^{3 , \infty} }. 
\end{equation}
\end{theorem}
\begin{theorem}[Stability of generalized mild solutions]\label{thm110}{ \ }\\
Assume that $w$ is as in Assumption \ref{assA} and that $|u_\infty| \leq \gamma$ for some $\gamma > 0$. Let $3 < p  <\infty$, and let $\delta_0$ the constant appearing in Theorem \ref{thm19}. Suppose that
\begin{equation*}
v_0 \in L^{3, \infty}_\sigma (\Omega),{ \ } \| w \|_{L^{3,\infty}} < \delta_0, { \ }\text{ and }{ \ }\| v_0 \|_{L^{3,\infty}} < \delta_0 .
\end{equation*}
Let $v$ be the global-in-time generalized mild solution $v$ of the system \eqref{eq15} with the initial datum $v_0$, obtained by Theorem \ref{thm19}. Let $3p/(p + 3) < \alpha < 3$. Then there is $\delta_1 = \delta_1 ( \gamma , p , \alpha) >0$ such that if
\begin{equation*}
\| w \|_{L^{3,\infty}} < \delta_1 { \ }\text{ and }{ \ }\| v_0 \|_{L^{3,\infty}} < \delta_1 ,
\end{equation*}
then the following three assertions hold:\\
\noindent $(\mathrm{i})$ Assume in addition that $v_0 \in L^{\alpha , \infty}_\sigma (\Omega)$. Then
\begin{align}
\sup_{t>0} \{ \| v (t) \|_{L^{\alpha, \infty}} \} \leq \text{ Const.} < + \infty, \notag\\
\sup_{t>0} \{ t^{ \frac{3}{2} \left( \frac{1}{\alpha} -\frac{1}{p} \right) } \| v ( t ) \|_{L^{p , \infty}} \} \leq \text{ Const.} < + \infty,\\
\sup_{t>0} \{ t^{ \frac{3}{2} \left( \frac{1}{\alpha} -\frac{1}{3} \right) }  \| v (t) \|_{L^{3,\infty}} \} \leq \text{ Const.}  < + \infty .
\end{align}
\noindent $(\mathrm{ii})$ Assume in addition that $v_0 \in L^{3,\infty}_{0 , \sigma } ( \Omega )$ and that $\| \nabla w \|_{L^2} < + \infty$. Then
\begin{align}
\lim_{t \to 0 +0} \| v (t) - v_0 \|_{L^{3,\infty}} =0, \notag\\
\lim_{t \to \infty} \| v (t ) \|_{L^{3,\infty}} = 0 .\label{eq18}
\end{align}
\noindent $(\mathrm{iii})$ Assume in addition that $v_0 \in L^{3, r }_{\sigma } ( \Omega )$ for some $1 < r <\infty$ and that $\| \nabla w \|_{L^2} < + \infty$. Then
\begin{align}
v \in BC ( [ 0 ,  \infty ) ; L^{3,r}_\sigma ( \Omega ) ),\notag \\ 
\lim_{t \to 0 +0} \| v (t) - v_0 \|_{L^{3, r}} =0, \notag \\
\lim_{t \to \infty} \| v (t ) \|_{L^{3,r}} = 0 .
\end{align}
\end{theorem}
\begin{theorem}[Unique mild solution and stability]\label{thm111}{ \ }\\
Assume that $w$ is as in Assumption \ref{assB} and $|u_\infty| \leq \gamma$ for some $\gamma > 0$. {\color{red}{Let $6 < p < \infty$ such that ${3 \beta}/(2 \beta -3)< p$}}. Then then there is $\delta_2 = \delta_2 ( \gamma , p ) > 0$ such that if
\begin{equation*}
v_0 \in L_{0 , \sigma}^{3,\infty} (\Omega), { \ }\| w \|_{L^{3,\infty}} < \delta_2, { \ }\text{ and }{ \ }\| v_0 \|_{L^{3,\infty}} < \delta_2 ,
\end{equation*}
then there exists a unique global-in-time mild solution $v$ in $X_3$ of the system \eqref{eq15} with the initial datum $v_0$, satisfying $v \in X_p$, \eqref{eq16}, \eqref{eq18},
\begin{align}
& v \in BC ([0,\infty); L^{3,\infty}_{0,\sigma} (\Omega)), \notag \\
& \text{ for each fixed } T>0, { \ }\sup_{0 < t <T}\{t^\frac{1}{2} \| \nabla v (t ) \|_{L^{3 , \infty} } \} < + \infty \label{eq1233}, \\
& \sup_{0 < t \leq 1} \{ t^\frac{1}{2}  \| v ( t ) \|_{L^\infty} \} + \sup_{t \geq 1} \{ t^\frac{3 - \beta}{2 \beta} \| v ( t ) \|_{L^\infty} \} \leq \text{ Const. } <+ \infty .\label{eq111}
\end{align} 
Moreover, the following two assertions hold:\\
\noindent $(\mathrm{i})$ Assume in addition that $v_0 \in L^{\alpha , \infty }_\sigma ( \Omega )$ for some $3p/(p + 3) < \alpha < 3$ and that 
\begin{equation*}
\| v_0 \|_{L^{3 , \infty}} < \delta_1 \text{ and } \| w \|_{L^{3 , \infty } } < \delta_1 .
\end{equation*}
Then
\begin{equation}\label{eq00}
\sup_{0 < t \leq 1} \{ t^\frac{3}{2 \alpha }  \| v ( t ) \|_{L^\infty} \} + \sup_{t \geq 1} \{ t^{ \frac{3}{2} \left( \frac{1}{\alpha} + \frac{1}{ \beta } \right) - 1} \| v ( t ) \|_{L^\infty} \} \leq \text{ Const. } <+ \infty .
\end{equation}
Here $\delta_1$ is the constant appearing in Theorem \ref{thm110}.\\
\noindent $( \mathrm{ii})$ Assume in addition that $w$ satisfies Assumption \ref{assC}. Then for each fixed $T >0$
\begin{equation*}
\nabla v (T) \in L^{3 ,\infty}_0( \Omega ) .
\end{equation*}
Here $L^{3 ,\infty}_0 ( \Omega ) : = \overline{ C_0^\infty ( \Omega ) }^{\| \cdot \|_{L^{3 , \infty}} }$.
\end{theorem}

\begin{theorem}[Unique strong solution and stability]\label{thm112}{ \ }\\
Assume that $w$ is as in Assumption \ref{assC} and that $|u_\infty| \leq \gamma$ for some $\gamma > 0$. {\color{red}{Let $6 < p < \infty$ such that ${3 \beta}/(2 \beta -3)< p$}}. Then then there is $\delta_3 = \delta_3 ( \gamma , p ) > 0$ such that if
\begin{equation*}
v_0 \in L_{0 , \sigma}^{3,\infty} (\Omega), { \ }\| w \|_{L^{3,\infty}} < \delta_3 , { \ }\text{ and }{ \ }\| v_0 \|_{L^{3,\infty}} < \delta_3 ,
\end{equation*}
then there exists a unique global-in-time strong solution $v$ in $X_3$ of the system \eqref{eq15} with the initial datum $v_0$, satisfying $v \in X_p$, \eqref{eq16}, \eqref{eq18}, \eqref{eq1233}, \eqref{eq111},
\begin{equation*}
v \in BC ( [0,\infty); L_{0 , \sigma}^{3,\infty} (\Omega)), { \ } L v, v_t \in C ( (0,\infty); L_{0 , \sigma}^{3,\infty} (\Omega)) .
\end{equation*}
\end{theorem}

\begin{remark}
$(\mathrm{i})$ We can choose $u_\infty =0$ in Theorems \ref{thm19}-\ref{thm112}. \\
$(\mathrm{ii})$ It is not clear whether $\| v (t) \|_{L^\infty} = O (t^{-\frac{1}{2}})$ as $t \to \infty$ under our assumptions.\\
$(\mathrm{iii})$ Shibata \cite{Shi01} gave the sketch of the proof of Theorem \ref{thm19}. This paper describes a detailed proof of Theorem \ref{thm19}.\\
$(\mathrm{iv})$ Under only the condition that $v_0 \in L^{3,\infty}_\sigma (\Omega )$ it is difficult to drive \eqref{eq18}. See Yamazaki \cite[Remark 1.4]{Yam00} for the reason.\\
$(\mathrm{v})$ In order to show that $\| v (t ) - v_0 \|_{L^{3,\infty}} \to 0$ as $t \to 0 + 0$, we make use of $H^1$-solutions of the system \eqref{eq15}. We need the assumption that $\| \nabla w \|_{L^2}$ is finite to construct a $H^1$-solution of \eqref{eq15}. See Koba \cite{Kob13} for details.\\
{\color{red}{$(\mathrm{vi})$ We use the restriction that ${3 \beta }/(2 \beta - 3 ) < p < \infty$ to construct a global-in-time mild solution of the system \eqref{eq15}, while we need the condition that $6 < p < \infty$ to derive $L^\infty$-decay for the mild solution. See Section \ref{sect4} for details.}}
\end{remark}

Let us state results related to the main results of this paper. We first introduce stability results for the system \eqref{eq15} when $u_\infty =0$ and $w \neq 0$. Kozono and Ogawa \cite{KOK94} studied \eqref{eq15} under the conditions that $v_0 \in L^3_\sigma (\Omega)$ and $w \in W^{1,\infty} (\Omega)$. They proved that \eqref{eq15} admits a unique global-in-time strong $L^3$-solution satisfying $L^3$-asymptotic stability if $\| v_0 \|_{L^3}$, $\| w \|_{L^3}$, and $\| \nabla w \|_{L^{3/2}}$ are sufficiently small. Borchers and Miyakawa \cite{BM95} established the existence of a unique global-in-time strong $L^{3,\infty}$-solution of \eqref{eq15} when $\| v_0 \|_{L^{3,\infty}}$ and $\sup \{ |x| |w(x)| \} + \sup \{ |x|^2 | \nabla w (x)| \}$ are sufficiently small. They derived $L^{3,\infty}$-asymptotic stability of their solution in the case when $v_0 \in L^{3,\infty}_{0,\sigma} (\Omega)$, and investigated $L^\infty$-decay of their solution when $v_0 \in L^{3,\infty}_\sigma ( \Omega ) \cap L^{p,\infty}_{\sigma} ( \Omega )$ for some $1 < p < 3$. Under the conditions that $w \in L^{3,\infty} \cap L^\infty$ and $\nabla w \in L^p$ for some $p>3$, Kozono and Yamazaki \cite{KY98B} proved that \eqref{eq15} admits a unique global-in-time strong $L^{3,\infty}$-solution satisfying $L^q(q >3)$-asymptotic stability if $\| v_0 \|_{L^{3,\infty}} $ and $\| w \|_{L^{3,\infty}}$ are sufficiently small. Therefore the main results of this paper improve their results. Note that they studied $n$-dimensional case. A key observation of their methods is to derive $L^{p,r}$-$L^{q,r}$ estimates for the semigroup generated by the main linear operator $- \mathcal{L}$ of \eqref{eq15} under their restrictions, where $\mathcal{L} f = P\{ -\Delta f +  ( f , \nabla )w + (w,  \nabla ) f  \}$. However, in our situation, it is not easy to derive $L^{p,r}$-$L^{q,r}$ estimates for the semigroup whose generator is the main linear operator $-\mathscr{L}$ of \eqref{eq15}. Here $\mathscr{L}f = P\{ -\Delta f + (u_\infty , \nabla ) f +  ( f , \nabla )w + (w,  \nabla ) f  \}$.

Next we introduce stability results for the system \eqref{eq15} when $u_\infty \neq 0$ and $w \neq 0$. Galdi-Heywood-Shibata \cite{GHS97} and Shibata \cite{Shi99} applied $L^p$-$L^q$ estimates for the Oseen semigroup, which was obtained by Kobayashi-Shibata \cite{KS98}, to prove that \eqref{eq15} has a unique global-in-time mild $L^3$-solution satisfying $L^q(q>3)$-asymptotic stability when $\| v_0 \|_{L^3}$ and $\sup \{ (1 + |x|)(1 + |x| - x \cdot u_\infty/|u_\infty|)^{\epsilon} | w (x) | \} + \sup \{ (1 + |x| )^{3/2}(1 + | x | - x \cdot u_\infty/|u_\infty|)^{1/2 + \epsilon} | \nabla w (x)|  \}$ are sufficiently small for some $\epsilon > 0$. Shibata \cite{Shi01} showed the existence of a unique global-in-time generalized mild $L^{ 3 , \infty}$-solution, satisfying $L^q (q>3)$-asymptotic stability, of \eqref{eq15} when $\| v_0 \|_{ L^{3,\infty}}$ and $\| w \|_{L^{3,\infty}}$ are sufficiently small by applying $L^p$-$L^q$ estimates for the Oseen semigroup and the real interpolation theory. Enomoto and Shibata \cite{ES05} studied \eqref{eq15} when $w \in L^{3/(1 + \epsilon )} \cap L^{3/(1 - \epsilon)} $ and $\nabla w \in L^{3/(2 + \epsilon)} \cap L^{3/(2-\epsilon)} $ for some small $\epsilon >0$. They showed that there exists a unique global-in-time mild $L^3$-solution of \eqref{eq15} satisfying both $L^3$-asymptotic stability and $L^\infty$-asymptotic stability if $\| v_0 \|_{L^3}$, $\| w \|_{L^{3/(1- \epsilon)}}$, $\| w \|_{L^{3/(1+ \epsilon)}}$, $\| \nabla w \|_{L^{3/(2+ \epsilon)}}$, and $\| \nabla w \|_{L^{3/(2 - \epsilon)}}$ are sufficiently small. Under the same assumptions in \cite{ES05}, Bae and Roh \cite{BR12} investigated the decay rate with respect to time of $L^3$-solutions of \eqref{eq15} when the initial datum is in a weighted Lebesgue space. Therefore the main results of this paper are the generalization of a part of their results. A key tool of their methods is $L^p$-$L^q$ estimates for the Oseen semigroup. However, since our assumptions on $w$ is weaker than those of \cite{GHS97}, \cite{Shi99}, \cite{ES05}, we cannot directly apply their method to show the stability \eqref{eq16}-\eqref{eq00}.

In order to overcome the difficulties mentioned above, we use an equivalent norm on the Lorentz norm:
\begin{align*}
\| f \|_{L^{q,r}} \leq C \sup_{\phi \in [ L^{q' , r'}( \Omega )]^3 , { \ }\| \phi \|_{L^{q',r'}} \leq 1 } |\dual{f , \phi }| { \ \ \ }\text{ for } f \in [ L^{q,r}(\Omega) ]^3,
\end{align*}
where $1 < q < \infty$, $1 \leq r \leq \infty$, $1/q+1/{q'}=1$, $1/r + 1/{r'}=1$, apply $L^p$-$L^q$ type estimates for the Oseen semigroup: 
\begin{equation*}
\| \mathrm{e}^{- t L} P \partial_j f \|_{L^{p,1}} \leq C t^{-\frac{1}{2} - \frac{3}{2} (\frac{1}{q} - \frac{1}{p})} \| f \|_{L^{q,\infty}},
\end{equation*}
and make use of the following properties of the Oseen semigroup:
\begin{align*}
& \| ( \mathrm{e}^{- t_2 L} - 1) \mathrm{e}^{- t_1 L } P \phi \|_{L^{q, 1}} \leq C t_2^\alpha \mathrm{e}^{\frac{t_1}{2} + t_2} t_1^{ -\alpha -\frac{3}{2} \left( 1 - \frac{1}{q} \right)} \| \phi \|_{L^1},\\
& \| ( \mathrm{e}^{- t_2 L} - 1) \mathrm{e}^{- t_1 L } P \partial_j \phi \|_{L^{q,1}} \leq C t_2^\alpha \mathrm{e}^{\frac{t_1}{2} + t_2} t_1^{- \alpha - \frac{3}{2} \left( 1 - \frac{1}{q} \right)} \| \phi \|_{L^{\frac{3}{2},1} } . 
\end{align*}
\noindent More precisely, we apply an equivalent norm on the Lorentz space and $L^p$-$L^q$ estimates for the Oseen semigroup to show the existence of a unique global-in-time generalized mild solution of the system \eqref{eq15} satisfying $L^{3,\infty}$-asymptotic stability. Using the $L^p$-$L^q$ type estimates for the Oseen semigroup, we establish the existence of a unique global-in-time mild solution of the system \eqref{eq15} satisfying $L^\infty$-asymptotic stability. We make use of both the above estimates for the Oseen semigroup and fractional powers of the Oseen operator to show that the mild solution is a strong solution of the system \eqref{eq15}.

Let us briefly explain our idea of deriving $L^\infty$-decay for mild solutions of the system \eqref{eq15}. We first apply a duality argument to have
\begin{multline*}
\| v (  t ) \|_{L^\infty} \leq \| \mathrm{e}^{- t L} a \|_{L^\infty} \\
+ C \sup_{\phi \in [ C_0^\infty ]^3, { \ }\| \phi \|_{L^1} \leq 1} \bigg| \int_0^t \dual{  \mathrm{e}^{- \frac{(t-s)}{2} L} P (v , \nabla ) v  , \mathrm{e}^{- \frac{(t-s)}{2} L^*} P \phi  } d s \bigg|\\
+ C\sup_{\phi \in [ C_0^\infty ]^3, { \ }\| \phi \|_{L^1} \leq 1} \bigg| \int_0^t \dual{  \mathrm{e}^{- \frac{(t-s)}{2} L} P \{ (v , \nabla ) w + ( w , \nabla ) v \} , \mathrm{e}^{- \frac{(t-s)}{2} L^*} P \phi  } d s \bigg| .
\end{multline*}
Next we use the Cauchy-Schwarz inequality and $L^p$-$L^q$ type estimates for the Oseen semigroup to derive $L^\infty$-decay for our solutions. This is one of the key ideas of this paper to show the stability.

Finally, we state some results on $L^{3,\infty}$-solutions of incompressible fluid systems. Yamazaki \cite{Yam00} showed the existence of a unique generalized mild $L^{n,\infty}$-solution of the Navier-Stokes system with time-dependent external force in unbounded domains in $\mathbb{R}^n$. He derived the decay rate with respect to time of his solution by studying the time-dependent external force. Hishida-Shibata \cite{HS09} considered the stability of the Navier-Stokes flow in an exterior domain to a rotating obstacle. They derived $L^p$-$L^q$ estimates for the semigroup whose generator is the Stokes operator with both the Coriolis and centrifugal forces to show the existence of a unique global-in-time generalized mild $L^{3, \infty}$-solution, satisfying $L^q(q>3)$-asymptotic stability, of their system. Kang-Miura-Tsai \cite{KMT12} showed the existence of a unique very weak $L^{3,\infty}$-solution of the Navier-Stokes system with non-decaying boundary data. They studied the stability of the solution when the initial datum is asymptotically self-similar and the external force is time periodic.

\subsection{Notation} 

Let us introduce some fundamental notation.

For $m \in \mathbb{N}$, and $1 \leq q, r \leq \infty$, the symbols $L^q(\Omega)$, $W^{m,q}(\Omega)$, $L^{q,r}(\Omega)$ denote the usual Lebesgue space, the Sobolev space, and the Lorentz space with norms $\| \cdot \|_{L^q}(=\| \cdot \|_{L^q(\Omega)})$, $\| \cdot \|_{W^{m,q}} (=\| \cdot \|_{W^{m,q}(\Omega)})$, and $\| \cdot \|_{L^{q,r}}(= \| \cdot \|_{L^{q,r}(\Omega )})$, respectively. See Subsection \ref{subsec21} for the Lorentz spaces.

By $P$, we denote the Helmholtz projection, that is, $P: [ L^{q,r}(\Omega)]^3 \to L^{q,r}_\sigma (\Omega)$. See Subsection \ref{subsec22}.

Let $X$ be a Banach space, and $\mathscr{A}$ a linear operator densely defined in $X$. The symbols $X^*$ and $\mathscr{A}^*$ represent the dual space of $X$ and the dual operator of $\mathscr{A}$, respectively.

Let $X$ be a Banach space, and $\mathscr{A}$ a linear operator on $X$. The symbols $D(\mathscr{A})$, $R(\mathscr{A})$, and $N ( \mathscr{A})$ represent the domain of $\mathscr{A}$, the range of $\mathscr{A}$, and the null space of $\mathscr{A}$, respectively. When $\mathscr{A}$ generates a semigroup on $X$, we write the semigroup as $\mathrm{e}^{ t \mathscr{A}}$. We also write the dual semigroup of $\mathrm{e}^{ t \mathscr{A}}$ as $\mathrm{e}^{ t \mathscr{A}^*}$.

We will use the symbol $C$ to denote a positive constant. We write $C( \eta_1, \eta_2)$ if the constant $C$ depends on certain quantities $\eta_1,\eta_2$. However, the dependency on $\Omega$ is usually omitted.

\subsection{Outline of this Paper} 
In Section \ref{sect2}, we study the Lorentz spaces $L^{q,r}$, the solenoidal spaces $L^{q,r}_\sigma$ and $L^{q,\infty}_{0,\sigma}$, the Helmholtz projection $P$, the Oseen operator $L$, and the Oseen semigroup $\mathrm{e}^{-t L}$. We first give fundamental properties of the Lorentz spaces, the solenoidal Lorentz spaces, and the Helmholtz projection. Next we study the dual semigroup of the Oseen semigroup, and derive $L^p$-$L^q$ type estimates for the Oseen semigroup by applying $L^p$-$L^q$ estimates for the Oseen semigroup, a duality argument, and the real interpolation theory. Finally we study fractional powers of the Oseen operator in $L^{q,r}_{\sigma}$ and derive key estimates for the Oseen semigroup .

In Section \ref{sect3}, we prove Theorems \ref{thm19} and \ref{thm110}. We first derive basic properties of the global-in-time generalized mild solutions of the system \eqref{eq15}. Applying an equivalent norm on the Lorentz space, $L^p$-$L^q$ estimates for the Oseen semigroup, and a contraction mapping theory to establish the existence of a unique global-in-time generalized mild solution of the system \eqref{eq15}. Moreover, we investigate $L^{3,r}$-asymptotic stability of the solution when the initial datum belongs to $L^{3,r}_\sigma ( \Omega )$.

In Section \ref{sect4}, we give the proof of Theorem \ref{thm111}. We first construct a unique local-in-time mild $L^{3,\infty}$-solution of the system \eqref{eq15} by using both an equivalent norm on the Lorentz norm and $L^p$-$L^q$ type estimates for the Oseen semigroup. Next we apply the uniqueness of the generalized mild solutions of \eqref{eq15} and the existence of a unique global-in-time generalized mild solution of \eqref{eq15} to construct a unique global-in-time mild solution of \eqref{eq15}. Finally, we investigate $L^\infty$-asymptotic stability of the global-in-time mild solution.

In Section \ref{sect5}, we show the existence of a unique global-in-time strong solution of the system \eqref{eq15} to prove Theorem \ref{thm112}. Applying fractional powers of the Oseen operator and some estimates for the Oseen semigroup, we prove that the mild solutions obtained by Section \ref{sect4} are strong solutions of \eqref{eq15}. 

In the Appendix, we characterize the Lorentz norm by using the measure theory and the definition of the Lorentz norm.

\section{Preliminaries}\label{sect2}

In this section we prepare key tools to prove Theorems \ref{thm19}-\ref{thm112}. We first recall fundamental properties of the Lorentz spaces. We state useful inequalities to analyze the system \eqref{eq15}. Secondly, we study solenoidal $L^{q,r}$-function spaces and the Helmholtz projection on the Lorentz spaces $(1<q<\infty, 1\leq r \leq \infty)$. Thirdly, we consider the Stokes operator and the Oseen operator in a solenoidal $L^{q , r}$-space. More precisely, we observe that the Oseen operator generates an analytic semigroup on $L^{q,r}_{\sigma } ( \Omega )$, and characterize the dual operator of the Oseen operator in $L^{q,r}_{\sigma }( \Omega ) $ $(r \neq \infty)$. Fourthly, we state well-known $L^p$-$L^q$ estimates for the Oseen semigroup. Using the $L^p$-$L^q$ estimates, the real interpolation theory, and a duality argument, we derive $L^{p}$-$L^q$ type estimates for the Oseen semigroup. Fifthly, we study fractional powers of the Oseen operator in $L^{q,r}_\sigma (\Omega)$. Especially, we deal with fractional powers of the Oseen operator in $L^{q,\infty}_\sigma (\Omega )$ and $L^{q,1}_\sigma (\Omega )$. Finally, we derive useful estimates for the Oseen semigroup to show the existence of a strong solution of \eqref{eq15}.

\subsection{Lorentz Spaces}\label{subsec21} 
Let us recall the Lorentz spaces. For $1 \leq q \leq \infty$ and $1 \leq r \leq \infty$
\begin{equation*}
L^{q,r} ( \Omega ) := \{ f \in L^1 (\Omega) + L^\infty ( \Omega ) ; { \ }\| f \|_{L^{q,r} (\Omega )} < +\infty  \}
\end{equation*}
with
\begin{equation*}
\| f \|_{L^{q,r}} = \| f \|_{L^{q,r}(\Omega)} := 
\begin{cases}
\left( \int_0^\infty (t^{1/q} f^{**} (t) )^r \frac{d t}{t} \right)^{1/r} & \text{ if } 1 \leq r < + \infty, \\
\sup_{t >0} \{ t^{1/q} f^{**} ( t ) \}  & \text{ if } r = \infty .
\end{cases}
\end{equation*}
Here 
\begin{align*}
&f^{**}(t) : = t^{-1}\int_0^t f^*(s) d s, { \ }t \geq 0,\\
&f^* (t) : = \inf \{ \sigma >0  ; { \ }  \mu \{ x \in \Omega; { \ } | f (x) | > \sigma \} \leq t \}, { \ } t \geq 0,
\end{align*}
where $\mu \{ \cdot \}$ denotes the $3$-dimensional Lebesgue measure. For all \\$f = (f^1, f^2,f^3 ) \in [L^{q,r}(\Omega)]^3$,
\begin{equation*}
\| f \|_{L^{q,r} (\Omega)} := \|{ \ } |f| { \ }\|_{L^{q,r}(\Omega)} = \| \sqrt{  ( f^1)^2 + (f^2)^2 + (f^3)^2} \|_{L^{q,r}( \Omega ) }.
\end{equation*}
Note that there are $C_1 = C_1 (q,r) >0$ and $C_2 = C_2 (q,r) >0$ such that\\
$\| f \|_{L^{q,r} (\Omega)} \leq C_1 ( \| f^1 \|_{L^{q,r} (\Omega)} + \| f^2 \|_{L^{q,r} (\Omega )} + \| f^3 \|_{L^{q,r} (\Omega )}) \leq C_2 \|  f \|_{L^{q,r} ( \Omega )}$. See the Appendix for details. From \cite[Section 3.3]{BB67}, \cite[Chapter 1]{BL76}, and \cite[IV Lemma 4.5]{BS88}, we see that $L^{q,r} (\Omega)$ is a Banach space when $1 < q < \infty$ and $1 \leq r \leq \infty$. We also see that $L^{q,q}(\Omega) = L^q (\Omega)$ when $1 < q \leq \infty$ and that $L^{1,\infty}(\Omega) = L^1 (\Omega)$. Furthermore, we find that $C_0^\infty (\Omega) $ is dense in $L^{q,r}(\Omega )$ if $1 < q < \infty$ and $1 \leq r < \infty$. From \cite[$\S$1.4.2 and Corollary 1.7]{Lun09}, we obtain
\begin{lemma}\label{lem21}
Let $1 \leq q_1 < q < q_2 \leq \infty$ and $1 \leq r \leq \infty$. Then there is $C = C ( q_1 , q , q_2 , r ) > 0 $ such that for all $f \in L^{q_1,\infty} ( \Omega )\cap L^{q_2,\infty}(\Omega)$
\begin{equation}
\| f \|_{L^{q,r}} \leq C \| f \|_{L^{q_1, \infty}}^{1 - \theta } \| f \|_{L^{q_2, \infty}}^{\theta} \label{eq21} .
\end{equation}
Here $\theta = \{ q_2(q-q_1)\}/\{q(q_2 - q_1)\}$.
\end{lemma}
\noindent From \cite[Chapters 2 and 4]{BS88}, we have the two lemmas.
\begin{lemma}\label{lem22}Let $1 < q < \infty$ and $1 \leq r_1 \leq r_2 \leq \infty$. Then there is $C = C(q,r_1,r_2) >0$ such that for all $f \in L^{q,r_1}(\Omega)$
\begin{equation}\label{eq22}
\| f \|_{L^{q,r_2}} \leq C \| f \|_{L^{q,r_1}} .
\end{equation}
\end{lemma}
\begin{lemma}\label{lem23}
Let $1 \leq q,q' \leq \infty$ and $1 \leq r,r' \leq \infty$ such that $1/q+1/{q'}=1$ and $1/r + 1/{r'}=1$. Then for all $f \in [ L^{q,r}(\Omega ) ]^3$ and $g \in [ L^{q',r'}(\Omega) ]^3$,
\begin{equation}\label{eq23}
|\dual{ f , g } | \leq \| f \|_{L^{q,r}} \| g \|_{L^{q',r'}} .
\end{equation}
Here $L^{1, \ell} (\Omega) := L^1 (\Omega )$ and $L^{\infty , \ell} (\Omega ) := L^\infty ( \Omega )$ for $1 \leq \ell \leq \infty$.
\end{lemma}
Next we state the weak H\"{o}lder inequality.
\begin{lemma}\label{lem24} $(\mathrm{i})$ {\color{red}{\cite[Proposition 2.1]{KY99}}}Let $1 < q, q_1, q_2 < \infty$ and $1 \leq r_1 , r_2 \leq \infty$ such that $1/q = 1/{q_1} + 1/{q_2}$. Then there is $C = C(q , q_1 , q_2,r_1,r_2)>0$ such that for all $f \in L^{q_1 , r_1} (\Omega)$ and $g \in L^{q_2,r_2}(\Omega)$
\begin{equation}\label{eq24}
\| f g \|_{L^{q,r} } \leq C \| f \|_{L^{q_1 , r_1}} \| g \|_{L^{q_2 , r_2}}, \text{ where }r := \min \{ r_1 , r_2 \}.
\end{equation}
\noindent $(\mathrm{ii})$ Let $1 < q < \infty$ and $1 \leq r \leq \infty$. Then for all $f \in L^{q, r} ( \Omega )$ and $g \in L^\infty (\Omega)$
\begin{equation}\label{eq25}
\| f g \|_{L^{q,r}} \leq \| f \|_{L^{q,r}} \| g \|_{L^{ \infty } }.
\end{equation}
\end{lemma}
\noindent Combining the usual H\"{o}lder inequality and the Marcinkiewicz interpolation theorem gives the assertion $(\mathrm{ii})$ of Lemma \ref{lem24}. {\color{red}{See also the proof of Lemma 8.2 in \cite{HS09}.}}

Finally, we describe a characterization of the Lorentz norm $\| \cdot \|_{L^{q,r}}$. 
\begin{lemma}[Equivalent norms for the Lorentz space]\label{lem25}{ \ }\\
\noindent $(\mathrm{i})$ Let $1 < q,q' < \infty$ and $1 \leq r,r' \leq \infty$ such that $1/q + 1/{q'}=1$ and $1/r + 1/r' =1$. For each $f = (f^1,f^2,f^3) \in [L^{q,r}(\Omega)]^3$,
\begin{align*}
& \| f \|_{X^{q,r} (\Omega)} := \| f \|_{L^{q,r}} \equiv \| \sqrt{  ( f^1)^2 + (f^2)^2 + (f^3)^2} \|_{L^{q,r}},\\
& \| f \|_{Y^{q,r} (\Omega)} := \max_{1 \leq j \leq 3} \{ \| f^j \|_{L^{q,r}} \},\\
& \| f \|_{Z^{q,r} (\Omega)} := \sup_{\phi \in [ L^{q',r'}(\Omega ) ]^3,{ \ } \| \phi \|_{L^{q',r'}} \leq 1 } |\dual{f , \phi}|.
\end{align*}
Then there is $C_\vartriangle = C_\vartriangle (q,r) >0$ such that for all $f \in [L^{q,r}(\Omega)]^3$
\begin{equation*}
\| f \|_{X^{q,r}} \leq C_\vartriangle \| f \|_{Y^{q,r}} \leq C_\vartriangle \| f \|_{Z^{q,r}} \leq C_\vartriangle \| f \|_{X^{q,r}} .
\end{equation*}
Moreover, assume in addition that $r \neq 1$. Then
\begin{equation*}
\| f \|_{Z^{q,r} (\Omega)}= \sup_{\phi \in [ C_0^\infty (\Omega ) ]^3,{ \ } \| \phi \|_{L^{q',r'}} \leq 1 } |\dual{f , \phi}|.
\end{equation*}
\noindent $(\mathrm{ii})$ There is $C>0$ such that for each $f = (f^1,f^2,f^3) \in [L^\infty ( \Omega ) ]^3$
\begin{equation}\label{eq26}
\| f \|_{L^\infty} \leq C \sup_{\phi \in [ C_0^\infty ( \Omega )]^3, { \ }\| \phi \|_{L^1} \leq 1} |\dual{f , \phi }| .
\end{equation}
\end{lemma}
The Appendix gives the proof of the assertion $(\mathrm{i})$ of Lemma \ref{lem25}.

\subsection{Solenoidal Spaces and Helmholtz Projection}\label{subsec22}

We first state the dual of the solenoidal space $L^{q,r}_\sigma (\Omega)$. From \cite[Section 3.7]{BL76} and \cite[Theorems 5.2 and 5.5]{BM95}, we obtain
\begin{lemma}\label{lem26}
Let $1 < q,q' < \infty$ and $1 \leq r,r' \leq \infty$ such that $1/q + 1/{q'}=1$ and $1/r + 1/{r'}=1$. Then,\\
\noindent $(\mathrm{i})$ If $1 \leq r < \infty$, then $C_{0,\sigma}^\infty (\Omega)$ is dense in $L^{q,r}_\sigma ( \Omega )$.\\
\noindent $(\mathrm{ii})$ If $1 \leq r < \infty$, then the dual of $L^{q,r}_\sigma (\Omega)$ is $L^{q',r'}_\sigma (\Omega)$.\\
\noindent $(\mathrm{iii})$ The dual of $L^{q, \infty}_{0,\sigma} (\Omega)$ is $L^{q', 1}_{\sigma}(\Omega)$.
\end{lemma}
\noindent See Definition \ref{def13} for the solenoidal spaces $C_{0 , \sigma}^\infty ( \Omega )$, $L^{q,r}_\sigma ( \Omega )$, and $L^{q,\infty}_{0,\sigma} ( \Omega )$.

Next we introduce the Helmholtz projection. Fix $1 < q < \infty$ and $1 \leq r \leq \infty$. From \cite[Theorem 5.2]{BM95}, we see that
\begin{equation*}
L^{q,r}_\sigma (\Omega) = \{ f \in L^{q,r}(\Omega )^3; { \ }\nabla \cdot f =0,{ \ } f \cdot \nu |_{\partial \Omega} =0 \}.
\end{equation*}
We also see that $L^{q,r}(\Omega)^3 = L^{q,r}_\sigma (\Omega) \oplus G_{q,r} (\Omega)$, where
\begin{equation*}
G_{q,r} (\Omega) := \{ \nabla p \in L^{q,r}(\Omega )^3 ;{ \ } p \in L^{q,r}_{ loc } (\overline{\Omega}) \} .
\end{equation*}

\begin{lemma}\cite[Theorem 5.2]{BM95} \label{lem27}
There exists a linear operator $P$ with the three properties:\\
\noindent $(\mathrm{i})$ For each $1 < q < \infty$ and $1 \leq r \leq \infty$ there is $C = C(q,r) >0$ such that for all $f \in L^{q,r}(\Omega)^3$
\begin{align}
& \| P f \|_{L^{q,r}} \leq C \| f \|_{L^{q,r}},\label{eq27}\\
& Pf \in L^{q,r}_\sigma ( \Omega ).\notag
\end{align}
\noindent $(\mathrm{ii})$ Let $1 < q < \infty$ and $1 \leq r \leq \infty$. Then for all $g \in L^{q,r}_\sigma ( \Omega )$
\begin{equation*}
P g = g .
\end{equation*}
\noindent $(\mathrm{iii})$ Let $1 < q,q' < \infty$ and $1 \leq r,r' \leq \infty$ such that $1/q + 1/{q'}=1$ and $1/r + 1/{r'}=1$. Then for every $f \in [ L^{q,r}(\Omega) ]^3$ and $g \in [ L^{q',r'}(\Omega ) ]^3$
\begin{equation*}
\dual{P f, g} = \dual{f , P g} .
\end{equation*}
\end{lemma}
\noindent See also \cite{FM77}, \cite{SS92}, and \cite[Lemma 1.1]{MY92} for the Helmholtz projection. Now we study the function spaces $L^{q,\infty}_\sigma (\Omega)$ and $L^{q,\infty}_{0 , \sigma} ( \Omega )$.
\begin{lemma}\label{lem28}Let $1 < q <  \infty$. Set
\begin{equation*}
L^{q,\infty}_0 (\Omega) = \overline{C_0^\infty (\Omega) }^{\| \cdot \|_{L^{q,\infty}}} .
\end{equation*}
Then
\begin{equation*}
P : [ L^{q,\infty}_0 (\Omega ) ]^3 \to L^{q,\infty}_{0,\sigma } ( \Omega ) .
\end{equation*}
\end{lemma}
\begin{proof}[Proof of Lemma \ref{lem28}]
Fix $q \in (1,\infty)$ and $\varepsilon >0$. Let $f \in [ L^{q ,\infty}_0 (\Omega ) ]^3$. By the definition of the function space $L^{q,\infty }_{0} ( \Omega )$, there is $f_0 \in [C_0^\infty (\Omega)]^3$ such that 
\begin{equation}\label{eq28}
\| f - f_0 \|_{L^{q,\infty}} < \varepsilon .
\end{equation}
Since $P f_0 \in L^q_\sigma (\Omega)$, we take $f_1 \in C_{0,\sigma}^\infty (\Omega)$ such that
\begin{equation}\label{eq29}
\| P f_0  - f_1 \|_{L^q} < \varepsilon .
\end{equation}
Using \eqref{eq27}, \eqref{eq22}, \eqref{eq28}, and \eqref{eq29}, we check that
\begin{align*}
\| Pf - f_1 \|_{L^{q,\infty} } \leq & \| Pf - P f_0 \|_{L^{q,\infty}} + \| P f_0 - f_1 \|_{L^{q,\infty}}\\
\leq & C \| f - f_0 \|_{L^{q,\infty}} + C \| P f_0 - f_1  \|_{L^q} < C \varepsilon .
\end{align*}
Since $f_1 \in C_{0,\sigma}^\infty (\Omega)$ and $\varepsilon$ is arbitrary, therefore the lemma follows.
\end{proof}
Applying Lemmas \ref{lem25}-\ref{lem27}, we have
\begin{lemma}\label{lem29}
Let $1 < q,q' < \infty$ and $1 \leq r, r' \leq \infty$ such that $1/q + 1/{q'}=1$ and $1/r + 1/{r'}=1$. Then there is $C=C (q,r) >0$ such that for every $f \in L^{q,r}_\sigma (\Omega)$
\begin{equation}\label{eq210}
 \| f \|_{L^{q,r}} \leq C \sup_{\varphi \in L^{q' ,r'}_\sigma (\Omega ) , { \ }\| \varphi \|_{L^{q',r'} }\leq 1} \left| \dual{f , \varphi} \right| .
\end{equation}
In particular, if $r \neq 1$, then
\begin{equation}\label{eq211}
 \| f \|_{L^{q,r}} \leq C \sup_{\varphi \in C_{0,\sigma}^\infty (\Omega ) , { \ }\| \varphi \|_{L^{q',r'} }\leq 1} \left| \dual{f , \varphi} \right| .
\end{equation}
\end{lemma}
\begin{proof}[Proof of Lemma \ref{lem29}]
Fix $f \in L^{q,r}_\sigma (\Omega)$. From Lemmas \ref{lem25}-\ref{lem27}, we observe that
\begin{align*}
\| f \|_{L^{q,r}} \leq & C \sup_{\phi \in [L^{q', r'}(\Omega ) ]^3, { \ }\| \phi \|_{L^{q',r'}} \leq 1 } |\dual{f , \phi} |\\
= &C \sup_{\phi = \phi_1 + \phi_2, { \ }\phi_1  \in L^{q', r'}_\sigma, { \ }\phi_2 \in G_{q',r'} { \ }\| \phi \|_{L^{q',r'}} \leq 1 } |\dual{f , P \phi} |\\
= &C \sup_{\phi_1 \in L^{q', r'}_\sigma, { \ }\| \phi_1 \|_{L^{q',r'}} \leq 1 } |\dual{f , \phi_1} | .
\end{align*}
Therefore we see that
\begin{equation*}
\| f \|_{L^{q,r}} \leq C \sup_{\varphi \in L^{q', r'}_\sigma ( \Omega ) , { \ }\| \varphi \|_{L^{q',r'}} \leq 1 } |\dual{f , \varphi} | .
\end{equation*}
Since $C_{0,\sigma}^\infty (\Omega ) $ is dense in $L^{q' , r'}_\sigma ( \Omega )$ if $r' \neq \infty$, we deduce \eqref{eq211}.
\end{proof}

\subsection{Oseen Semigroup and Its Dual Semigroup}\label{subsec23} 
Let us now study the Stokes operator and the Oseen operator in a solenoidal $L^{q,r}$-space. Fix $1 < q < \infty$, $1 \leq r \leq \infty$, and $u_\infty \in \mathbb{R}^3$.
We define the three linear operators $A$, $L$, and $L^*$ in $L^{q,r}_\sigma ( \Omega )$ by
\begin{align*}
&\begin{cases}
A = A_{q , r} = P (- \Delta ) f,\\
D (A) = D (A_{q,r}) = \{ f \in L^{q,r}_\sigma (\Omega ) ; { \ } f|_{\partial \Omega} =0,{ \ } \| \nabla f \|_{L^{q,r}} + \| \nabla^2 f \|_{L^{q,r}} < + \infty \},
\end{cases}\\
&\begin{cases}
L= L_{q,r}  = P\{ (- \Delta) f + (u_\infty, \nabla ) f \},\\
D (L) = D (L_{q,r}) = D (A_{q,r}) .
\end{cases}\\
&\begin{cases}
L^*= L^*_{q,r}  = P\{ (- \Delta) f - (u_\infty, \nabla ) f \},\\
D (L^*) = D (L^*_{q,r}) = D (A_{q,r}) .
\end{cases}
\end{align*}
We call $A_{q,r}$ the \emph{Stokes operator} in $L^{q,r}_\sigma (\Omega)$ and $L_{q,r}$ the \emph{Oseen operator} in $L^{q,r}_\sigma(\Omega)$. From the argument in \cite[Sect.5]{BM95}, \cite{Shi01}, and the resolvent estimates for the Oseen operator in \cite{ES04,ES05}, we have
\begin{lemma}\label{lem210}
Let $1 < q < \infty$ and $1 \leq r \leq \infty$. Then
\begin{enumerate}
\renewcommand{\labelenumi}{\rm{(\roman{enumi})}}
\item There is $C = C (\Omega , q,r) >0$ such that for all $f \in D (A_{q,r})$
\begin{equation*}
\| \nabla f \|_{L^{q,r}} + \| \nabla^2 f \|_{L^{q,r}} \leq C (\| A_{q,r} f \|_{L^{q,r}} + \| f \|_{L^{q,r}}).
\end{equation*}
\item Each operator $- A_{q,r}$, $- L_{q,r}$ and $- L^*_{q,r}$ generates an analytic semigroup on $L^{q,r}_\sigma (\Omega )$.
\item If $r \neq \infty$, then each operator $-A_{q,r}$, $- L_{q,r}$, and $- L^*_{q,r}$ generates a $C_0$-semigroup on $L^{q,r}_\sigma (\Omega )$.
\end{enumerate}
\end{lemma}
\noindent See also \cite{Iwa89}, \cite{KS98} for the Oseen semigroup in exterior domains, and \cite{BS87}, \cite{GS89} for the Stokes semigroup in exterior domains.

We characterize the dual semigroup of the semigroup $\mathrm{e}^{- t L_{q,r}}$.
\begin{lemma}[Characterization of the dual semigroup of $\mathrm{e}^{- t L_{q,r}}$]\label{lem211}{ \ }\\
Let $1 < q < \infty$ and $1 \leq r < \infty$. Let $(L_{q,r})'$ be the dual operator of the operator $L_{q,r}$ in $L^{q,r}_{\sigma}(\Omega)$. Then $(L_{q,r})' = L^*_{q',r'}$ and $D ((L_{q,r})') = D (L^*_{q', r'})$. Here $1/q + 1/{q'}=1$ and $1/r + 1/{r'} = 1$, where $1/\infty := 0$. Moreover, for all $f \in L^{q,r}_{\sigma} (\Omega )$, $g \in L^{q',r'}_{\sigma}(\Omega)$, and $t>0$
\begin{equation}\label{c1}
\dual{ \mathrm{e}^{-t L_{q,r}} f , g } = \dual{f , \mathrm{e}^{-t L^*_{q'. r'}} g}.  
\end{equation}
\end{lemma}

\begin{proof}[Proof of Lemma \ref{lem211}]Let $1 < q \leq q' < \infty$ and $1 \leq r \leq r' \leq \infty$ such that $1/q + 1/{q'}=1$ and $1/r +1/{r'}=1$, where $r' := \infty$ if $r =1$. Fix $q,q',r,r'$. Since $\mathrm{e}^{-t L_{q,r}}$ is a $C_0$-semigroup on $L_{\sigma}^{q,r}$, it follows from the Hille-Yosida theorem or the Lumer-Phillips theorem that $D (L_{q,r})$ is dense in $L^{q,r}_{\sigma}(\Omega)$. Therefore one can define the dual operator of the operator $L_{q,r}$. It is easy to check that for all $f \in D(A_{q,r})$ and $g \in D (A_{q',r'})$
\begin{multline*}
\dual{L_{q,r} f , g} = \dual{P \{ - \Delta f + (u_\infty, \nabla ) \} f , g} = \dual{ f , P \{ - \Delta f - (u_\infty, \nabla )\} g}\\
 = \dual{f , L^*_{q',r'} g} = \dual{f , (L_{q,r})' g}.
\end{multline*}
This implies that $(L_{q,r})'= L^*_{q',r'}$ on $D (A_{q',r'})$. Fix $\phi_1 \in D( (L_{q,r})')$. By the definition of the dual operator, there is $h \in L^{q',r'}_{\sigma } ( \Omega )$ such that for all $f \in D (A_{q,r})$
\begin{equation*}
\dual{f , h} =  \dual{f , (L_{q,r})' \phi_1 } = \dual{L_{q,r} f , \phi_1}.
\end{equation*}
Now we show that $\phi_1 \in D(A_{q',r'})$. Since $\mathrm{e}^{-t L_{q,r}}$ and $\mathrm{e}^{-t L^*_{q',r'}}$ are analytic semigroups, we find $\eta >0$ such that $\eta \in \rho (- L_{q,r}) \cap \rho (- L^*_{q',r'})$. Since $\eta \in \rho (- L^*_{q',r'})$ and $h \in L^{q',r'}_{\sigma }( \Omega )$, there is $\phi_2 \in D (A_{q',r'})$ such that
\begin{equation*}
(L^*_{q',r'} + \eta ) \phi_2 = h .
\end{equation*}
Next we prove that $(L^*_{q',r'} + \eta )$ is injective. Let $\phi_3 \in D (A_{q',r'})$ such that for all $f \in D (A_{q,r})$
\begin{equation*}
0 = \dual{ f , (L^*_{q',r'} + \eta )\phi_3} = \dual{ (L_{q,r} + \eta )f ,  \phi_3 } .
\end{equation*}
Since $\eta \in \rho ( - L_{q,r})$, we find that $\phi_3 = 0$. Therefore we see that $(L^*_{q',r'} + \eta )$ is one-to-one. It is easy to check that for all $f \in D (A_{q,r})$
\begin{align*}
0 = & \dual{f , h} - \dual{f , h} \\
= & \dual{f , (L_{q,r})' \phi_1} - \dual{f , L^*_{q',r'} \phi_2 }\\
= & \dual{f , L^*_{q',r'} (\phi_1 - \phi_2 )}.
\end{align*}
Thus, we see that $\phi_1 = \phi_2$. This implies that $\phi_1 \in D (A_{q',r'})$, that is, $D ((L_{q,r})') \subset D (A_{q',r'})$. Therefore we conclude that $(L_{q,r})' = L^*_{q',r'}$ and $D ((L_{q,r})')= D (A_{q',r'})$. 

Finally, we deduce \eqref{c1}. We see at once that for $f \in L^{q,r}_{\sigma} (\Omega )$ and $g \in L^{q',r'}_{\sigma}(\Omega)$
\begin{align*}
\dual{ ( \lambda +  L_{q,r} )^{-1} f , g } = &\dual{ ( \lambda +  L_{q,r} )^{-1} f , (\lambda + L^*_{q',r'})(\lambda + L^*_{q',r'})^{-1} g } \\
=&  \dual{ f , (\lambda + L^*_{q',r'})^{-1} g }
\end{align*}
if $\lambda \in \rho ( - L_{q,r}) \cap \rho ( - L^*_{ q' , r' })$. Therefore we see \eqref{c1}. Remark that $\lambda \in \rho (- L_{q,r})$ then $\lambda \in \rho (-L_{q',r'}^*)$ and $(\lambda + L_{q',r'}^*)^{-1} = ((\lambda + L_{q,r})^{-1})^*$. See \cite[Lemma 10.2 in Chapter 1]{Paz83}.
\end{proof}
From Lemma \ref{lem211}, we observe that for all $f \in L^{q,1}_{\sigma} (\Omega )$, $g \in L^{q',\infty}_{\sigma}(\Omega)$, and $t >0$
\begin{equation*}
\dual{ \mathrm{e}^{-t L_{q,1}} f , g } = \dual{f , \mathrm{e}^{-t L^*_{q', \infty}} g}.  
\end{equation*}
Therefore we consider $\mathrm{e}^{- t L_{q,\infty}}$ as the dual semigroup of $\mathrm{e}^{- t L^*_{q',1}}$.

Next we study the Oseen operator in $L^{q, \infty}_{0,\sigma}(\Omega)$.
\begin{lemma}\label{lem212}
Let $1 < q < \infty$. Set
\begin{equation*}
\begin{cases}
\mathscr{L}= \mathscr{L}_{q,\infty}  = P\{ (- \Delta) f + (u_\infty, \nabla ) f \},\\
D (\mathscr{L}_{q,\infty}) = \{ f \in L^{q,\infty}_{0, \sigma } (\Omega) ; { \ }f|_{\partial \Omega} =0,{ \ \ } \| \nabla f \|_{L^{q,\infty}} + \| \nabla^2 f \|_{L^{q,\infty}} < + \infty \} .
\end{cases}
\end{equation*}
Then $- \mathscr{L}_{q,\infty}$ generates an analytic $C_0$-semigroup on $L^{q,\infty}_{0, \sigma} (\Omega )$. Moreover, if $1 < q \leq 3$, then $\nabla \mathrm{e}^{- t \mathscr{L}_{q , \infty}} f \in L^{q, \infty}_0 ( \Omega )$ for each $t>0$ and $f \in L^{q,\infty}_{0 , \sigma } ( \Omega )$.
\end{lemma}
\begin{proof}[Proof of Lemma \ref{lem212}]
Fix $1 < q < \infty$. We first show that $- \mathscr{L}_{q,\infty}$ generates an analytic semigroup on $L^{q,\infty}_{0, \sigma } (\Omega )$. Since $\mathrm{e}^{- t L_{q,\infty} } $ is an analytic semigroup on $L^{q,\infty}_\sigma ( \Omega)$, we see that there are $C = C(\Omega ,\gamma, q ) >0$ and $\eta = \eta (\Omega ,\gamma , q) >0 $ such that for all $f \in L^{q,\infty}_\sigma (\Omega)$ and $t>0$,
\begin{align}
& \mathrm{e}^{- t L_{q,\infty}} f \in L^{q,\infty}_{\sigma } (\Omega ), \\
&\| \mathrm{e}^{- t L_{q,\infty} } f \|_{L^{q,\infty}} \leq C \mathrm{e}^{\eta t } \| f \|_{L^{q,\infty}} , \label{c2}\\
&\| \mathscr{L}_{q,\infty} \mathrm{e}^{- t L_{q,\infty} } f \|_{L^{q,\infty}} \leq C t^{-1} \mathrm{e}^{\eta t } \| f \|_{L^{q,\infty}} \label{c222}.
\end{align}
Since $\mathscr{L}_{q,\infty} = L_{q , \infty}$ on $L^{q,\infty}_{0, \sigma} (\Omega)$, we only have to prove that for $t>0$
\begin{align*}
& \mathrm{e}^{- t \mathscr{L}_{q,\infty}} f \in L^{q,\infty}_{0, \sigma } (\Omega ) , \\
& L_{q,\infty} \mathrm{e}^{- t \mathscr{L}_{q,\infty}} f \in L^{q,\infty}_{0, \sigma } (\Omega ) 
\end{align*}
if $f \in L^{q,\infty}_{0,\sigma} (\Omega)$. Fix $t >0$ and $f \in L^{q,\infty}_{0,\sigma} (\Omega)$. By the definition, there are $f_m \in C_{0,\sigma}^\infty (\Omega)$ such that
\begin{equation*}
\lim_{m \to \infty} \| f -f_m \|_{L^{q,\infty}} = 0.
\end{equation*}
Since $C_{0, \sigma}^\infty (\Omega) \subset L^q_\sigma (\Omega) \subset L^{q,\infty}_{0,\sigma} (\Omega)$, we find that for each $m \in \mathbb{N}$
\begin{align*}
& \mathrm{e}^{- t \mathscr{L}_{q,\infty}} f_m \in L^{q,\infty}_{0, \sigma } (\Omega ) , \\
& L_{q,\infty} \mathrm{e}^{- t \mathscr{L}_{q,\infty}} f_m \in L^{q,\infty}_{0, \sigma } (\Omega ) .
\end{align*}
By \eqref{c2} and \eqref{c222}, we observe that as $m \to \infty$
\begin{align*}
& \| \mathrm{e}^{- t \mathscr{L}_{q,\infty}} f - \mathrm{e}^{- t L_q} f_m \|_{L^{q,\infty}} \leq C \mathrm{e}^{\eta t}\| f - f_m \|_{L^{q,\infty}} \to 0 ,\\
& \| \mathscr{L}_{q,\infty} \mathrm{e}^{- t \mathscr{L}_{q,\infty}} f - L_{q,\infty} \mathrm{e}^{- t L_q} f_m \|_{L^{q,\infty}} \leq C t^{-1} \mathrm{e}^{\eta t}\| f - f_m \|_{L^{q,\infty}} \to 0 . 
\end{align*}
Therefore we see that $- \mathscr{L}_{q,\infty}$ generates an analytic semigroup on $L^{q,\infty}_{0, \sigma } (\Omega )$.

Next we prove that $- L_{q,\infty}$ generates a $C_0$-semigroup on $L^{q,\infty}_{0, \sigma } (\Omega )$. Fix $f \in L^{q,\infty}_{0,\sigma} (\Omega)$ and $\varepsilon >0$. By definition, there is $f_0 \in C_{0,\sigma}^\infty (\Omega)$ such that
\begin{equation*}
\| f -f_0 \|_{L^{q,\infty}} < \varepsilon.
\end{equation*}
Using \eqref{c2}, we check that for $0<t <1$
\begin{multline*}
\| \mathrm{e}^{- t \mathscr{L}_{q,\infty}} f - f \|_{L^{q,\infty}} \\
 \leq \| \mathrm{e}^{- t \mathscr{L}_{q,\infty}} f - \mathrm{e}^{- t \mathscr{L}_{q,\infty}} f_0 \|_{L^{q,\infty}} + \| \mathrm{e}^{- t \mathscr{L}_{q,\infty}} f_0 - f_0 \|_{L^{q,\infty}} + \| f_0 - f  \|_{L^{q,\infty}} \\
 \leq C \| f - f_0 \|_{L^{q,\infty}} + C \| \mathrm{e}^{- t L_q} f_0 - f_0 \|_{L^q}\leq C \varepsilon \text{ as }t \to 0 + 0.
\end{multline*}
Here we used the fact that $\mathrm{e}^{- t L}$ is a $C_0$-semigroup on $L^q_\sigma ( \Omega )$. Since $\varepsilon$ is arbitrary, we find that
\begin{equation*}
\lim_{t \to 0 + 0} \| \mathrm{e}^{- t \mathscr{L}_{q,\infty}} f - f \|_{L^{q,\infty}} = 0.
\end{equation*}
Therefore we conclude that $\mathrm{e}^{- t \mathscr{L}_{q,\infty}}$ is a $C_0$-semigroup on $L^{q,\infty}_{0,\sigma} ( \Omega )$.

 We now assume that $1 < q \leq 3$. From Lemma \ref{lem213}, we find that for each $t>0$ and $f \in L^{q , \infty }_{0 , \sigma} ( \Omega ) $ 
\begin{equation*}
\| \nabla \mathrm{e}^{- t \mathscr{L}_{q , \infty}} f \|_{L^{q,\infty}} \leq C t^{- \frac{1}{2}} \| f \|_{L^{q,\infty}} .
\end{equation*}
By the previous argument, we see that $\nabla \mathrm{e}^{- t \mathscr{L}_{ q , \infty} } f \in L^{q , \infty}_0 ( \Omega ) $ for each $t>0$ and $f \in L^{q , \infty }_{0 , \sigma} ( \Omega ) $. Therefore the lemma follows.
\end{proof}

\subsection{$L^p$-$L^q$ Type Estimates for the Oseen Semigroup}\label{subsec24} 

Let us derive $L^p$-$L^q$ type estimates for the Oseen semigroup. Applying $L^p$-$L^q$ estimates for the Oseen kernel, the real interpolation theory, and {\color{red}{a}} duality argument, we obtain $L^p$-$L^q$ type estimates for the Oseen semigroup.\\

Throughout this subsection we assume that $|u_\infty| \leq \gamma$ for some $\gamma > 0$, and the symbols $L$ and $L^*$ represent the two linear operators defined by Subsection \ref{subsec23}.\\

We first state well-known $L^p$-$L^q$ estimates for the Oseen semigroup. Applying the real interpolation theory and the estimates for the Oseen semigroup obtained by \cite{KS98} and \cite{ES04}, we have 
\begin{lemma}[$L^{p}$-$L^{q}$ estimates for the Oseen semigroup (I)]\label{lem213} { \ }
\begin{enumerate}
\renewcommand{\labelenumi}{\rm{\rm{(\roman{enumi})}}}
\item Let $1 < p , q < \infty$ and $1 \leq r \leq \infty$ such that $p \leq q$. Then there is $C= C( \gamma , p,q,r )>0$ such that for all $t >0$ and $f \in L_\sigma^{p,r} ( \Omega )$
\begin{align}
& \| \mathrm{e}^{-t L} f \|_{L^{q,r} } \leq C t^{-\frac{3}{2}(\frac{1}{p} - \frac{1}{q})} \| f \|_{L^{p,r} } \label{e1} ,\\
& \| \mathrm{e}^{-t L^*} f \|_{L^{q,r} } \leq C t^{-\frac{3}{2}(\frac{1}{p} - \frac{1}{q})} \| f \|_{L^{p,r} }.\label{e2}
\end{align}
\item Let $1 < p , q \leq 3$ and $1 \leq r < \infty$ such that $p \leq q$. Then there is $C= C( \gamma , p,q,r )>0$ such that for all $t >0$ and $f \in L_\sigma^{p,r} ( \Omega )$
\begin{align}
& \| \nabla \mathrm{e}^{-t L} f \|_{L^{q,r} } \leq C t^{-\frac{1}{2} - \frac{3}{2}(\frac{1}{p} - \frac{1}{q})} \| f \|_{L^{p,r}},  \label{e3}\\
& \| \nabla \mathrm{e}^{-t L^*} f \|_{L^{q,r} } \leq C t^{-\frac{1}{2} - \frac{3}{2}(\frac{1}{p} - \frac{1}{q})} \| f \|_{L^{p,r}}.  \label{e4}
\end{align}
\item Let $1 < p , q \leq 3$ such that $p \leq q$. Then there is $C= C( \gamma , p,q )>0$ such that for all $t >0$ and $f \in L_{0 , \sigma}^{p,\infty} ( \Omega )$
\begin{align}
& \| \nabla \mathrm{e}^{-t L} f \|_{L^{q,\infty} } \leq C t^{-\frac{1}{2} - \frac{3}{2}(\frac{1}{p} - \frac{1}{q})} \| f \|_{L^{p,\infty}},  \label{eE3}\\
& \| \nabla \mathrm{e}^{-t L^*} f \|_{L^{q,\infty} } \leq C t^{-\frac{1}{2} - \frac{3}{2}(\frac{1}{p} - \frac{1}{q})} \| f \|_{L^{p,\infty}}.  \label{eE4}
\end{align}
\item Let $1 < p < \infty$. Then there is $C= C( \gamma , p )>0$ such that for all $t >0$ and $f \in L_\sigma^p ( \Omega )$
\begin{align}
& \| \mathrm{e}^{-t L} f \|_{L^\infty } \leq C t^{-\frac{3}{2 p} } \| f \|_{L^p} ,\label{e5}\\
& \| \mathrm{e}^{-t L^*} f \|_{L^\infty } \leq C t^{-\frac{3}{2 p} } \| f \|_{L^p} \label{e6}.
\end{align}
\item Let $1 < p , q \leq 3$ such that $1/p - 1/q=1/3$. Then there is $C = C( \gamma , p ) >0$ such that for all $f \in L^{p,1}_\sigma ( \Omega)$
\begin{align}
& \int_0^\infty \| \nabla \mathrm{e}^{- t L} f \|_{L^{q,1}} d t\leq C \| f \|_{L^{p,1}},\label{e7}\\
& \int_0^\infty \| \nabla \mathrm{e}^{- t L^* } f \|_{L^{q,1}} d t \leq C \| f \|_{L^{p,1}}\label{e8}.
\end{align}
\end{enumerate}
\end{lemma}
\noindent Combining the local energy decay of the Oseen semigroup derived by \cite{KS98} and \cite{ES04}, $L^p$-$L^q$ estimates for the Oseen semigroup, the interpolation theory, and an argument similar to that in \cite{ES05} and \cite[Sections 7,8]{HS09}, we deduce the properties $(\mathrm{i})$-$(\mathrm{iv})$ of Lemma \ref{lem213}. See \cite[Theorem 2.1]{Mar13} for the property $(\mathrm{iii})$ when $u_\infty = 0$. Note that $C_{0,\sigma}^\infty (\Omega)$ is dense in $L^{p,\infty}_{0,\sigma} (\Omega)$. See \cite[Section 5]{Shi01} and \cite[Lemma 8.4]{HS09} for the proof of the assertion $(\mathrm{v})$ of Lemma \ref{lem213}. See also \cite{Yam00} and \cite[Section 5]{Shi01}.

 Now we use Lemma \ref{lem213} and the real interpolation theory to prove the following lemma.
\begin{lemma}[$L^{p}$-$L^{q}$ estimates for the Oseen semigroup (II)]\label{lem214}{ \ }
\begin{enumerate}
\renewcommand{\labelenumi}{\rm{\rm{(\roman{enumi})}}}
\item Let $1 < p , q < \infty$ such that $p < q$. Then for each $1 \leq r_1 , r_2 \leq \infty$, there is $C= C( \gamma , p,q , r_1 ,r_2 )>0$ such that for all $t >0$ and $f \in [ L^{p, r_2} ( \Omega )]^3$
\begin{align}
& \| \mathrm{e}^{-t L} P f \|_{L^{q,r_1} } \leq C t^{-\frac{3}{2}(\frac{1}{p} - \frac{1}{q})} \| f \|_{L^{p,r_2} },\label{pp1}\\
& \| \mathrm{e}^{-t L^*} P f \|_{L^{q,r_1} } \leq C t^{-\frac{3}{2}(\frac{1}{p} - \frac{1}{q})} \| f \|_{L^{p,r_2} } \label{pp2}.
\end{align}
\item Let $1 < p , q < 3$ such that $p < q$. Then for each $1 \leq r_1,r_2 \leq \infty$, there is $C= C(\gamma , p,q, r_1 ,r_2  )>0$ such that for all $t >0$ and $f \in [ L^{p,r_2} ( \Omega ) ]^3$
\begin{align}
& \| \nabla \mathrm{e}^{-t L} P f \|_{L^{q,r_1} } \leq C t^{ - \frac{1}{2} -\frac{3}{2}(\frac{1}{p} - \frac{1}{q})} \| f \|_{L^{p, r_2} } \label{pp3} ,\\
& \| \nabla \mathrm{e}^{-t L^*} P f \|_{L^{q,r_1} } \leq C t^{-\frac{1}{2} - \frac{3}{2}(\frac{1}{p} - \frac{1}{q})} \| f \|_{L^{p,r_2} }.\label{pp4}
\end{align}
\end{enumerate}
\end{lemma}

\begin{proof}[Proof of Lemma \ref{lem214}]
Let $1 < p, q < \infty$ such $p < q$. We first show that there is $C = C (p,q) >0$ such that for every $f \in L^{q,\infty}_\sigma ( \Omega )$
\begin{align}
& \| \mathrm{e}^{-t L} P f \|_{L^{q,1} } \leq C t^{-\frac{3}{2}(\frac{1}{p} - \frac{1}{q})} \| f \|_{L^{p,\infty} },\label{pp5}\\
& \| \mathrm{e}^{-t L^*} P f \|_{L^{q,1} } \leq C t^{-\frac{3}{2}(\frac{1}{p} - \frac{1}{q})} \| f \|_{L^{p,\infty} } \label{pp6}.
\end{align}
Set $\theta=( q - p)/(q-1)$. It is easy to check that $0 < \theta <1$ and that $1/q=(1 - \theta )/p + \theta/(pq)$. Since
\begin{equation*}
(L^{p,\infty} (\Omega) , L^{p q , \infty} ( \Omega ))_{\theta, 1} = L^{q,1} (\Omega),
\end{equation*}
we use \eqref{e1} to see that
\begin{align*}
\| \mathrm{e}^{-t L} P f \|_{L^{q,1}} \leq & C \| \mathrm{e}^{- t L} P f \|_{L^{p,\infty}}^{1-\theta} \| \mathrm{e}^{- t L} P f \|_{L^{p q,\infty}}^{\theta}\\
\leq & C \| P f \|_{L^{p,\infty}}^\frac{p-1}{q-1} ( t^{-\frac{3}{2}(\frac{1}{p} - \frac{1}{p q})} \| P f \|_{L^{p , \infty}})^\frac{q-p}{q-1} \\
\leq & C t^{-\frac{3}{2} (\frac{1}{p} - \frac{1}{q} )}\| f \|_{L^{p,\infty}}.
\end{align*}
Thus, we see \eqref{pp5}. Similarly, we have \eqref{pp6}. Combining \eqref{eq22}, \eqref{pp5}, and \eqref{pp6}, we obtain \eqref{pp1} and \eqref{pp2}.

Since $(L^{p,\infty}(\Omega) , L^{3,\infty} ( \Omega ) )_{\theta , 1}  = L^{q ,1} (\Omega)$ when $\theta = \{3 (q-p) \}/\{ (3-p)q \}$, we repeat an argument similar to show $(\mathrm{i})$ to see $(\mathrm{ii})$. 
\end{proof}

\begin{lemma}[$L^{p}$-$L^{q}$ estimates for the Oseen semigroup (III)]\label{lem215}{ \ }
\begin{enumerate}
\renewcommand{\labelenumi}{\rm{(\roman{enumi})}}
\item Let $1 < p, q < \infty$ such that $p \leq q$. Then there is $C= C( \gamma , q)>0$ such that for all $t >0$ and $f \in [ L^p ( \Omega ) \cap L^1 ( \Omega )]^3$
\begin{align}
& \| \mathrm{e}^{-t L} P f \|_{L^q } \leq C t^{-\frac{3}{2} \left( 1 - \frac{1}{q} \right) } \| f \|_{L^1} ,\label{p1}\\
& \| \mathrm{e}^{-t L^*} P f \|_{L^q } \leq C t^{-\frac{3}{2} \left( 1 - \frac{1}{q} \right) } \| f \|_{L^1} \label{p2}.
\end{align}
\item Let $1 < p, q < \infty$ and $1 \leq r_1, r_2 \leq \infty$ such that $p <q$. Then there is $C= C( \gamma , q, r_1 )>0$ such that for all $t >0$ and $f \in [ L^{p, r_2 } ( \Omega ) \cap L^1 ( \Omega )]^3$
\begin{align}
& \| \mathrm{e}^{-t L} P f \|_{L^{q,r_1} } \leq C t^{-\frac{3}{2} \left( 1 - \frac{1}{q} \right) } \| f \|_{L^1} ,\label{p3}\\
& \| \mathrm{e}^{-t L^*} P f \|_{L^{q,r_1} } \leq C t^{-\frac{3}{2} \left( 1 - \frac{1}{q} \right) } \| f \|_{L^1} \label{p4}.
\end{align}
\item Let $1 < p < \infty$ and $1 \leq r \leq \infty$. Then there is $C= C( \gamma , p, r)>0$ such that for all $t >0$ and $f \in [ L^{p, r } ( \Omega ) ]^3$
\begin{align}
& \| \mathrm{e}^{-t L} P f \|_{L^\infty } \leq C t^{-\frac{3}{2 p} } \| f \|_{L^{p,r}} ,\label{p5}\\
& \| \mathrm{e}^{-t L^*} P f \|_{L^\infty } \leq C t^{-\frac{3}{2 p} } \| f \|_{L^{p,r}} \label{p6}.
\end{align}
\end{enumerate}
\end{lemma}
\begin{proof}[Proof of Lemma \ref{lem215}]
We first show $(\mathrm{i})$ and $(\mathrm{ii})$. Fix $t > 0$. Using \eqref{eq211}, the H\"{o}lder inequality, and \eqref{e6}, we see that
\begin{align*}
\| \mathrm{e}^{- t L} P f \|_{L^q} & \leq C \sup_{\varphi \in C^\infty_{0,\sigma} , { \ }\| \varphi \|_{L^{q'}} \leq 1} | \dual{ f , \mathrm{e}^{- t L^*} \varphi }  |\\
& \leq C \| f \|_{L^1} \sup_{\| \varphi \|_{L^{q'}} \leq 1} \| \mathrm{e}^{-t L^*} \varphi \|_{L^\infty}\\
& \leq C \| f \|_{L^1} t^{ -\frac{3}{2} \left( 1 -\frac{1}{q} \right) } .
\end{align*}
Here $1/q + 1/{q'} =1$. Thus, we have \eqref{p1}. Similarly, we see \eqref{p2}. Applying \eqref{p1}, \eqref{p2}, and \eqref{eq21}, we deduce $(\mathrm{ii})$.

Next we show $(\mathrm{iii})$. Fix $t >0$. By \eqref{eq26}, \eqref{eq23}, and \eqref{p4}, we observe that
\begin{align*}
\| \mathrm{e}^{- t L} P f \|_{L^\infty} & \leq C \sup_{\phi \in [C^\infty_0]^3 , { \ }\| \phi \|_{L^1} \leq 1} | \dual{  f , \mathrm{e}^{- t L^*} P \phi }  |\\
& \leq C \| f \|_{L^{p,r}} \sup_{\| \phi \|_{L^1} \leq 1} \| \mathrm{e}^{-t L^*} P \phi \|_{L^{p',r'}}\\
& \leq C \| f \|_{L^{p,r}} t^{ -\frac{3}{2} \left( 1 -\frac{1}{p'} \right) } ,
\end{align*}
which is \eqref{p5}. Similarly, we have \eqref{p6}. Therefore the lemma follows.
\end{proof}

Now we derive $L^p$-$L^q$ type estimates for the Oseen semigroup.
\begin{lemma}[$L^{p}$-$L^q$ type estimates for the Oseen semigroup]\label{lem216}{ \ }\\
Let $3/2 \leq p , q < \infty$ and $1 \leq r_1 , r_2 \leq \infty$ such that $p < q$. Let $j \in \{ 1,2,3 \}$. Then there is $C= C( \gamma , p , q , r_1,r_2 ) >0$ such that for all $f \in [ C_0^\infty ( \Omega ) ]^3$ and $t >0$
\begin{align}
& \| \mathrm{e}^{-t L} P \partial_j f \|_{L^{q,r_1} } \leq C t^{-\frac{1}{2} - \frac{3}{2}(\frac{1}{p} - \frac{1}{q})} \| f \|_{L^{p,r_2}} \label{e15},\\
& \| \mathrm{e}^{-t L^*} P \partial_j f \|_{L^{q,r_1} } \leq C t^{-\frac{1}{2} - \frac{3}{2}(\frac{1}{p} - \frac{1}{q})} \| f \|_{L^{p,r_2}}.  \label{e16}
\end{align}
\end{lemma}
\begin{proof}[Proof of Lemma \ref{lem216}]
Set $\ell = (p+q)/2$. Applying \eqref{eq22}, \eqref{eq210}, \eqref{eq23}, and \eqref{pp2}, we check that
\begin{align}
\| \mathrm{e}^{- t L} P \partial_j f \|_{ L^{q , r_1} } & \leq C \| \mathrm{e}^{- t L} P \partial_j f \|_{L^{q,1}} \notag \\
& \leq C \sup_{ \varphi \in L^{ q' , \infty }_\sigma , { \ } \| \varphi \|_{ L^{q' , \infty}}  \leq 1 } | \dual{ \mathrm{e}^{ - \frac{t}{2} L} P \partial_j f  , \mathrm{e}^{- \frac{t}{2} L^*} \varphi  } | \notag \\
& \leq C \| \mathrm{e}^{ - \frac{t}{2} L } P \partial_j f \|_{L^{\ell , \infty}} \sup_{\| \varphi \|_{L^{q' , \infty}} \leq 1 } \| \mathrm{e}^{- \frac{t}{2} L^* } \varphi \|_{L^{\ell' , 1}} \notag \\
& \leq C t^{- \frac{3}{2} \left( \frac{1}{q'} - \frac{1}{\ell'} \right) } \| \mathrm{e}^{- \frac{t}{2}L } P \partial_j f \|_{L^{\ell , \infty }} \label{abc1}.
\end{align}
Here $1/\ell + 1/{\ell'} =1$ and $1/q + 1/{q'} =1$. Note that $\ell' >q'$ since $q > \ell $. Using \eqref{eq210}, \eqref{eq23}, and \eqref{e4}, we observe that
\begin{align}
\| \mathrm{e}^{- \frac{t}{2} L} P \partial_j f \|_{ L^{\ell , \infty} } & \leq C \sup_{ \varphi \in C^\infty_{0,\sigma} , { \ } \| \varphi \|_{ L^{\ell' , 1}}  \leq 1 } | \dual{ f  , - \partial_j \mathrm{e}^{- \frac{t}{2} L^*} \varphi  } | \notag \\
& \leq C \| f \|_{L^{p , r_2}} \sup_{\| \varphi \|_{L^{\ell' , 1}} \leq 1 } \| \nabla \mathrm{e}^{- \frac{t}{2} L^* } \varphi \|_{L^{p' , r_2'}} \notag \\
& \leq C \| f \|_{L^{p , r_2}} \sup_{\| \varphi \|_{L^{\ell' , 1}} \leq 1 } \| \nabla \mathrm{e}^{- \frac{t}{2} L^* } \varphi \|_{L^{p' , 1}} \notag \\
& \leq C t^{- \frac{1}{2} - \frac{3}{2} \left( \frac{1}{\ell'} - \frac{1}{p'} \right) } \| f \|_{L^{p , r_2}} \label{abc2}.
\end{align}
Note that $p' \leq 3$ since $3/2 \leq p$. Combining \eqref{abc1} and \eqref{abc2} gives
\begin{equation*}
\| \mathrm{e}^{- t L} P \partial_j f \|_{ L^{q , r_1} } \leq C t^{- \frac{1}{2} - \frac{3}{2} \left( \frac{1}{p} - \frac{1}{q} \right) } \| f \|_{L^{p , r_2}}.
\end{equation*}
Thus, we have \eqref{e15}. Similarly, we obtain \eqref{e16}.
\end{proof}

\subsection{Fractional Powers of the Oseen Operator in the Lorentz Space}\label{subsec25} 
In this subsection we study fractional powers of the Oseen operator. Let $L_{q,r}$ and $L_{q,r}^*$ be the two operators defined by Subsection \ref{subsec23}. Since $\mathrm{e}^{-t L_{q,r}}$ is a uniformly bounded $C_0$-semigroup on $L^{q,r}_\sigma (\Omega)$ if $1 <q < \infty$ and $1 \leq r < \infty$, one can define fractional powers of the operator $(L_{q,r} + \eta)$ for $\eta \geq 0$. More precisely, for each $0< \zeta <1$, we define $L_{q,r}^\zeta f := \lim_{\eta \to 0+0}(L_{q,r} + \eta )^\zeta f$, $D (L_{q,r}^\zeta):= D ( ( L_{q,r} + \eta)^\zeta)$. Since $\mathrm{e}^{- t L_{q,r}}$ is an analytic semigroup on $L^{q,r}_\sigma (\Omega)$, we apply \cite[Theorem 6.13 in Chapter 2]{Paz83}, we obtain
\begin{lemma}\label{lem217}
Let $1 <q < \infty$ and $1 \leq r < \infty$. Then for each $0 < \zeta \leq 1$ and $\eta >0$ there is $C = C (\gamma, q, r, \zeta , \eta ) >0$ such that for each $t > 0$, $f \in L^{q,r}_\sigma ( \Omega )$, $g \in D (L_{q,r}^\zeta )$, and $h \in D ( (L^*_{q,r} )^\zeta)$,
\begin{align}
& \| (L_{q,r} + \eta )^\zeta \mathrm{e}^{- t L_{q,r}} f \|_{L^{q,r}} \leq C t^{- \zeta} \mathrm{e}^{\eta t} \| f \|_{L^{q,r}}, \label{kk1}\\
& \| (L_{q,r}^* + \eta )^\zeta \mathrm{e}^{- t L^*_{q,r}} f \|_{L^{q,r}} \leq C t^{- \zeta} \mathrm{e}^{\eta t} \| f \|_{L^{q,r}}, \label{kk2}\\
& \| \mathrm{e}^{ - t L_{q,r} } g - g \|_{L^{q,r}} \leq C t^\zeta \mathrm{e}^{ \eta t} \| (L_{q,r} + \eta )^\zeta g \|_{L^{q,r}}, \label{kk3}\\
& \| \mathrm{e}^{ - t L^*_{q,r} } h - h \|_{L^{q,r}} \leq C t^\zeta \mathrm{e}^{\eta t} \| (L^*_{q,r} + \eta )^\zeta h \|_{L^{q,r}}. \label{kk4}
\end{align}
\end{lemma}
\noindent See also \cite[Lemma 3.13 and the proof]{Kob13} for details.

We introduce fractional powers of the Oseen operator $L_{q,\infty}$ in $L_{\sigma}^{q,\infty} (\Omega)$. Let $1 < q,q' < \infty$ such that $1/q + 1/{q'}=1$. For each $0 < \zeta <1$ and $\eta > 0$, we define $(L_{q,\infty} + \eta )^\zeta$ as follows: for all $g \in D ( (L_{q' ,1}^* )^\zeta)$
\begin{equation*}
\dual{ (L_{q,\infty} + \eta )^\zeta f , g } := \dual{ f , (L_{q',1}^*  + \eta )^\zeta g }
\end{equation*}
and
\begin{equation*}
D ( ( L_{q,\infty} + \eta)^\zeta ) := \{ f \in L^{q,\infty}_\sigma ( \Omega ) ; \| ( L_{q,\infty}  + \eta )^\zeta f \|_{L^{q,\infty}}  < + \infty \} . 
\end{equation*}
Moreover, for $0 < \zeta < 1$, we define $L_{q,\infty}^\zeta$ in $L^{q,\infty}_\sigma$ as follows: $g \in D ( (L_{q' ,1}^* )^\zeta)$
\begin{equation*}
\dual{ L_{q,\infty}^\zeta f , g } := \lim_{\eta \to 0 + 0} \dual{ f , (L_{q',1}^*  + \eta )^\zeta g } \text{ and } D (L^\zeta_{q,\infty}) : = D( (L_{q,\infty} + 1 )^\zeta ) .
\end{equation*}
Note that $D ( ( L_{q,\infty} + \eta )^\zeta)$ is not empty. In fact, by definition and Lemma \ref{lem29}, we check that for each $f \in D(L_{q,\infty} )$
\begin{multline*}
\| ( L_{q,\infty}  + \eta ) ^\zeta f \|_{L^{q,\infty}} \leq C \sup_{\varphi \in C_{0,\sigma}^\infty , { \ }\| \varphi \|_{L^{q',1}} \leq 1} |\dual{ f , (L_{q',1}^*  + \eta )^\zeta \varphi }|\\
= C \sup_{\| \varphi \|_{L^{q',1}} \leq 1} |\dual{ f , (L_{q',1}^*  + \eta ) (L_{q',1}^* + \eta )^{ - (1 - \zeta )} \varphi }|\\
= C \sup_{\| \varphi \|_{L^{q',1}} \leq 1} |\dual{ (L_{q,\infty}  + \eta) f , (L_{q',1}^*  + \eta )^{ - ( 1 - \zeta)} \varphi }|\\
\leq C (\eta ) \| ( L_{q,\infty} + \eta ) f \|_{L^{q,\infty}} .
\end{multline*}
Here we used that fact that $\| (L^*_{q',1}+ \eta )^{- (1 -\zeta)} \varphi \|_{L^{q',1}} \leq C (\zeta, \eta) \| \varphi \|_{L^{q',1}}$. 

Now we study fractional powers of the operator $L_{q,\infty}$.
\begin{lemma}\label{lem218}
Let $1 <q < \infty$. Then for each $0 < \zeta \leq 1$ and $\eta >0$ there is $C = C (\gamma, q, \zeta , \eta ) >0$ such that for each $t > 0$, $f \in L^{q,\infty}_\sigma ( \Omega )$, $g \in D (L_{q,\infty}^\zeta )$, and $h \in D ( (L^*_{q,\infty} )^\zeta)$,
\begin{align}
& \| (L_{q,\infty} + \eta )^\zeta \mathrm{e}^{- t L_{q,\infty}} f \|_{L^{q,\infty}} \leq C t^{- \zeta} \mathrm{e}^{\eta t} \| f \|_{L^{q,\infty}}, \label{kk5}\\
& \| (L_{q,\infty} + \eta )^\zeta \mathrm{e}^{- t L^*_{q,\infty}} f \|_{L^{q,\infty}} \leq C t^{- \zeta} \mathrm{e}^{\eta t} \| f \|_{L^{q,\infty}}, \label{kk6}\\
& \| \mathrm{e}^{ - t L_{q,\infty} } g - g \|_{L^{q,\infty}} \leq C t^\zeta \mathrm{e}^{ \eta t} \| (L_{q,\infty} + \eta )^\zeta g \|_{L^{q,\infty}}, \label{kk7}\\
& \| \mathrm{e}^{ - t L^*_{q,\infty} } h - h \|_{L^{q,\infty}} \leq C t^\zeta \mathrm{e}^{\eta t} \| (L^*_{q,\infty} + \eta )^\zeta h \|_{L^{q,\infty}} \label{kk8}.
\end{align}
\end{lemma}
\begin{proof}[Proof of Lemma \ref{lem218}]
We first show \eqref{kk5}. By definition, \eqref{eq211}, \eqref{eq23}, and \eqref{kk2}, we see that
\begin{align*}
\| ( L_{q,\infty}  + \eta ) ^\zeta \mathrm{e}^{- t L_{q,\infty}} f \|_{L^{q,\infty}} \leq C \sup_{\varphi \in C_{0,\sigma}^\infty , { \ }\| \varphi \|_{L^{q',1}} \leq 1} |\dual{ f , \mathrm{e}^{- t L^*_{q',1} } (L_{q',1}^*  + \eta )^\zeta \varphi }|\\
\leq C \| f \|_{L^{q,\infty}} \sup_{\| \varphi \|_{L^{q',1}} \leq 1} \| (L_{q',1}^*  + \eta )^\zeta \mathrm{e}^{- t L^*_{q',1}} \varphi  \|_{L^{q',1}} \\
\leq C t^{- \zeta}\mathrm{e}^{\eta t} \| f \|_{L^{q,\infty}},
\end{align*}
where $1/q+ 1/{q'}=1$. Thus, we see \eqref{kk5}, Similarly, we have \eqref{kk6}. 

Next we show \eqref{kk7}. Using Lemmas \ref{lem29} and \ref{lem217}, we check that
\begin{multline*}
\| \mathrm{e}^{- t L_{q,\infty}} f -f \|_{L^{q,\infty}} \leq C \sup_{\varphi \in C_{0,\sigma}^\infty , { \ }\| \varphi \|_{L^{q',1}} \leq 1} |\dual{ f , ( \mathrm{e}^{- t L^*_{q',1} } -1 )  \varphi }|\\
= C \sup_{\| \varphi \|_{L^{q',1}} \leq 1} |\dual{ ( L_{q,\infty} + \eta )^\zeta f , ( \mathrm{e}^{- t L^*_{q',1} } -1 ) (L^*_{q,1} + \eta )^{-\zeta}  \varphi }|\\
\leq C t^\zeta \mathrm{e}^{\eta t} \| (L_{q,\infty} + \eta )^\zeta f \|_{L^{q,\infty}} .
\end{multline*}
Therefore we see \eqref{kk7}. Similar, we obtain \eqref{kk8}.
\end{proof}

Finally, we prove a key lemma.
\begin{lemma}\label{lem219}
Let $1 < q \leq 3$, $\eta >0 $, and $1/2 < \zeta < 1$. Then there is $ C = C (q, \eta , \zeta) >0$ such that for all $f \in D (A_{q,1})$
\begin{align}
\| \nabla f \|_{L^{q,1}} \leq C \| ( L_{q,1} + \eta )^{\zeta} f \|_{L^{q,1}} \label{tt1},\\
\| \nabla f \|_{L^{q,1}} \leq C \| ( L^*_{q,1} + \eta )^{\zeta} f \|_{L^{q,1}} .\label{tt2}
\end{align}
\end{lemma}
To prove Lemma \ref{lem219}, we prepare the following two lemmas.
\begin{lemma}\cite[Theorem 2.1]{Yam00} \label{ss2}
Let $1 < q \leq 3$. Then there is $ C = C (q) >0$ such that for all $f \in D (A_{q,1}^{1/2})$
\begin{equation}\label{tt3}
\| \nabla f \|_{L^{q,1}} \leq C \| A_{q,1}^{1/2} f \|_{L^{q,1}}.
\end{equation}
\end{lemma}
\begin{lemma}\cite[Theorem 6.12 in Chap.2]{Paz83}\label{ss4}
Let $X$ be a Banach space and $\| \cdot \|_X$ its norm. Let $\mathscr{A} : D(\mathscr{A}) (\subset X ) \to X$ be a closed linear operator densely defined in $X$. Let $\mathscr{B} : D ( \mathscr{B}) (\subset X) \to X$ be a closed linear operator on $X$ such that $D(\mathscr{A}) \subset D ( \mathscr{B})$. Let $0 < \alpha , \zeta < 1$ such that $\alpha < \zeta$. Suppose that $0 \in \rho ( \mathscr{A})$. Assume that there is $C >0$ such that for all $\varepsilon >0$ and $f \in D (\mathscr{A})$ 
\begin{equation*}
\| \mathscr{B} f \|_X \leq \varepsilon^{1 - \alpha} \| \mathscr{A} f \|_X + C \varepsilon^{ - \alpha} \| f \|_X .
\end{equation*}
Then there is $C = C (\zeta) >0$ such that for all $f \in D ( \mathscr{A})$
\begin{equation*}
\| \mathscr{B} f \|_X \leq C \| \mathscr{A}^\zeta  f \|_X .
\end{equation*}
\end{lemma}
Let us attack Lemma \ref{lem219}.
\begin{proof}[Proof of Lemma \ref{lem219} ]
We only show \eqref{tt1} since the proof to derive \eqref{tt2} is similar.

We first show that for each $\eta >0$ there is $C = C (\eta ) >0$ such that for all $f \in D (A_{q,1})$
\begin{equation}
\| A_{q,1} f  \|_{L^{q,1}} \leq C \| (L_{q,1} + \eta ) f \|_{L^{q,1}} \label{tt5} .
\end{equation}
Fix $f \in D (A_{q,1})$. Applying \eqref{tt3}, the moment inequality, and the Young inequality, we check that
\begin{align*}
\| A_{q,1} f \|_{L^{p,1}} & \leq \| (L_{q,1} + \eta ) f \|_{L^{q,1}} + C \| \nabla f \|_{L^{q,1}} + C \| f \|_{L^{}q,1}\\
& \leq C \| (L_{q,1} + \eta ) f \|_{L^{q,1}} + C \| A_{q,1}^{1/2} f \|_{L^{q,1}} + C \| f \|_{L^{q,1}}\\
& \leq C \| (L_{q,1} + \eta ) f \|_{L^{q,1}} + \frac{1}{2} \| A_{q,1} f \|_{L^{q,1}} + C \| f \|_{L^{q,1}}.
\end{align*}
Since $\| f \|_{L^{q,1}} \leq C \| (L_{q,1} + \eta  ) f \|_{L^{q,1}}$, we have \eqref{tt5}.

Next we prove that there is $C >0$ such that for all $ f \in D (A_{q,1})$
\begin{equation*}
\| A^{1/2}_{q,1} f \|_{L^{q,1}} \leq \varepsilon^{1/2} \| (L_{q,1} + \eta ) f \|_{L^{q,1}} + C\frac{1}{\varepsilon^{1/2}} \| f \|_{L^{q,1}} .
\end{equation*}
Fix $f \in D (A_{q,1})$. Using \eqref{tt3} the moment inequality, we see that
\begin{align*}
\| A_{q,1}^{1/2} f \|_{L^{q,1}} & \leq C \| A_{q,1} f \|_{L^{q,1}}^\frac{1}{2} \| f \|_{L^{q,1}}^\frac{1}{2} \\
& \leq C \| (L_{q,1} + \eta ) f \|_{L^{q,1}}^\frac{1}{2} \| f \|_{L^{q,1}}^\frac{1}{2} = : \text{ (RHS)}.
\end{align*}
The Young inequality shows that for each $\varepsilon > 0$ 
\begin{equation*}
\text{(RHS)} \leq \varepsilon^\frac{1}{2} \| (L_{q,1} + \eta ) f \|_{L^{q,1}} + C \varepsilon^{- \frac{1}{2}} \| f \|_{L^{q,1}} .
\end{equation*}
Applying Lemma \ref{ss4}, we see that for each $1/2 < \zeta < 1$ there is $C = C (\zeta) >0$ such that for all $f \in D (A_{q,1})$
\begin{equation*}
\| \nabla f \|_{L^{q,1}} \leq C \| ( L_{q,1} + \eta )^{\zeta} f \|_{L^{q,1}}.
\end{equation*}
Therefore the lemma follows.
\end{proof}

\subsection{Key estimates for the Oseen semigroup}\label{subsec26} 
Throughout this subsection we assume that $|u_\infty| \leq \gamma$ for some $\gamma > 0$, and the symbols $L$ and $L^*$ represent the two linear operators defined by Subsection \ref{subsec23}. The aim of this subsection is to prove the following lemma.
\begin{lemma}\label{lem223}{ \ }\\
\noindent $(\mathrm{i})$ Let $1/2 < \zeta < 1$ and $0 < \alpha < 1/2$ such that $\zeta + \alpha <1$. Let $1 <p , q \leq 3$ and $1 \leq r \leq \infty$ such that $p < q$. Then there is $C = C (\gamma, \zeta, \alpha, p , q ,r ) >0$ such that for all $t_1,t_2 >0$ and $\phi \in [L^{p,r} ( \Omega )]^3 $ 
\begin{align}
& \| \nabla \{ (\mathrm{e}^{-t_2 L} -1) \mathrm{e}^{-t_1 L} P \phi \} \|_{L^{q,1}} \leq C t_2^\alpha \mathrm{e}^{\frac{t_1}{2} + t_2} {(t_1)}^{- \zeta - \alpha - \frac{3}{2} \left( \frac{1}{p} - \frac{1}{q} \right)}  \| \phi \|_{L^{p, r }} ,\label{6a}\\
& \| \nabla \{ (\mathrm{e}^{-t_2 L^*} -1) \mathrm{e}^{-t_1 L^*} P \phi \} \|_{L^{q,1}} \leq Ct_2^\alpha \mathrm{e}^{\frac{t_1}{2} + t_2 } {(t_1)}^{- \zeta - \alpha - \frac{3}{2} \left( \frac{1}{p} - \frac{1}{q} \right)}  \| \phi \|_{L^{p, r }}.\label{6b}
\end{align}
\noindent $(\mathrm{ii})$ Let $j \in \{ 1 , 2 , 3 \}$, $0 < \alpha <1$, and $3/2 \leq q \leq 3$. Then there is $C = C (\gamma, \alpha, q) >0$ such that for all $t_1 , t_2 >0$ and $\phi \in [ C_0^\infty (\Omega) ]^3$
\begin{align}
& \| (\mathrm{e}^{- t_2 L} - 1 ) \mathrm{e}^{- t_1 L} P \partial_j \phi \|_{L^{q ,1}} \leq C t_2^\alpha \mathrm{e}^{ \frac{t_1}{2} + t_2} {(t_1)}^{- \alpha - \frac{3}{2}\left( 1 -\frac{1}{q} \right) } \| \phi \|_{L^{\frac{3}{2} , 1}} ,\label{6c}\\
& \| (\mathrm{e}^{- t_2 L^*} - 1 ) \mathrm{e}^{- t_1 L^*} P \partial_j \phi \|_{L^{q ,1}} \leq C t_2^\alpha \mathrm{e}^{ \frac{t_1}{2} + t_2} {(t_1)}^{- \alpha - \frac{3}{ 2}\left( 1 - \frac{1}{q} \right)} \| \phi \|_{L^{\frac{3}{2} , 1}} .\label{6d}
\end{align}
\noindent $(\mathrm{iii})$ Let $0 < \alpha <1$ and $1 <q < \infty$. There is $C = C (\gamma , \alpha , q )  >0$ such that for all $t_1,t_2 >0$ and $\phi \in [C_0^\infty ( \Omega )]^3$
\begin{align}
& \| ( \mathrm{e}^{- t_2 L} - 1 ) \mathrm{e}^{- t_1 L} P \phi \|_{L^{q,1}} \leq  \ C t_2^\alpha \mathrm{e}^{\frac{t_1}{2} +t_2} {(t_1)}^{- \alpha -\frac{3}{2} \left( 1 -\frac{1}{q} \right) } \| \phi \|_{L^1} \label{6e},\\
& \| ( \mathrm{e}^{- t_2 L^*} - 1 ) \mathrm{e}^{- t_1 L^*} P \phi \|_{L^{q,1}} \leq C t_2^\alpha \mathrm{e}^{\frac{t_1}{2} +t_2 } { (t_2)}^{- \alpha -\frac{3}{2}\left( 1 - \frac{1}{q} \right) } \| \phi \|_{L^1} \label{6f}.
\end{align}
\end{lemma}

\begin{proof}[Proof of Lemma \ref{lem223}]
We first derive \eqref{6a} to prove $(\mathrm{i})$. Fix $t_1,t_2 >0$ and $\phi \in [L^{p,r} ( \Omega )]^3 $. Using \eqref{tt1}, \eqref{kk3}, \eqref{eq210}, \eqref{pp1}, and \eqref{kk2}, we see that
\begin{multline*}
\| \nabla \{ (\mathrm{e}^{-t_2 L } -1) \mathrm{e}^{- t_1 L} P \phi \} \|_{L^{q,1}} \leq C \| (\mathrm{e}^{-t_2 L } -1) \mathrm{e}^{- \frac{t_1}{2} L} (L + 1)^{\zeta} \mathrm{e}^{-\frac{t_1}{2} L} P \phi \|_{L^{q,1}}\\
\leq C t_2^\alpha \mathrm{e}^{t_2} \|  \mathrm{e}^{- \frac{t_1}{2} L} (L + 1)^{\zeta + \alpha} \mathrm{e}^{- \frac{t_1}{2} L} P \phi \|_{L^{q,1}}\\
\leq C t_2^\alpha \mathrm{e}^{t_2} \sup_{\varphi \in L^{q',\infty}_\sigma, { \ }\| \varphi \|_{L^{q',\infty} } \leq 1} | \dual{ \mathrm{e}^{- \frac{t_1}{2} L} P \phi , (L^* + 1)^{\zeta + \alpha} \mathrm{e}^{- \frac{t_1}{2} L^*} \varphi} | \\
\leq C t_2^\alpha \mathrm{e}^{t_2} \| \mathrm{e}^{- \frac{t_1}{2} L}  P \phi \|_{L^{q,1}}  \sup_{\| \varphi \|_{L^{q',\infty} } \leq 1} \| (L^* + 1)^{\zeta + \alpha} \mathrm{e}^{- \frac{t_1}{2} L^*} \varphi \|_{L^{q',\infty}}\\
\leq C t_2^\alpha \mathrm{e}^{t_2} t_1^{ - \frac{3}{2} \left( \frac{1}{p} - \frac{1}{q} \right)} \| \phi \|_{L^{p,r}} t_1^{- \zeta - \alpha} \mathrm{e}^{ \frac{t_1}{2}} . 
\end{multline*}
Here $1/q + 1/{q'}=1$. Thus, we have \eqref{6a}. Similarly, we obtain \eqref{6b}.

Secondly, we prove \eqref{6c} to show $(\mathrm{ii})$. Using \eqref{kk3}, \eqref{eq210}, \eqref{e15}, and \eqref{kk1}, we check that
\begin{multline*}
\| ( \mathrm{e}^{- t_2 L} - 1 ) \mathrm{e}^{- t_1 L} P \partial_j \phi \|_{L^{q , 1}} \leq C t_2^\alpha \mathrm{e}^{t_2} \| \mathrm{e}^{ - \frac{t_1}{2} L} (L + 1)^\alpha \mathrm{e}^{- \frac{t_1}{2} L} P \partial_j \phi  \|_{L^{q,1}}\\
\leq C t_2^\alpha \mathrm{e}^{t_2} \sup_{\varphi \in L^{q',\infty}_\sigma, { \ }\| \varphi \|_{L^{q',\infty}} \leq 1 } |\dual{ \mathrm{e}^{- \frac{t_1}{2} L} P \partial_j \phi , (L^* + 1)^\alpha \mathrm{e}^{- \frac{t_1}{2} L^*} \varphi }|\\
\leq C t_2^\alpha \mathrm{e}^{t_2} \| \mathrm{e}^{- \frac{t_1}{2} L } P \partial_j \phi \|_{L^{q, 1}} \sup_{\| \varphi \|_{L^{q',\infty}} \leq 1} \| (L^* +  1)^\alpha \mathrm{e}^{ - \frac{t_1}{2} L^*} \varphi \|_{L^{q',\infty}}\\
\leq C t_2^\alpha \mathrm{e}^{t_2} t_1^{ -\frac{1}{2} -\frac{3}{2}\left( \frac{2}{3} -\frac{1}{q} \right) } \| \phi \|_{L^{\frac{3}{2}, 1}} t_1^{-\alpha} \mathrm{e}^{\frac{t_1}{2}} . 
\end{multline*}
Therefore we have \eqref{6c}. Similarly, we see \eqref{6d}.

Finally, we attack \eqref{6e} to prove $(\mathrm{iii})$. By \eqref{kk3}, \eqref{eq210}, \eqref{pp1}, \eqref{p3}, and \eqref{kk2}, we check that
\begin{multline*}
\| ( \mathrm{e}^{- t_2 L} - 1 ) \mathrm{e}^{- t_1 L} P \phi \|_{L^{q,1}} \leq C t_2^\alpha \mathrm{e}^{t_2} \| \mathrm{e}^{- \frac{t_1}{2} L } (L + 1)^\alpha \mathrm{e}^{- \frac{t_1}{2} L }  P \phi \|_{L^{q,1}}\\
\leq C t_2^\alpha \mathrm{e}^{t_2} \sup_{\varphi \in L_\sigma^{q',\infty},{ \ }\| \varphi \|_{L^{q',\infty}} \leq 1 } | \dual{ \mathrm{e}^{ - \frac{t_1}{2} L } P \phi , (L^* + 1)^\alpha \mathrm{e}^{- \frac{t_1}{2} L^*} \varphi } |\\
\leq C t_2^\alpha \mathrm{e}^{t_2} \| \mathrm{e}^{ - \frac{t_1}{2} L } P \phi \|_{L^{q,1}}\sup_{\| \varphi \|_{L^{q',\infty}} \leq 1 } \| (L^* + 1)^\alpha \mathrm{e}^{- \frac{t_1}{2} L^*} \varphi \|_{L^{q',\infty}}\\
\leq C t_2^\alpha \mathrm{e}^{t_2} t_1^{ - \frac{3 (q - 1)}{2 q} } \| \phi \|_{L^1} t_1^{- \alpha} \mathrm{e}^{ \frac{t_1}{2}}.
\end{multline*}
Thus, we see \eqref{6e}. Similarly, we have \eqref{6f}. Therefore the lemma follows.
\end{proof}

 \section{Global-in-Time Generalized Mild Solution}\label{sect3}

In this section we establish the existence of a unique global-in-time generalized mild solution of the system \eqref{eq15} when $v_0 \in L^{3,\infty}_\sigma (\Omega)$ and both $\| v_0 \|_{L^{3,\infty}}$ and $\| w \|_{L^{3,\infty}}$ are sufficiently small. Moreover, we investigate the asymptotic stability of the generalized mild solution. To this end, we first discuss both the uniqueness and the weak continuity of the generalized mild solutions of \eqref{eq15}, and we show the existence of a unique global-in-time generalized mild solution of \eqref{eq15} to prove Theorem \ref{thm19}. Secondly, we derive the asymptotic stability of the generalized mild solution to prove Theorem \ref{thm110}.

Throughout this section we assume that $w$ is as in Assumption \ref{assA} and that $|u_\infty| \leq \gamma$ for some $\gamma > 0$. Let $3 < p < \infty$ and $1 < p' < 3/2$ such that $1/p + 1/{p'}=1$. Fix $p$ and $p'$. The symbols $L$ and $L^*$ represent the two linear operators defined by Subsection \ref{subsec23}.

\subsection{Global-in-Time Generalized Mild Solution}\label{subsec31} 
We first study {\color{red}{the uniqueness and the weak continuity}} of the generalized mild solutions of the system \eqref{eq15}.
\begin{lemma}[Uniqueness (I)]\label{lem31}
Let $T \in (0, \infty]$. Suppose that $a, b \in L^{3,\infty}_{\sigma}( \Omega )$ and that
\begin{equation*}
v,u \in BC ((0, T) ; L^{3,\infty}_{\sigma} (\Omega ) ) .
\end{equation*}
Set
\begin{equation}\label{eq31}
\delta : = \max \{ \sup_{0 < t < T}\{ \|  v (t)\|_{L^{3,\infty}} \} , { \ }\sup_{0 < t < T}\{ \| u (t)\|_{L^{3,\infty}} \} ,{ \ } \| w \|_{L^{3,\infty}} \}.
\end{equation}
Assume that for every $0 <t <T$ and $\varphi \in C_{0,\sigma}^\infty (\Omega)$
\begin{align*}
& \dual{v (t) , \varphi} = \dual{a , \mathrm{e}^{- t L^*}\varphi} + \int_0^t \dual{v \otimes v + v\otimes w + w \otimes v , \nabla \mathrm{e}^{-(t-s) L^*} \varphi}  d s,\\
& \dual{u (t) , \varphi} = \dual{b , \mathrm{e}^{- t L^*} \varphi} + \int_0^t \dual{u \otimes u + u\otimes w + w \otimes u , \nabla \mathrm{e}^{-(t-s) L^*} \varphi}  d s.
\end{align*}
Then there are $\delta_\dagger = \delta_\dagger (\gamma )>0$ and {\color{red}{$K_\dagger = K_\dagger (\gamma) >0$}} such that if $\delta \leq \delta_\dagger$ then
\begin{align}
\sup_{0<t <T}\{ \| v(t) - u(t) \|_{L^{3,\infty}} \} \leq & K_\dagger \| a - b \|_{L^{3,\infty}},\label{eq32}\\
\sup_{0 < t <T } \{ \| v ( t ) \|_{L^{3, \infty}} \}  \leq & K_\dagger \| a \|_{L^{3,\infty}}, \label{eq33}\\
\sup_{0 < t <T } \{ \| u ( t ) \|_{L^{3, \infty}} \}  \leq & K_\dagger \| b \|_{L^{3,\infty}}. \label{eq34}
\end{align}
\end{lemma}
\begin{lemma}[Weak continuity]\label{lem32}
Under the same assumptions of Lemma \ref{lem31}, for each $\psi \in L^{3/2 , 1}_\sigma ( \Omega )$
\begin{align}
\lim_{t \to 0 + 0} \dual{v(t) , \psi } = \dual{a , \psi},\label{eq35}\\
\lim_{t \to 0 + 0} \dual{u(t) , \psi } = \dual{b , \psi}.\label{eq36}
\end{align}
\end{lemma}
\begin{lemma}[Uniqueness (II)]\label{lem33}
Under the same hypotheses of Lemma \ref{lem31}, we assume in addition that
\begin{equation*}
v,u \in C ((0,T); L^{p,\infty}_\sigma (\Omega )).
\end{equation*}
Then there are $\delta_{\dagger \dagger} = \delta_{\dagger \dagger} (\gamma ,p )>0$ and $K_{\dagger \dagger} = K_{\dagger \dagger} (\gamma , p) >0$ such that if $\delta \leq \delta_{ \dagger \dagger}$ then
\begin{align}
\sup_{0<t <T}\{ \| v(t) - u(t) \|_{L^{3,\infty}} \} \leq & K_{\dagger \dagger} \| a - b \|_{L^{3,\infty}},\label{eq37}\\
\sup_{0<t <T} \{t^\frac{p -3}{2 p} \| v(t) - u (t) \|_{L^{p,\infty}} \} \leq & K_{\dagger \dagger} \| a - b \|_{L^{3,\infty}},\label{eq38}\\
\sup_{0 < t <T } \{ \| v ( t ) \|_{L^{3, \infty}} \} + \sup_{0<t <T} \{t^\frac{p -3}{2 p} \| v(t) \|_{L^{p,\infty}} \} \leq & K_{\dagger \dagger} \| a \|_{L^{3,\infty}}, \label{eq39}\\
\sup_{0 < t <T } \{ \| u ( t ) \|_{L^{3, \infty}} \} + \sup_{0<t <T} \{t^\frac{p -3}{2 p} \| u (t) \|_{L^{p,\infty}} \} \leq & K_{\dagger \dagger} \| b \|_{L^{3,\infty}}. \label{eq310}
\end{align}
Here $\delta$ is {\color{red}{the constant}} defined by \eqref{eq31}. Moreover, assume in addition that $a , b \in L^{3,r}_\sigma (\Omega )$ for some $1 < r < \infty$. Then the two assertions hold:\\
\noindent $(\mathrm{i})$
\begin{equation*}
v,u \in B C ( (0 , T ) ; L^{3,r}_\sigma ( \Omega ) ).
\end{equation*}
\noindent $(\mathrm{ii})$ There is $K_{\dagger \dagger \dagger} = K_{\dagger \dagger \dagger}  ( \gamma ,p, r ) >0$ such that
\begin{align}
\sup_{0<t <T} \{ \| v(t) - u (t) \|_{L^{3,r}} \} \leq & K_{\dagger \dagger \dagger} \| a - b \|_{L^{3, r }}, \label{eq311}\\
\sup_{0 < t < T} \{  \| v (t) \|_{L^{3,r}} \} \leq & K_{\dagger \dagger \dagger} \| a \|_{L^{3, r}},\label{eq312}\\
\sup_{0 < t <T} \{ \| u (t) \|_{L^{3,r}} \} \leq &  K_{\dagger \dagger \dagger} \|  b \|_{L^{3,r}} \label{eq313}.
\end{align}
\end{lemma}
\begin{proof}[Proof of Lemma \ref{lem31}]
Let us derive \eqref{eq32}. Using \eqref{eq211}, \eqref{eq23}, \eqref{eq24}, \eqref{e2}, and \eqref{e8}, we check that for $0 < t < T$
\begin{multline*}
\| v(t) - u (t) \|_{L^{3,\infty}} \leq C \sup_{\varphi \in C_{0 , \sigma }^\infty, { \ } \| \varphi \|_{L^{\frac{3}{2},1}}\leq 1 }|\dual{v(t) - u (t), \varphi }|\\
\leq C \sup_{ \| \varphi \|_{L^{\frac{3}{2} , 1}} \leq 1  }| \dual{ a - b , \mathrm{e}^{- t L^*} \varphi }| \\
+ C \sup_{ \| \varphi \|_{L^{\frac{3}{2},1}}\leq 1 } \bigg| \int_0^t \langle ( v(s) - u(s) )\otimes v(s) + u(s) \otimes (v (s) - u (s)) \\
+ (v(s) - u (s) ) \otimes w + w \otimes (v(s) - u(s)) , \nabla \mathrm{e}^{-(t-s) L^*} \varphi \rangle  d s\bigg|\\
\leq C \| a - b \|_{L^{3,\infty}} + C \delta \sup_{ \| \varphi \|_{L^{\frac{3}{2},1}}\leq 1 } \int_0^t \| v(s) - u (s) \|_{L^{3,\infty}} \| \nabla \mathrm{e}^{-(t-s)L^*} \varphi \|_{L^{3,1}} d s\\
\leq C_1 \| a - b \|_{L^{3,\infty}} + C_1 \delta \sup_{0 < s <t} \{ \| v(s) - u (s) \|_{L^{3,\infty}} \}.
\end{multline*}
Note that {\color{red}{$C_1 = C_1 (\gamma)>0$}} does not depend on both $t$ and $T$. Therefore we see {\color{red}{that}}
\begin{equation}\label{EQ90}
\sup_{0 < t <T} \{ \| v(t) - u (t) \|_{L^{3,\infty}} \} \leq 2 C_1 \| a - b \|_{L^{3,\infty}}
\end{equation}
if $C_1 \delta \leq 1/2$. This gives \eqref{eq32}. Similarly, we have \eqref{eq33} and \eqref{eq34}. Remark that $\delta_\dagger := 1/(2C_1)$ and $K_\dagger := 2 C_1$.
\end{proof}

\begin{proof}[Proof of Lemma \ref{lem32}]
We first show \eqref{eq35}. Fix $\psi \in L^{3/2 , 1}_\sigma ( \Omega )$ and $\varepsilon >0$. Since $C_{0 , \sigma}^\infty ( \Omega )$ is dense in $L^{ 3/2 , 1 }_\sigma (\Omega)$, there is $\psi_0 \in C_{0,\sigma}^\infty ( \Omega )$ such that
\begin{equation}\label{eq3132}
\| \psi - \psi_0 \|_{L^{\frac{3}{2}, 1}} < \min \left\{ \frac{ \varepsilon }{ 3 \delta } , \frac{\varepsilon }{ 3 \| a \|_{L^{3,\infty}}} \right\} .
\end{equation}
Here $\delta$ is defined by \eqref{eq31}. By \eqref{eq23} and \eqref{eq3132}, we see that for $t>0$
\begin{align}
|  \dual{ v (t) , \psi } - \dual{ a , \psi } | \leq  & | \dual{ v (t), \psi - \psi_0 } | + |\dual{v (t) - a , \psi_0 }| + | \dual{a , \psi_0 - \psi} | \notag \\
\leq & ( \| v (t) \|_{L^{3,\infty}} + \| a \|_{L^{3,\infty}} ) \| \psi - \psi_0 \|_{L^{\frac{3}{2}, 1}} + |\dual{v (t) - a  , \psi_0 }| \notag \\
< & \frac{2}{3} \varepsilon + |\dual{v (t) - a, \psi_0 }|.\label{eq314}
\end{align}
Using \eqref{eq23}, \eqref{eq24}, and \eqref{e4}, we check that for $t>0$
\begin{multline*}
| \dual{v (t) - a , \psi_0 } | \\
\leq | \dual{ a ,( \mathrm{e}^{- t L^*} - 1) \psi_0 } | + \bigg| \int_0^t \dual{v \otimes v + v \otimes w + w \otimes v , \nabla \mathrm{e}^{-(t-s) L^*} \psi_0}  d s \bigg| \\
\leq | \dual{a, (\mathrm{e}^{-t L^*} - 1) \psi_0 }|  + C \delta^2 \int_0^t \| \nabla \mathrm{e}^{- (t - s) L^*} \psi_0 \|_{L^{3, 1}}  d s \\
\leq \| a \|_{L^{3,\infty}} \| (\mathrm{e}^{- t L^*} - 1) \psi_0 \|_{L^{\frac{3}{2} ,1}} +C \delta^2 \int_0^t \frac{\| \psi_0 \|_{L^{2,1}}}{ (t - s)^\frac{3}{4}}  d s\\
\leq \| a \|_{L^{3,\infty}} \| (\mathrm{e}^{- t L^*} - 1) \psi_0 \|_{L^{\frac{3}{2} ,1}}  + C t^\frac{1}{4} \delta^2 \| \psi_0 \|_{L^{2,1}} .
\end{multline*}
Since $\mathrm{e}^{- t L^*}$ is a $C_0$-semigroup on $L^{ 3/2 , 1}_\sigma (\Omega)$, there is $T_0 >0$ such that for $0 < t < T_0$
\begin{equation}\label{eq315}
| \dual{v (t) - a , \psi_0 } | < \frac{\varepsilon}{3} .
\end{equation}
Combining \eqref{eq314} and \eqref{eq315} implies that for $0 < t < T_0$
\begin{equation*}
| \dual{v(t) - a, \psi } | < \varepsilon .
\end{equation*}
Since $\varepsilon$ is arbitrary, we see \eqref{eq35}. Similarly, we have \eqref{eq36}.
\end{proof}

\begin{proof}[Proof of Lemma \ref{lem33}]
Let $C_1$ be the constant appearing in  \eqref{EQ90}. Assume that $\delta \leq 1/(2C_1)$. We first derive \eqref{eq38}. Let $t \in (0,T)$. Using \eqref{eq211}, we check that
\begin{multline*}
t^\frac{p -3}{2 p} \| v(t) - u (t) \|_{L^{ p , \infty }} \leq C t^\frac{p -3}{2 p} \sup_{\varphi \in C_{0 , \sigma }^\infty, { \ } \| \varphi \|_{L^{p',1}}\leq 1 }| \dual{ a - b , \mathrm{e}^{- t L^*} \varphi }| \\
+ C t^\frac{p -3}{2 p} \sup_{\varphi \in C_{0 , \sigma }^\infty, { \ } \| \varphi \|_{L^{p',1}}\leq 1 } \bigg| \int_0^{\frac{t}{2}} \langle ( v(s) - u(s) )\otimes v(s) + u(s) \otimes (v (s) - u (s)) \\
+ (v(s) - u (s) ) \otimes w + w \otimes (v(s) - u(s)) , \nabla \mathrm{e}^{-(t-s) L^*} \varphi \rangle  d s\bigg|\\
+ C t^\frac{p -3}{2 p} \sup_{\varphi \in C_{0 , \sigma }^\infty, { \ } \| \varphi \|_{L^{p',1}} \leq 1 } \bigg| \int_{\frac{t}{2}}^t \langle ( v(s) - u(s) )\otimes v(s) + u(s) \otimes (v (s) - u (s)) \\
+ (v(s) - u (s) ) \otimes w + w \otimes (v(s) - u(s)) , \nabla \mathrm{e}^{-(t-s) L^*} \varphi \rangle  d s\bigg|\\
=: J_1(t) + J_2(t) + J_3 (t).
\end{multline*}
By \eqref{eq23} and \eqref{e2}, we have
\begin{equation*}
J_1 (t) \leq C \| a - b \|_{L^{3,\infty}} .
\end{equation*}
Note that $p' < 3/2$. Using \eqref{eq23}, \eqref{eq24}, and \eqref{e4}, we see that
\begin{multline*}
J_2 (t) \leq C \delta t^\frac{p -3}{2 p} \int_0^{\frac{t}{2} } \frac{ (s^\frac{p - 3}{2 p} \| v( s )- u (s) \|_{L^{p,\infty}} )}{(t -s )s^\frac{p -3}{2 p}} d s\\
\leq C \delta \int_0^{\frac{t}{2} } \frac{ ( s^\frac{p - 3}{2 p} \| v( s )- u (s) \|_{L^{p,\infty}} )}{(t -s )^{\frac{p + 3}{2 p}}s^\frac{p -3}{2 p}} d s\\
\leq C \delta \sup_{0< s < t}\{ s^\frac{p -3}{2 p} \| v (s) - u (s) \|_{L^{p,\infty}} \}.
\end{multline*}
Note that $p' < 3p/(2 p - 3) < 3$. From \eqref{eq23}, \eqref{eq24}, and \eqref{e8}, we observe that
\begin{multline*}
J_3 (t) \leq C \delta t^\frac{p -3}{2 p} \sup_{ \| \varphi \|_{L^{p',1}} \leq 1 } \int_{\frac{t}{2}}^t \frac{(s^\frac{p - 3}{2 p} \| v ( s ) - u (s) \|_{L^{p , \infty }} )}{s^\frac{p - 3}{2 p}} \| \nabla \mathrm{e}^{- (t -s ) L^*} \varphi \|_{L^{\frac{3 p}{2 p -3},1}} d s \\
\leq C \delta \sup_{ \| \varphi \|_{L^{p',1}} \leq 1 } \int_{\frac{t}{2}}^t (s^\frac{p - 3}{2 p} \| v ( s ) - u (s) \|_{L^{p , \infty }} ) \| \nabla \mathrm{e}^{- (t -s ) L^*} \varphi \|_{L^{\frac{3 p}{2 p -3},1}} d s \\
\leq C \delta \sup_{0<s <t}\{ s^\frac{p - 3}{2 p} \| v ( s ) - u (s) \|_{L^{p , \infty }} \} .
\end{multline*}
Consequently, there is $C_2 = C_2 ( \gamma,p) >0$ such that for $0 < t < T$
\begin{equation*}
t^\frac{p -3}{2 p} \| v (t) - u (t) \|_{L^{p,\infty}} \leq C_2 \| a - b \|_{L^{3,\infty}} + C_2 \delta \sup_{0<s<t}\{ s^\frac{p-3}{2 p}\| v (s) - u (s) \|_{L^{p,\infty}} \}.
\end{equation*}
Therefore we find that
\begin{equation*}
\sup_{0< t <T}\{ t^\frac{p -3}{2 p} \| v (t) - u (t) \|_{L^{p,\infty}} \} \leq 2 C_2 \| a - b \|_{L^{3,\infty}}
\end{equation*}
if $C_2 \delta \leq 1/2$. This gives \eqref{eq38}. Similarly, we have \eqref{eq39} and \eqref{eq310}. Remark that $\delta_{ \dagger \dagger }:= \min \{ 1/(2C_1) , 1/(2 C_2) \}$ and $K_{\dagger \dagger} := 2 C_1 + 2 C_2$. By the same argument to derive \eqref{eq32}, we have \eqref{eq37}.

Next we show $(\mathrm{i})$ and $(\mathrm{ii})$. We now assume that $\delta \leq \delta_\dagger$ and that $a,b \in L^{3,r}_\sigma (\Omega)$ for some $1 < r < \infty$. Let $r' \in (1, \infty )$ such that $1/r + 1 /{r'}=1$.

We shall show that $v \in BC ((0, T) ; L^{3,r}_\sigma ( \Omega ))$. Let $t \in ( 0 , T )$. Using \eqref{eq211}, \eqref{eq23}, \eqref{eq24}, \eqref{e2}, \eqref{eq33}, and \eqref{pp4}, we check that
\begin{multline*}
\| v ( t ) \|_{L^{3 , r}} \leq C \sup_{\varphi \in C_{ 0 , \sigma }^\infty { \ } \| \varphi \|_{L^{\frac{3}{2} , r'}} \leq 1 } | \dual{a , \mathrm{e}^{- t L^* } \varphi } | \\
+C \delta \sup_{\varphi \in C_{0 , \sigma}^\infty { \ } \| \varphi \|_{L^{\frac{3}{2} , r'}} \leq 1 } \int_0^t \frac{(s^\frac{p - 3}{2 p} \| v (s) \|_{L^{p,\infty}})}{s^\frac{p - 3}{2 p}} 
\| \nabla \mathrm{e}^{- (t -s) L^*} \varphi \|_{L^{\frac{3 p}{2 p - 3} , 1}} d s \\
\leq C \| a \|_{L^{3,r}} + C \delta K_{\dagger \dagger} \| a \|_{L^{3,\infty}} \int_0^t \frac{1}{s^\frac{p - 3}{2 p}(t-s)^\frac{p + 3}{2 p}} d s \\
\leq C \| a \|_{L^{3,r}} + C \delta K_{\dagger \dagger} \| a \|_{L^{3,\infty}} .
\end{multline*}
This shows that
\begin{equation}\label{eq318}
\sup_{0 < t < T} \{ \| v ( t ) \|_{L^{3,r}} \} \leq C \| a \|_{L^{3,r}} + C \delta K_\dagger \| a \|_{L^{3,\infty}} < + \infty .
\end{equation}
Fix $\varepsilon >0$ such that $\varepsilon < T$. Let $t_1 , t_2 >0$ such that $\varepsilon \leq t_1 \leq t_2 < T$. By \eqref{eq211}, we see that
\begin{multline*}
\| v (t_2) - v( t_1 ) \|_{L^{3,r}} \leq C \sup_{\varphi \in C_{0,\sigma}^\infty, { \ }\| \varphi \|_{L^{\frac{3}{2}, r'}} \leq 1} |\dual{a , ( \mathrm{e}^{- (t_2 -t_1) L^*} - 1) \mathrm{e}^{- t_1 L^*} \varphi } |  \\
+ C \sup_{\| \varphi \|_{L^{\frac{3}{2}, r'}} \leq 1} \bigg| \int_0^{t_1} \dual{v \otimes v + v \otimes w + w \otimes v , \nabla ( \mathrm{e}^{-(t_2 -s) L^*} - \mathrm{e}^{- (t_1 - s) L^*} )  \varphi } d s \bigg|\\
+ C \sup_{\| \varphi \|_{L^{\frac{3}{2}, r'}} \leq 1} \bigg| \int_{t_1}^{t_2} \dual{ v(s) \otimes v (s) + v (s) \otimes w + w \otimes v (s) , \nabla \mathrm{e}^{-(t_2 -s) L^*} \varphi } d s \bigg|\\
=: J_4 (t_1,t_2) + J_5 (t_1,t_2) + J_6 (t_1,t_2).
\end{multline*}
Using \eqref{eq23}, \eqref{kk4}, and \eqref{kk2}, we find that
\begin{multline*}
J_4 (t_1,t_2) \leq C \| a \|_{L^{3,r}}(t_2-t_1)^{\frac{1}{4}} \mathrm{e}^{(t_2 - t_1)} \sup_{\| \varphi \|_{L^{\frac{3}{2} , r'} } \leq 1} \| (L^* + 1)^{\frac{1}{4}} \mathrm{e}^{-  t_1 L} \varphi \|_{L^{\frac{3}{2}, r'}} \\
\leq C (T) (t_2 -t_1)^{\frac{1}{4}} \mathrm{e}^{t_1} t_1^{- \frac{1}{4} } \leq C (\varepsilon , T) (t_2 - t_1 )^{\frac{1}{4}}.
\end{multline*}
Applying \eqref{eq23}, \eqref{eq24}, \eqref{eq33}, and \eqref{6b}, we check that
\begin{multline}\label{eqj5}
J_5 (t_1,t_2) \leq C \sup_{\| \varphi \|_{L^{\frac{3}{2}, r'}} \leq 1} \bigg\{ \int_0^{t_1} \frac{ ( \| v (s) \|_{L^{3,\infty}} + \| w \|_{L^{3,\infty}} ) ( s^\frac{p - 3}{2 p} \| v (s) \|_{L^{ p,\infty}} )}{s^\frac{p - 3}{2 p}}\\
\cdot \| \nabla \{ ( \mathrm{e}^{- (t_2 -t_1) L^*} - 1 ) \mathrm{e}^{- (t_1-s) L^*} \varphi \} \|_{L^{\frac{3 p }{2 p - 3},1}} d s \bigg\}\\
\leq C \delta K_{\dagger \dagger} \| a \|_{L^{3 , \infty}} \int_0^{ t_1} \frac{1}{s^{\frac{p -3 }{2 p}}} \cdot (t_2-t_1)^{ \frac{1}{4 p} } \mathrm{e}^{\frac{t_1 - s}{2} + t_2 - t_1} (t_1 - s)^{- \left( \frac{1}{2} + \frac{1}{4 p} \right) - \frac{1}{4 p} -\frac{3}{2 p} } d s\\
\leq C ( T ) t_1^{- \frac{1}{2 p} } (t_2 - t_1)^{\frac{1}{4 p} } \leq C (\varepsilon , T) (t_2 - t_1 )^{ \frac{1}{4 p}} .
\end{multline}
By \eqref{eq23}, \eqref{eq24}, \eqref{eq33}, and \eqref{pp4}, we observe that
\begin{multline*}
J_6 (t_1,t_2) \leq C \sup_{ \| \varphi \|_{L^{\frac{3}{2}, r'}} \leq 1} \bigg\{ \int_{t_1}^{t_2} \frac{ ( \| v (s) \|_{L^{3,\infty}} + \| w \|_{L^{3,\infty}} ) ( s^\frac{p -3}{2 p} \| v (s) \|_{L^{p,\infty}} ) }{s^\frac{p -3}{2 p}}\\
\cdot \| \nabla \mathrm{e}^{-  (t_2-s) L^*} \varphi \|_{L^{\frac{3 p }{2 p - 3},1}} d s \bigg\} \\
\leq C \delta K_{\dagger \dagger} \| a \|_{L^{3,\infty}} t_1^{- \frac{p - 3}{2 p} } \int_{t_1}^{t_2} \frac{1}{ (t_2 - s)^\frac{p + 3}{2 p} } d s \leq C ( \varepsilon )( t_2 - t_1)^\frac{ p - 3}{2 p } . 
\end{multline*}
Therefore we conclude that
\begin{equation*}
\| v (t_2) - v (t_1) \|_{L^{3 , r} } \to 0 \text{ as } t_2 \to t_1.
\end{equation*}
This implies that
\begin{equation*}
v \in C([\varepsilon , T ) ; L^{3, r }_\sigma ( \Omega ) ).
\end{equation*}
Since $\varepsilon$ is arbitrary, it follows from \eqref{eq318} that
\begin{equation*}
v \in BC( ( 0, T ) ; L^{3, r }_\sigma ( \Omega ) ).
\end{equation*}
Similarly, we check that $u \in BC( ( 0, T ) ; L^{3, r }_\sigma ( \Omega ) )$. Therefore we see $(\mathrm{i})$.

Finally, we show $(\mathrm{ii})$. By the same argument to derive \eqref{eq318}, we see that for $0 < t <T$
\begin{multline*}
\| v ( t ) - u (t) \|_{L^{3 , r}} \leq C \sup_{\varphi \in C_{0,\sigma}^\infty , { \ } \| \varphi \|_{L^{\frac{3}{2} , r'}} \leq 1 } | \dual{a - b , \mathrm{e}^{- t L^*} \varphi } | \\
+C \delta \sup_{\varphi \in C_{0,\sigma}^\infty , { \ } \| \varphi \|_{L^{\frac{3}{2} , r'}} \leq 1 } \int_0^t \frac{ (s^\frac{p - 3}{2 p} \| v (s) - u (s) \|_{L^{p,\infty}})}{s^\frac{p - 3}{2 p} (t - s)^{\frac{p + 3}{2 p}}} d s \\
\leq C \| a - b \|_{L^{3,r}} + C \delta K_{\dagger \dagger} \| a - b  \|_{L^{3,\infty}} .
\end{multline*}
Since $\| \cdot \|_{L^{3,\infty}} \leq C \| \cdot \|_{L^{3,r}}$, we have
\begin{equation*}
\| v ( t ) - u (t) \|_{L^{3 , r}} \leq C \| a - b \|_{L^{3,r}} + C \delta K_{\dagger \dagger} \| a - b  \|_{L^{3,r}} .
\end{equation*}
Therefore we see that there is $C_3 = C_3 (\gamma , p , r) >0$ such that
\begin{equation*}
\sup_{0 < t <T} \{ \| v ( t ) - u (t) \|_{L^{3 , r}} \} \leq C_3 \| a - b \|_{L^{3,r}},
\end{equation*}
which is \eqref{eq311}. Similarly, we have \eqref{eq312} and \eqref{eq313}. Remark that $K_{\dagger \dagger \dagger} := 2 C_3$. Therefore the lemma follows.
\end{proof}

{\color{red}{[[}}
Let us show the existence of a unique global-in-time generalized mild solution of the system \eqref{eq15}. We apply a contraction mapping theory to obtain a unique global-in-time generalized mild solution. However, the argument for the application of the contraction mapping theory is not simple, so we give a detailed proof of Theorem \ref{thm19}.

\begin{proof}[Proof of Theorem \ref{thm19}]
Let $\delta >0$. Define
\begin{equation*}
X_p^\delta := \{ f \in BC ((0 , \infty ); L^{3,\infty}_\sigma ( \Omega ) ) \cap C( (0, \infty ); L^{p,\infty}_\sigma ( \Omega ) ) ; { \ } \| f \|_{X_p} \leq \delta \},
\end{equation*}
where $\displaystyle{ \| f \|_{X_p} = \sup_{t>0}\{ \| f (t) \|_{L^{3,\infty }}\} + \sup_{t>0}\{ t^\frac{p -3}{2 p} \| f (t) \|_{L^{p,\infty}} \}}$.\\
For each $u \in X_p^\delta$, $a \in L^{3,\infty}_\sigma ( \Omega)$, $\psi \in L^{3/2 , 1}_\sigma ( \Omega )$, and $t>0$, we set
\begin{multline*}
\mathcal{M}(u , a ,\psi , t ) \\:= \dual{ a , \mathrm{e}^{-t L^*} \psi } + \int_0^t \dual{u(s) \otimes u(s) + u(s) \otimes w + w \otimes u(s) , \nabla \mathrm{e}^{-(t-s) L^*} \psi } d s.
\end{multline*}
By \eqref{eq23}, \eqref{eq24}, \eqref{e2}, and \eqref{e8}, we check that that for all $u \in X_p^\delta$, $a \in L^{3,\infty}_\sigma ( \Omega)$, $\psi \in L^{3/2 , 1}_\sigma ( \Omega )$, and $t>0$,
\begin{multline}\label{eq0421}
| \mathcal{M}(u, a ,\psi , t)  | \leq \| a \|_{L^{3,\infty}}\| \mathrm{e}^{- t L^*} \psi \|_{L^{\frac{3}{2},1}} \\
+ C \int_0^t ( \| u(s) \|_{L^{3,\infty}} + \| w \|_{L^{3,\infty}}) \| u (s) \|_{L^{3,\infty}} \| \nabla \mathrm{e}^{-(t-s) L^*} \psi \|_{L^{3,1}} d s\\
\leq C \| a \|_{L^{3,\infty}} \| \psi \|_{L^{\frac{3}{2},1}} + C (\delta + \| w \|_{L^{3,\infty} }) \delta \| \psi \|_{L^{ \frac{3}{2},1}}.
\end{multline}
Note that the above constant $C >0$ does not depend on $t$.\\
Let us attack Theorem \ref{thm19}. Fix $u \in X_p^\delta$ and $a \in L^{3,\infty}_\sigma ( \Omega)$. To prove Theorem \ref{thm19}, we pass through four steps. In the first step, we show the existence of a $L^{3,\infty}_\sigma$-valued function $v$ such that $\dual{v(t), \psi} = \mathcal{M}(u , a ,\psi , t )$ for $\psi \in L^{3/2 , 1}_\sigma ( \Omega )$ and $t>0$. In the second step, we prove that the function $v$ belongs to $X_p$. In the third step, we show that the map $\mathcal{M}$ is a contraction mapping such that $\mathcal{M} : X_p^\delta \to X_p^\delta$ if $\| a \|_{L^{3,\infty}}$, $\| w \|_{L^{3,\infty}}$, and $\delta$ are sufficiently small. In the final step, we apply a contraction mapping theory to show the existence of a unique global-in-time generalized mild solution of the system \eqref{eq15} when both $\| a \|_{L^{3,\infty}}$ and $\| w \|_{L^{3,\infty}}$ are sufficiently small.\\

\noindent Step 1: \lq\lq Existence of a function $v$ such that $\dual{v(t), \psi} = \mathcal{M}(u , a ,\psi , t )$."

Let $u \in X_p^\delta$ and $a \in L^{3,\infty}_\sigma ( \Omega)$. We fix $u$ and $a$. We now prove that there exists a $L^{3,\infty}_\sigma$-valued function $v: (0,\infty) \ni t \to v (t) \in L^{3,\infty}_\sigma (\Omega)$ such that for each $\psi \in L^{3/2 , 1}_\sigma ( \Omega )$ and $t>0$,
\begin{equation*}
\dual{v (t) , \psi } = \dual{ a , \mathrm{e}^{-t L^*} \psi } + \int_0^t \dual{u \otimes u(s) + u(s) \otimes w + w \otimes u(s) , \nabla \mathrm{e}^{-(t-s) L^*} \psi } d s.
\end{equation*}
To this end, we show that $\mathcal{M}$ is a bounded linear functional on $L^{3/2 , 1}_\sigma ( \Omega )$ when $u,a,t$ are fixed. From \eqref{eq0421}, we check that 
\begin{equation*}
| \mathcal{M}(u, a ,\psi , t)  | 
\leq C \| a \|_{L^{3,\infty}} + C (\delta + \| w \|_{L^{3,\infty} }) \delta < + \infty
\end{equation*}
for every $t >0$ and $\psi \in L^{3/2 , 1}_\sigma ( \Omega )$ satisfying $\| \psi \|_{L^{\frac{3}{2}, 1}} \leq 1$. Now we fix $t>0$. It is clear that for every $\psi , \psi_1, \psi_2 \in L^{3/2 , 1}_\sigma ( \Omega )$ and $\lambda \in \mathbb{R}$,
\begin{align*}
\mathcal{M}(u , a , \psi_1 + \psi_2  , t ) & = \mathcal{M}(u , a ,\psi_1 , t ) + \mathcal{M}(u , a ,\psi_2  , t ),\\
\mathcal{M}(u , a , \lambda \psi , t ) & = \lambda \mathcal{M}(u , a ,\psi , t ).
\end{align*}
Therefore we see that $\mathcal{M}$ is a bounded linear functional on $L^{3/2 , 1}_\sigma ( \Omega )$ when $u,a,t$ are fixed. Since the dual space of $L^{3/2 , 1}_\sigma ( \Omega )$ is $L^{3 , \infty}_\sigma ( \Omega )$, we apply the Hahn-Banach theorem to find that there is $v (t) \in L^{3 , \infty}_\sigma ( \Omega )$ such that for every $\psi \in L^{3/2 , 1}_\sigma ( \Omega )$ 
\begin{equation*}
\dual{v (t) , \psi } = \mathcal{M} (u,a, \psi , t).
\end{equation*}
Since the above argument is valid for each $t>0$, we conclude that for each fixed $u \in X_p^\delta$ and $a \in L^{3,\infty}_\sigma (\Omega)$ there exists a $L^{3,\infty}_\sigma$-valued function $v: (0,\infty) \ni t \to v (t) \in L^{3,\infty}_\sigma (\Omega)$ such that for each $\psi \in L^{3/2 , 1}_\sigma ( \Omega )$ and $t>0$,
\begin{equation*}
\dual{v (t) , \psi } \\= \dual{ a , \mathrm{e}^{-t L^*} \psi } + \int_0^t \dual{u \otimes u(s) + u(s) \otimes w + w \otimes u(s) , \nabla \mathrm{e}^{-(t-s) L^*} \psi } d s.
\end{equation*}

\noindent Step 2: \lq\lq Properties of $v$ such that $\dual{v(t), \psi} = \mathcal{M}(u , a ,\psi , t )$."

Let $u \in X_p^\delta$ and $a \in L^{3,\infty}_\sigma ( \Omega)$. We fix $u$ and $a$. We assume that there exists a $L^{3,\infty}_\sigma$-valued function $v: (0,\infty) \ni t \to v (t) \in L^{3,\infty}_\sigma (\Omega)$ such that for each $\psi \in L^{3/2 , 1}_\sigma ( \Omega )$ and $t>0$,
\begin{equation*}
\dual{v (t) , \psi }= \dual{ a , \mathrm{e}^{-t L^*} \psi } + \int_0^t \dual{u(s) \otimes u(s) + u(s) \otimes w + w \otimes u(s) , \nabla \mathrm{e}^{-(t-s) L^*} \psi } d s.
\end{equation*}
We now prove that $v \in X_p$. We first estimate $\| v \|_{X_p}$. Using \eqref{eq211}, \eqref{eq23}, \eqref{eq24}, \eqref{e2}, and \eqref{e8}, we observe that for $t>0$
\begin{multline*}
\| v (t) \|_{L^{3,\infty}} \leq C \sup_{\varphi \in C_{0,\sigma}^\infty, { \ }\| \varphi \|_{L^{\frac{3}{2}, 1}} \leq 1} |\dual{ a , \mathrm{e}^{-t L^*} \varphi } | \\
+ C \sup_{\varphi \in C_{0,\sigma}^\infty, { \ }\| \varphi \|_{L^{\frac{3}{2}, 1}} \leq 1} \bigg| \int_0^t \dual{u(s) \otimes u(s) + u(s) \otimes w + w \otimes u(s) , \nabla \mathrm{e}^{-(t-s) L^*} \varphi } d s \bigg|\\
\leq C \| a \|_{L^{3 , \infty}} \\
+ C \sup_{\| \varphi \|_{L^{\frac{3}{2}, 1}} \leq 1} \int_0^t (\| u ( s) \|_{L^{3,\infty} } + \| w \|_{L^{3,\infty}} ) \| u ( s ) \|_{L^{3,\infty}} \| \nabla \mathrm{e}^{- (t-s) L^*} \varphi \|_{L^{3,1}} d s\\
\leq  C ( \| a \|_{L^{3,\infty}}  +  \| w \|_{L^{3,\infty}} \delta  + \delta^2 ) < + \infty  .
\end{multline*}
Applying \eqref{eq211}, \eqref{eq23}, \eqref{eq24}, \eqref{e2}, \eqref{e4}, and \eqref{e8}, we check that for $t>0$
\begin{multline*}
t^\frac{p- 3}{2 p} \| v (t) \|_{L^{p,\infty}} \leq C t^\frac{p - 3}{2 p} \sup_{\varphi \in C_{0,\sigma}^\infty, { \ }\| \varphi \|_{L^{p', 1}} \leq 1} |\dual{ a , \mathrm{e}^{-t L^*} \varphi } | \\
+ C t^\frac{p-3}{2 p} \sup_{\varphi \in C_{0,\sigma}^\infty, { \ }\| \varphi \|_{L^{p', 1}} \leq 1} \bigg| \int_0^t \dual{u \otimes u + u \otimes w + w \otimes u , \nabla \mathrm{e}^{-(t-s) L^*} \varphi } d s \bigg|\\
\leq C \| a \|_{L^{3 , \infty}} + C t^\frac{p - 3}{2 p} \int_0^{\frac{t}{2}} \frac{(\| u ( s) \|_{L^{3,\infty} } + \| w \|_{L^{3,\infty}} ) (s^\frac{p - 3}{2 p} \| u ( s ) \|_{L^{p,\infty}} )}{s^\frac{p - 3}{2 p} (t - s) } d s
\\ + C \delta t^\frac{p - 3}{2 p} \sup_{\| \varphi \|_{L^{p', 1}} \leq 1} \int_{\frac{t}{2}}^t \frac{(\| u ( s) \|_{L^{3,\infty} } + \| w \|_{L^{3,\infty}} ) }{s^\frac{p - 3}{2 p} } \| \nabla \mathrm{e}^{- (t-s)L^*} \varphi \|_{L^{\frac{3 p}{2 p - 3},1}} d s\\
\leq  C ( \| a \|_{L^{3,\infty}}  +  \| w \|_{L^{3,\infty}} \delta  + \delta^2 ) < + \infty .
\end{multline*}
Therefore, we conclude that there is $C_\dagger = C_\dagger (\gamma, p) >0$ such that
\begin{equation}\label{eqm3}
\| v \|_{X_p} \leq C_\dagger ( \| a \|_{L^{3,\infty}}  +  \| w \|_{L^{3,\infty}} \delta  + \delta^2 ) < + \infty .
\end{equation}
Next we show that $v \in C ((0, \infty ) ; L^{3,\infty}_\sigma ( \Omega )) \cap C ((0, \infty ) ; L^{p,\infty}_\sigma ( \Omega )) $. Let $\varepsilon , T > 0$ such that $\varepsilon < T$. Fix $\varepsilon$ and $T $. Let $t_1 , t_2 >0$ such that $\varepsilon \leq t_1 \leq t_2 < T$. By \eqref{eq211}, we see that
\begin{multline*}
\| v (t_2) - v( t_1 ) \|_{L^{3,\infty}} \leq C \sup_{\varphi \in C_{0,\sigma}^\infty, { \ }\| \varphi \|_{L^{\frac{3}{2}, 1}} \leq 1} |\dual{a, (\mathrm{e}^{- (t_2 -t_1 ) L^*} - 1) \mathrm{e}^{- t_1 L^*} \varphi } | \\
+ C \sup_{\| \varphi \|_{L^{\frac{3}{2}, 1}} \leq 1} \bigg| \int_0^{t_1} \dual{u \otimes u + u \otimes w + w \otimes u, \nabla \{ ( \mathrm{e}^{-(t_2 -s) L^*} - \mathrm{e}^{- (t_1 - s) L^*} ) \varphi \}} d s \bigg|\\
+ C \sup_{\| \varphi \|_{L^{\frac{3}{2}, 1}} \leq 1} \bigg| \int_{t_1}^{t_2} \dual{u(s) \otimes u(s) + u(s) \otimes w + w \otimes u(s) , \nabla \mathrm{e}^{-(t_2 -s) L^*} \varphi } d s \bigg|\\
=: I_1 (t_1,t_2) + I_2 (t_1,t_2) + I_3 (t_1,t_2).
\end{multline*}
From \eqref{eq23}, \eqref{kk4}, and \eqref{kk2}, we observe that
\begin{align*}
I_1 (t_1,t_2) \leq & C\|  a \|_{L^{3,\infty}} \sup_{\| \varphi \|_{L^{ \frac{3}{2} , 1 }} \leq 1} \| ( \mathrm{e}^{- (t_2 - t_1) L^*} - 1)  \mathrm{e}^{-  t_1 L^*} \varphi \|_{L^{\frac{3}{2}, 1}} \\
\leq & C (t_2- t_1)^{\frac{1}{4}} \mathrm{e}^{t_2 - t_1 } \sup_{\| \varphi \|_{ L^{\frac{3}{2},1}} \leq 1} \| (L^* + 1)^\frac{1}{4} \mathrm{e}^{- t_1 L} \varphi \|_{L^{ \frac{3}{2}, 1}}\\
\leq & C \varepsilon^{-\frac{1}{4}} (t_2 - t_1)^\frac{1}{4} \mathrm{e}^{2 T} .
\end{align*}
By the same argument to derive \eqref{eqj5}, we check that
\begin{multline*}
I_2 (t_1,t_2) \leq C \sup_{\| \varphi \|_{L^{\frac{3}{2}, 1}} \leq 1} \bigg\{ \int_0^{t_1} \frac{ ( \| u (s) \|_{L^{3,\infty}} + \| w \|_{L^{3,\infty}} ) ( s^\frac{p - 3}{2 p} \| u (s) \|_{L^{ p,\infty}} )}{s^\frac{p - 3}{2 p}}\\
\cdot \| \nabla \{ ( \mathrm{e}^{- (t_2 -t_1) L^*} - 1 ) \mathrm{e}^{- (t_1-s)  L^*} \varphi \} \|_{L^{\frac{3 p }{2 p - 3},1}} d s \bigg\}\\
\leq C ( \varepsilon ,T ) (t_2 - t_1)^\frac{1}{4 p} .
\end{multline*}
Using \eqref{eq23}, \eqref{eq24}, and \eqref{e4}, we find that
\begin{multline*}
I_3 \leq C \delta \sup_{ \| \varphi \|_{L^{\frac{3}{2}, 1}} \leq 1} \int_{t_1}^{t_2} \frac{ ( \| u (s) \|_{L^{3,\infty}} + \| w \|_{L^{3,\infty}} ) }{s^\frac{p -3}{2 p}}  \| \nabla \mathrm{e}^{-  (t_2-s) L^*} \varphi \|_{L^{\frac{3 p }{2 p - 3},1}}d s\\
\leq C ( \delta + \| w \|_{L^{3,\infty}}) \delta t_1^{- \frac{p-3}{2 p}} \int_{t_1}^{t_2} \frac{1}{ (t_2 - s)^\frac{p + 3}{2 p} } d s \\
\leq C \varepsilon^{- \frac{p -3}{2 p}} ( t_2 - t_1)^\frac{p - 3}{2 p } . 
\end{multline*}
Therefore we conclude that
\begin{equation*}
\| v (t_2) - v (t_1) \|_{L^{3,\infty}} \to 0 \text{ as } t_2 \to t_1.
\end{equation*}
This implies that
\begin{equation}\label{aeq}
v \in C([\varepsilon , T ) ; L^{3, \infty}_\sigma ( \Omega ) ).
\end{equation}
Next we prove that
\begin{equation}\label{beq}
v \in C([\varepsilon , T ) ; L^{p, \infty}_\sigma ( \Omega ) ).
\end{equation}
By \eqref{eq211}, we have
\begin{multline*}
\| v (t_2) - v( t_1 ) \|_{L^{p,\infty}} \leq C \sup_{\varphi \in C_{0,\sigma}^\infty, { \ }\| \varphi \|_{L^{ p', 1}} \leq 1} |\dual{a, \mathrm{e}^{-\frac{t_1}{2} L^*} (\mathrm{e}^{- (t_2 -t_1) L^*} -1) \mathrm{e}^{- \frac{t_1}{2} L^*} \varphi } | \\
+ C \sup_{\| \varphi \|_{L^{ p' , 1}} \leq 1} \bigg| \int_0^{t_1} \dual{u \otimes u + u \otimes w + w \otimes u, \nabla \{ ( \mathrm{e}^{-(t_2 -s) L^*} - \mathrm{e}^{-(t_1- s)L^*} ) \varphi \} } d s \bigg|\\
+ C \sup_{\| \varphi \|_{L^{p', 1}} \leq 1} \bigg| \int_{t_1}^{t_2} \dual{u(s) \otimes u(s) + u(s) \otimes w + w \otimes u(s) , \nabla \mathrm{e}^{-(t_2 -s) L^*} \varphi } d s \bigg|\\
=: I_4 ( t_1,t_2 ) + I_5 ( t_1,t_2 ) + I_6 ( t_1 ,t_2 ).
\end{multline*}
Applying \eqref{eq23}, \eqref{e2}, \eqref{kk4}, and \eqref{kk2}, we see that
\begin{multline*}
I_4 ( t_1,t_2) \\\leq  C \| a \|_{L^{3 , \infty}} t_1^{-\frac{3}{2} \left( \frac{1}{p'} - \frac{2}{3} \right)} (t_2 - t_1)^{\frac{1}{4}} \mathrm{e}^{(t_2-t_1)} \sup_{\| \varphi \|_{L^{p',1}}\leq 1} \| (L^* + 1)^{\frac{1}{4}} \mathrm{e}^{- t_1 L^*}\varphi \|_{L^{p',1}}\\
\leq C(\varepsilon , T) (t_2 - t_1)^\frac{1}{4} .
\end{multline*}
Set $\beta_1 = 3 \beta /(3 - \beta)$. Using \eqref{eq23}, \eqref{eq24}, and \eqref{6b}, we check that
\begin{multline*}
I_5 (t_1 , t_2 ) \leq C \sup_{\| \varphi \|_{L^{p', 1}} \leq 1} \\
\bigg\{ \int_0^{t_1} \frac{ ( s^\frac{p - 3}{2 p} \| u (s) \|_{L^{ p, \infty}} )^2 }{s^\frac{p - 3}{p}} \| \nabla \{  ( \mathrm{e}^{- (t_2 -t_1) L^*} - 1 ) \mathrm{e}^{-  (t_1-s) L^*} \varphi \} \|_{L^{\frac{p }{p - 2},1}} d s \bigg\}\\
+ C \sup_{\| \varphi \|_{L^{p', 1}} \leq 1} \bigg\{ \int_0^{t_1} \frac{ \| w \|_{L^{\beta_1 ,\infty}} ( s^\frac{p - 3}{2 p} \| u (s) \|_{L^{ p,\infty}} )}{s^\frac{p - 3}{2 p}}\\
\cdot \| \nabla \{ ( \mathrm{e}^{- (t_2 -t_1) L^*} - 1 ) \mathrm{e}^{-  (t_1-s)  L^*} \varphi \} \|_{L^{\frac{p \beta_1 }{p \beta_1 - p - \beta_1},1}} d s \bigg\}\\
\leq C \int_0^{t_1} \frac{1}{s^{\frac{p -3}{p} }} (t_2 - t_1 )^{\frac{1}{4 p}} \mathrm{e}^{\frac{t_1 -s}{2} + t_2 -t_1} (t_1 - s)^{- \left( \frac{1}{2} + \frac{1}{4 p} \right) - \frac{1}{4 p} - \frac{3}{2 p}} d s\\
+ C \int_0^{t_1} \frac{1}{s^{\frac{p -3}{2 p} }} (t_2 - t_1 )^{\frac{1}{4 \beta_1}} \mathrm{e}^{\frac{t_1 -s}{2} + t_2 -t_1} (t_1 - s)^{- \left( \frac{1}{2} + \frac{1}{4 \beta_1} \right) - \frac{1}{4 \beta_1} - \frac{3}{2 \beta_1}} d s \\
\leq C ( \varepsilon ,T) \{ (t_2 -t_1)^{\frac{1}{4 p} } + (t_2 -t_1)^{\frac{1}{4 \beta_1}} \} .
\end{multline*}
Note that $w \in L^{\beta_1 , \infty} ( \Omega )$ by Assumption \ref{assA}. Note also that $p' \leq p/(p-2) < 3$ and that $p' \leq p\beta_1/(p \beta_1 - p - \beta_1) \leq 3$ since $3/2 \leq \beta/ (\beta - 1)\leq 3 <p$, $1/p+1/p'=1$, and $\beta_1 = 3 \beta /(3 - \beta)$. \\
\noindent From \eqref{eq23}, \eqref{eq24}, and \eqref{e4}, we observe that
\begin{multline*}
I_6 (t_1 , t_2) \leq C \sup_{\| \varphi \|_{L^{p',1}} \leq 1} \int_{t_1}^{t_2} \frac{(s^\frac{p-3}{2 p} \| u (s) \|_{L^{p,\infty}})^2}{s^\frac{p-3}{p}} \| \nabla \mathrm{e}^{- (t_2 - s)L^*} \varphi \|_{L^{\frac{p}{p - 2} , 1}} d s \\
+ C \sup_{\| \varphi \|_{L^{p',1}} \leq 1} \int_{t_1}^{t_2} \frac{\| w \|_{L^{\beta_1 , \infty}}(s^\frac{p-3}{2 p} \| u (s) \|_{L^{p,\infty}})}{s^\frac{p-3}{2 p}} \| \nabla \mathrm{e}^{- (t_2 - s)L^*} \varphi \|_{L^{\frac{p \beta_1}{p \beta_1 - p - \beta_1} , 1}} d s \\
\leq C \delta^2 t_1^{- \frac{p-3}{p}} \int_{t_1}^{t_2} \frac{1}{ (t_2 - s)^\frac{p + 3}{2 p}} d s + C \| w \|_{L^{\beta_1 , \infty}} \delta t_1^{- \frac{p- 3}{2 p}}\int_{t_1}^{t_2}\frac{1}{(t_2 - s)^\frac{\beta_1 + 3}{2 \beta_1}} d s\\
\leq C \varepsilon^{- \frac{p-3}{p}} (t_2 - t_1)^\frac{p -3}{2 p} + C \varepsilon^{- \frac{ p -3}{2 p}} (t_2 - t_1)^\frac{\beta_1 -3}{2 \beta_1} .
\end{multline*}
Therefore we conclude that
\begin{equation*}
\| v (t_2) - v (t_1) \|_{L^{p,\infty}} \to 0 \text{ as }t_2 \to t_1 .
\end{equation*}
This implies that \eqref{beq}. From \eqref{aeq} and \eqref{beq}, we find that
\begin{equation*}
v \in C ([\varepsilon , T) ; L^{3,\infty}_\sigma ( \Omega ) ) \cap C ([\varepsilon , T) ; L^{p,\infty}_\sigma ( \Omega ) ).
\end{equation*}
Since $\varepsilon$ and $T$ are arbitrary, we conclude that
\begin{equation}\label{eqm2}
v \in C( ( 0, \infty ) ; L^{3,\infty}_\sigma ( \Omega ) ) \cap C( ( 0, \infty ) ; L^{p,\infty}_\sigma ( \Omega ) ).
\end{equation}
From \eqref{eqm3} and \eqref{eqm2}, we see that
\begin{equation*}
v \in X_p .
\end{equation*}
Consequently, we find that if there exists a $L^{3,\infty}$-valued function $v$ such that for all $\psi \in L^{3/2 , 1}_\sigma ( \Omega )$ and $t>0$,
\begin{equation*}
\dual{v (t) , \psi } = \mathcal{M}(u, a ,\psi , t )
\end{equation*}
then $v \in X_p$. This implies that we can consider $\mathcal{M}$ as the map $\mathcal{M} : X_p^\delta \ni u \to v \in X_p$ if $a$ is fixed. \\

\noindent Step 3: \lq\lq Contraction mapping $\mathcal{M}$."

We now prove that the map $\mathcal{M}$ is a contraction mapping such that $\mathcal{M} : X_p^\delta \to X_p^\delta$ if $\| a \|_{L^{3,\infty}}$, $\| w \|_{L^{3,\infty}}$, and $\delta$ are sufficiently small. Let $a \in L^{3,\infty}_\sigma (\Omega)$. Now we assume that
\begin{align}
\delta \leq & \min \left\{ \frac{1}{4}, \frac{1}{2 C_\dagger} , \delta_\dagger , \delta_{\dagger \dagger} \right\} ,\label{eq0421a}\\
\| a \|_{L^{3,\infty}} \leq & \min \left\{ \frac{\delta}{4 C_\dagger}, \delta \right\}, \\
\| w \|_{L^{3,\infty}} \leq & \min \left\{ \frac{1}{4 C_\dagger} , \delta \right\} .\label{eq0421c}
\end{align}
Here $C_\dagger$ is the constant appearing in \eqref{eqm3} and, $\delta_\dagger$ and $\delta_{\dagger \dagger}$ are the constants appearing in Lemmas \ref{lem31} and \ref{lem33}, respectively. Remark that $C_\dagger \geq 1$ in general.

We fix $a$. Let $u \in X_p^\delta$. From the argument in Step 1, we see that there exists a $L^{3,\infty}_\sigma$-valued function $v: (0,\infty) \ni t \to v (t) \in L^{3,\infty}_\sigma (\Omega)$ such that for each $\psi \in L^{3/2 , 1}_\sigma ( \Omega )$ and $t>0$,
\begin{equation*}
\dual{v (t) , \psi } = \dual{ a , \mathrm{e}^{-t L^*} \psi } + \int_0^t \dual{u \otimes u(s) + u(s) \otimes w + w \otimes u(s) , \nabla \mathrm{e}^{-(t-s) L^*} \psi } d s.
\end{equation*}
We also see that $v \in X_p$ from the argument in Step 2. By \eqref{eq0421a}-\eqref{eq0421c}, and \eqref{eqm3}, we have
\begin{equation}\label{e13199}
\| v \|_{X_p} \leq \delta .
\end{equation}
Thus, we see that $v \in X_p^\delta$. Therefore we conclude that the map $\mathcal{M}$ is a mapping such that $\mathcal{M} : X_p^\delta \to X_p^\delta$. 

Next we show that the map $\mathcal{M}$ is a contraction mapping such that $\mathcal{M} : X_p^\delta \to X_p^\delta$. Let $u_1,u_2 \in X_p^\delta$. From the previous argument, we see that there are $v_1,v_2 \in X_p^\delta$ such that for each $\psi \in L^{3/2 , 1}_\sigma ( \Omega )$ and $t>0$,
\begin{multline*}
\dual{v_1 (t) , \psi } \\= \dual{ a , \mathrm{e}^{-t L^*} \psi } + \int_0^t \dual{u_1 \otimes u_1(s) + u_1(s) \otimes w + w \otimes u_1(s) , \nabla \mathrm{e}^{-(t-s) L^*} \psi } d s
\end{multline*}
and
\begin{multline*}
\dual{v_2 (t) , \psi } \\= \dual{ a , \mathrm{e}^{-t L^*} \psi } + \int_0^t \dual{u_2 \otimes u_2(s) + u_2(s) \otimes w + w \otimes u_2(s) , \nabla \mathrm{e}^{-(t-s) L^*} \psi } d s.
\end{multline*}
In the same manner to derive \eqref{eqm3} with \eqref{eq0421a} and \eqref{eq0421c}, we see that
\begin{align*}
\| v_1 - v_2 \|_{X_p} & \leq C_\dagger (\| w \|_{L^{3,\infty}} + 2 \delta ) \| u_1 - u_2 \|_{X_p}\\
& \leq \frac{3}{4} \| u_1 - u_2 \|_{X_p}.
\end{align*}
Therefore we conclude that the map $\mathcal{M}$ is a contraction mapping such that $\mathcal{M} : X_p^\delta \to X_p^\delta$.

\noindent Step 4: \lq\lq Existence of $v$ such that $\dual{v(t), \psi} = \mathcal{M}(v , a ,\psi , t )$."

Let $a \in L^{3,\infty}_\sigma (\Omega)$. We assume that \eqref{eq0421a}-\eqref{eq0421c}. We fix $a$. From the argument in Step 3, we see that $\mathcal{M}$ is a contraction mapping. Applying a contraction mapping theory, we obtain a unique function $v \in X_p^\delta$ satisfying 
\begin{multline*}
\dual{v (t) , \psi} = \dual{ a , \mathrm{e}^{-t L^*} \psi} + \int_0^t \dual{v \otimes v(s) + v(s) \otimes w + w \otimes v(s) , \nabla \mathrm{e}^{-(t-s) L^*} \psi } d s
\end{multline*}
for all $\psi \in L^{3/2,1}_\sigma ( \Omega )$ and $t>0$. From Lemma \ref{lem31}, we see that the solution $v$ is unique in $X_3$. We also see that there is $K_0 = K_0 (\gamma , p ) >0$ such that
\begin{equation*}
\| v \|_{X_p} \leq  K_0 \| a \|_{L^{3,\infty}} .
\end{equation*}
Therefore Theorem \ref{thm19} is proved. Note that $K_0:= K_{\dagger \dagger}$, where $K_{\dagger \dagger}$ is the constant appearing in Lemma \ref{lem33}. 
\end{proof}{\color{red}{]]}}

\subsection{$L^{3,\infty}$-Asymptotic Stability and $L^{3,r}$-Asymptotic Stability}\label{subsec32} 
In this subsection we prove Theorem \ref{thm110}. Let $\delta_0$ and $K_0$ be the two constants appearing in Theorem \ref{thm19}. Our first task is to prove 
\begin{lemma}\label{lem34}
Suppose that $\| w \|_{L^{3,\infty}} < \delta_0$. Assume that $b \in L^{3,\infty}_{\sigma} ( \Omega ) $ and that $\| b \|_{L^{3,\infty}} < \delta_0$. Let $u$ be the global-in-time generalized mild solution of the system \eqref{eq15} with the initial datum $v_0 = b$ satisfying $\| u \|_{X_p} \leq K_0 \| b \|_{L^{3,\infty}}$, obtained by Theorem \ref{thm19}. Let $3 p/(p + 3) < \alpha < 3$. Then there is $\delta_1 = \delta_1 ( \gamma , p , \alpha ) >0$ such that if 
\begin{equation*}
b \in L^{\alpha, \infty}_\sigma (\Omega), { \ }\| b \|_{L^{3,\infty}} < \delta_1, \text{ and } \| w \|_{L^{3,\infty } } < \delta_1 ,
\end{equation*}
then the following three assertions hold:\\
\noindent $(\mathrm{i})$ The solution $u$ satisfies
\begin{align}
& \sup_{ t > 0 } \{ t^{\frac{3}{2}(\frac{1}{\alpha} - \frac{1}{p})} \| u (t) \|_{L^{p , \infty } } \} \leq \text{Const.} < + \infty \label{r2},\\
& \sup_{ t > 0 } \{  \| u (t) \|_{L^{\alpha , \infty } } \} \leq \text{Const.} < + \infty \label{r4},\\
& \sup_{t>0}\{ t^{\frac{3}{2} (\frac{1}{\alpha} - \frac{1}{3})} \| u (t) \|_{L^{3,\infty}} \} \leq \text{Const.} < + \infty \label{r5} .
\end{align}
\noindent $(\mathrm{ii})$ Assume in addition that $b \in L^{3,\infty}_{0, \sigma} (\Omega )$ and that $\| \nabla w \|_{L^2} < + \infty$. Then
\begin{equation}
\lim_{t \to 0 +0} \| u ( t ) - b \|_{L^{3,\infty}} = 0. \label{r3}\\
\end{equation}
\noindent $(\mathrm{iii})$ Assume in addition that $b \in L^{3, r }_{\sigma} (\Omega )$ for some $1 < r < \infty$ and that $\| \nabla w \|_{L^2} < + \infty$. Then
\begin{equation}
\lim_{t \to 0 +0} \| u ( t ) - b \|_{L^{3, r }} = 0. \label{r333}\\
\end{equation}

\end{lemma}
\begin{proof}[Proof of Lemma \ref{lem34}]
Let $3 p /(p+3) < \alpha < 3$. Fix $\alpha$. Assume that $b \in L^{\alpha , \infty}_\sigma ( \Omega )$.

We first derive \eqref{r2}. Fix $T>0$. Since $\| u \|_{X_p} \leq K_0 \| b \|_{L^{3 , \infty}}$, we check that for $0 < t <T$
\begin{align*}
t^{\frac{3}{2} \left( \frac{1}{\alpha} - \frac{1}{p} \right)} \| u (t) \|_{L^{p,\infty}} & \leq (t^\frac{p-3}{2 p} \| u (t) \|_{L^{p,\infty}} ) t^{ \frac{3}{2} \left( \frac{1}{\alpha} - \frac{1}{3} \right) }\\
& \leq K_0 \| b \|_{L^{3 ,\infty}} T^{\frac{3}{2} \left( \frac{1}{\alpha} - \frac{1}{3} \right)} < \infty .
\end{align*}
Let $0< t <T$. Lemma \ref{lem29} shows
\begin{multline}\label{eqs1}
\| u (t) \|_{L^{p , \infty} } \leq C \sup_{\varphi \in C_{0,\sigma}^\infty, { \ }\| \varphi \|_{L^{p',1}} \leq 1}|\dual{ b ,\mathrm{e}^{-t L^*} \varphi }| \\
+ C \sup_{\varphi \in C_{0,\sigma}^\infty, { \ }\| \varphi \|_{L^{p',1}} \leq 1}|\dual{ u (s) \otimes u (s) + u(s) \otimes w + w \otimes u(s) , \nabla \mathrm{e}^{- (t-s)L^*}\varphi }|\\
=: D_1 (t) + D_2 (t).
\end{multline}
By \eqref{eq23} and \eqref{e2}, we have
\begin{equation}\label{eqs2}
D_1 (t) \leq C t^{- \frac{3}{2}(\frac{1}{\alpha} - \frac{1}{p})}\| b \|_{L^{\alpha, \infty}} .
\end{equation}
Note that $p' \leq \alpha' \leq 3$. Applying \eqref{eq23}, \eqref{eq24}, and \eqref{e4}, we see that
\begin{multline}\label{eqs3}
D_2 (t) \leq C \int_0^{\frac{t}{2} } \frac{ ( \| u (s) \|_{L^{3,\infty}} + \| w \|_{L^{3, \infty}} ) (s^{\frac{3}{2}(\frac{1}{\alpha} - \frac{1}{p})} \| u ( s ) \|_{L^{p,\infty}} )}{(t -s )s^{ \frac{3}{2} (\frac{1}{\alpha} - \frac{1}{p}) }} d s\\
 + C \sup_{ \| \varphi \|_{L^{p',1} }\leq 1 } \bigg\{ 
\int_{\frac{t}{2}}^t \frac{( \| u (s) \|_{L^{3,\infty}} + \| w \|_{L^{3 , \infty} } )(s^{ \frac{3}{2}(\frac{1}{\alpha} - \frac{1}{p}) }\| u ( s ) \|_{L^{p , \infty }} )}{s^{\frac{3}{2}(\frac{1}{\alpha} - \frac{1}{p})}} \\
\cdot \| \nabla \mathrm{e}^{- (t -s ) L^*} \varphi \|_{L^{\frac{3 p}{2 p -3},1}} d s \bigg\} \\
\leq C t^{- \frac{3}{2}(\frac{1}{\alpha} - \frac{1}{p})} \int_0^{\frac{t}{2} } \frac{ ( \| u(s) \|_{L^{3,\infty}} + \| w \|_{L^{3, \infty}} ) (s^{\frac{3}{2}(\frac{1}{\alpha} - \frac{1}{p})} \| u ( s ) \|_{L^{p,\infty}} )}{(t -s )^{1- \frac{3}{2} (\frac{1}{\alpha} - \frac{1}{p})}s^{ \frac{3}{2} (\frac{1}{\alpha} - \frac{1}{p}) }} d s\\
 + C t^{- \frac{3}{2}(\frac{1}{\alpha} - \frac{1}{p})}\sup_{ \| \varphi \|_{L^{p',1}} \leq 1 } \bigg\{
\int_{\frac{t}{2}}^t ( \| u (s) \|_{L^{3,\infty}} + \| w \|_{L^{3 , \infty} } ) \\
\cdot (s^{ \frac{3}{2}(\frac{1}{\alpha} - \frac{1}{p}) }\| u ( s ) \|_{L^{p , \infty }} ) \| \nabla \mathrm{e}^{- (t -s ) L^*} \varphi \|_{L^{\frac{3 p}{2 p -3},1}} d s \bigg\}.
\end{multline}
Using \eqref{eqs1}-\eqref{eqs3}, and \eqref{e8}, we find that there is $C_1 = C_1 (\gamma, p , \alpha ) > 0$ such that for $0<t <T$
\begin{multline*}
t^{\frac{3}{2}(\frac{1}{\alpha} - \frac{1}{p} )} \| u (t) \|_{L^{p,\infty}} \leq t^{\frac{3}{2}(\frac{1}{\alpha} - \frac{1}{p} )}(D_1 (t) + D_2 (t)) \\
\leq C \| b \|_{L^{\alpha , \infty}} + C (K_0 \| b \|_{L^{3,\infty}} + \| w \|_{L^{3, \infty}}) \int_0^t \frac{ (s^{\frac{3}{2}(\frac{1}{\alpha} - \frac{1}{p})} \| u ( s ) \|_{L^{p,\infty}} )}{(t -s )^{1- \frac{3}{2} (\frac{1}{\alpha} - \frac{1}{p})}s^{ \frac{3}{2} (\frac{1}{\alpha} - \frac{1}{p}) }} d s\\
 + C (K_0 \| b \|_{L^{3,\infty}} + \| w \|_{L^{3, \infty}}) \\
\cdot \sup_{\| \varphi \|_{L^{p',1}} \leq 1} \int_0^t (s^{ \frac{3}{2}(\frac{1}{\alpha} - \frac{1}{p}) }\| u ( s ) \|_{L^{p , \infty }} ) \| \nabla \mathrm{e}^{- (t -s ) L^*} \varphi \|_{L^{\frac{3 p}{2 p -3},1}} d s\\
\leq C_1 \| b \|_{L^{\alpha , \infty}}  + C_1 (K_0 \| b \|_{L^{3,\infty}} + \| w \|_{L^{3, \infty}}) \sup_{0< s <t}\{ s^{\frac{3}{2} (\frac{1}{\alpha} - \frac{1}{p}) } \| u ( s ) \|_{L^{p,\infty}} \}.
\end{multline*}
Here we used the fact that $\| u \|_{X_p} \leq K_0 \| b \|_{L^{3,\infty}}$. Since $C_1$ does not depend on both $t$ and $T$, we have
\begin{multline*}
\sup_{ 0<t <T}\{ t^{\frac{3}{2}(\frac{1}{\alpha} - \frac{1}{p} )}\| u (t) \|_{L^{p,\infty}} \} \\
\leq C_1 \| b \|_{L^{\alpha , \infty}} + C_1 (K_0 \| b \|_{L^{3,\infty}} + \| w \|_{L^{3,\infty}}) \sup_{ 0<t <T}\{ t^{\frac{3}{2}(\frac{1}{\alpha} - \frac{1}{p} )}\| u (t) \|_{L^{p,\infty}} \} .
\end{multline*}
Therefore we see that
\begin{equation*}
\sup_{ 0<t <T}\{ t^{\frac{3}{2}(\frac{1}{\alpha} - \frac{1}{p} )}\| u (t) \|_{L^{p,\infty}} \} \leq 2 C_{1} \| b \|_{L^{\alpha , \infty}} < + \infty
\end{equation*}
if $C_1 (K_0 \| b \|_{L^{3, \infty}} + \| w \|_{L^{3, \infty}}) \leq 1/2$. Since $C_{1}$ does not depend on both $t$ and $T$, we conclude that
\begin{equation*}
\sup_{t >0}\{ t^{\frac{3}{2}(\frac{1}{\alpha} - \frac{1}{p} )}\| u (t) \|_{L^{p,\infty}} \} \leq 2 C_1 \| b \|_{L^{\alpha , \infty } }< + \infty .
\end{equation*}
Thus, we see \eqref{r2}. From now on we assume that $C_1 (K_0 \| b \|_{L^{3, \infty} } + \| w \|_{L^{3 , \infty}}) \leq 1 /2$.

Next we prove \eqref{r4}. Using \eqref{eq211}, \eqref{eq23}, \eqref{eq24}, \eqref{e2}, \eqref{e4}, and \eqref{r2}, we see that
\begin{multline*}
\| u (t) \|_{L^{\alpha , \infty}} \leq C \sup_{ \varphi \in C_{0 , \sigma }^\infty, { \ }\| \varphi \|_{L^{\alpha',1} } \leq 1 }| \dual{ b , \mathrm{e}^{-t L^*} \varphi }|\\
+ C \sup_{ \varphi \in C_{0 , \sigma }^\infty, { \ }\| \varphi \|_{L^{\alpha',1} } \leq 1}  \left| \int_0^t \dual{u \otimes u + u \otimes w + u \otimes w , \nabla \mathrm{e}^{- (t-s )L^*} \varphi  } d s \right|  \\
\leq C \| b \|_{L^{\alpha, \infty}} +\\
\sup_{\| \varphi \|_{L^{\alpha',1} }\leq 1 } \int_0^t \| u \otimes u + u \otimes w + u \otimes w \|_{L^{\frac{3 p}{ p + 3}, \infty}} \| \nabla \mathrm{e}^{- (t-s )L^*} \varphi \|_{L^{\frac{3 p}{ 2 p - 3}, 1}} d s \\
\leq C \| b \|_{L^{\alpha , \infty}} + \int_0^t \frac{ ( \| u (s) \|_{L^{3,\infty} } + \| w \|_{L^{3, \infty}} ) (s^{\frac{3}{2}(\frac{1}{\alpha} - \frac{1}{p})} \| u ( s ) \|_{L^{p,\infty}} )}{(t -s )^{1 - \frac{3}{2}(\frac{1}{\alpha} - \frac{1}{p})}s^{ \frac{3}{2}(\frac{1}{\alpha} - \frac{1}{p})} } d s \leq C .
\end{multline*}
Thus, we have \eqref{r4}. Here $1/\alpha + 1/{\alpha'}=1$. Note that $3 p /(2 p - 3) > \alpha'$ since $3 p / (p + 3) < \alpha < 3$.

Thirdly, we derive \eqref{r5}. Using \eqref{eq21}, \eqref{r4}, and \eqref{r2}, we check that for $t >0$
\begin{equation*}
\| u (t) \|_{L^{3 , \infty}} \leq C \| u(t ) \|_{L^{\alpha , \infty} }^{ 1 - \frac{p (3 - \alpha) }{3 (p - \alpha)}} \| u (t) \|_{L^{p , \infty}}^{ \frac{p (3 - \alpha) }{3 (p - \alpha)}} \leq Ct^{-\frac{3}{2} (\frac{1}{\alpha} - \frac{1}{3})} ,
\end{equation*}
which is \eqref{r5}.

Finally, we attack $(\mathrm{ii})$ and $(\mathrm{iii})$. Assume that $\| \nabla w \|_{L^2} < + \infty$. Set
\begin{equation*}
\delta_1 := \min \bigg\{ \frac{\delta_\dagger}{2} , \frac{1}{4 C_1 K_0}, \frac{1}{4 C_1} \bigg\} .
\end{equation*}
We assume that $\| b \|_{L^{3,\infty}} < \delta_1$ and that $\| w \|_{L^{3,\infty}}< \delta_1$. Here $\delta_\dagger$ is the constant appearing in Lemma \ref{lem31}.

We shall derive \eqref{r3}. Assume that $b \in L^{3 , \infty}_{0, \sigma} ( \Omega )$. Fix $\varepsilon >0$. By the definition of $L^{3,\infty}_{0 , \sigma} (\Omega )$, there is $b_0 \in C_{0,\sigma}^\infty (\Omega)$ such that
\begin{equation*}
\| b - b_0 \|_{L^{3,\infty}} < \varepsilon \text{ and }\| b_0 \|_{L^{3,\infty}} < \min \{ \delta_0, \delta_1 \}.
\end{equation*}
From Lemma \ref{lem34}, we see that the system \eqref{eq15} admits a unique global-in-time generalized mild solution $u_*$ satisfying \eqref{r2}-\eqref{r5}. We may assume that $\| b_0 \|_{L^{3,\infty}} < \delta_{\dagger }/2$. From \cite[Theorem 2.8]{Kob13}, we see that there is a local-in-time strong $L^2$-solution $U$ of the system \eqref{eq15} with the initial datum $v_0 = b_0$, satisfying
\begin{equation*}
U \in C ([0,T_*); H^1_{0,\sigma} ( \Omega )) \cap C ((0,T_*); W^{2,2}( \Omega )) \cap C^1 ((0,T_*); L^2_\sigma ( \Omega ))
\end{equation*}
for some $T_* >0$. Here $H^1_{0,\sigma}(\Omega) = \overline{ C_{0,\sigma}^\infty (\Omega ) }^{\| \cdot \|_{W^{1,2}}}$. Note that $T_*$ depends only on $(\gamma, \Omega , w , \| b_0 \|_{W^{1,2}})$. See also \cite[Chapter 7]{Kob12}. Since $H_{0, \sigma}^1 (\Omega) \subset L^{3,\infty}_{0, \sigma}(\Omega)$, we see that there is $T_{**}>0$ such that
\begin{equation}\label{h1}
\sup_{0 \leq t \leq T_{**}}\{ \| U (t) \|_{L^{3,\infty}}\} \leq \delta_{ \dagger} .
\end{equation}
Since $\mathrm{e}^{- t L}$ is an analytic semigroup on $L^2_\sigma (\Omega)$, we easily check that for $0<t <T_{**}$
\begin{equation*}
U (t) = \mathrm{e}^{- t L}b_0 - \int_0^t \mathrm{e}^{-( t -s ) L} P \{ (U, \nabla ) U + (U, \nabla ) w + (w , \nabla) U \} d s
\end{equation*}
and that for each $\varphi \in C^\infty_{0 , \sigma } ( \Omega)$ and $0 < t < T_{**}$
\begin{equation*}
\dual{U (t) , \varphi } = \dual{ b_0 , \mathrm{e}^{-t L^*} \varphi } + \int_0^t \dual{U \otimes U(s) + U(s) \otimes w + w \otimes U(s) , \nabla \mathrm{e}^{-(t-s) L^*} \varphi } d s.
\end{equation*}
From \eqref{h1} and Lemma \ref{lem31}, we find that $u_* = U$ on $(0,T_{**})$. Applying \eqref{eq32}, we observe that
\begin{multline*} 
\| u(t) - b \|_{L^{3,\infty}} \leq \| u(t) - u_*(t) \|_{L^{3,\infty}} + \| u_*(t) - b_0 \|_{L^{3,\infty}} + \| b - b_0 \|_{L^{3,\infty}}\\
\leq (K_\dagger + 1) \| b - b_0 \|_{L^{3,\infty}} + \| U(t) - b_0 \|_{L^{3,\infty}} \leq C \varepsilon \text{ as }t \to 0+ 0.
\end{multline*}
Since $\varepsilon$ is arbitrary, we see \eqref{r3}. Similarly, we deduce $(\mathrm{iii})$.
\end{proof}
Finally we derive the stability and continuity of the generalized mild solutions.
\begin{proof}[Proof of Theorem \ref{thm110}]
Combining Theorem \ref{thm19} and Lemma \ref{lem34}, we see the assertion $(\mathrm{i})$ of Theorem \ref{thm110}. Now we prove the assertions $(\mathrm{ii})$ and $(\mathrm{iii})$ of Theorem \ref{thm110}. Let $3p/(p+3) < \alpha < 3$, and let $\delta_0$, $\delta_1$ be the constants appearing in Theorem \ref{thm19} and Lemma \ref{lem34}, respectively. Assume that
\begin{equation*}
a \in L^{3, \infty}_\sigma (\Omega) , { \ }\| a \|_{L^{3,\infty}} < \min \{ \delta_0,  \delta_1 \}, \text{ and } \| w \|_{L^{3,\infty } } < \min\{ \delta_0 ,  \delta_1 \} .
\end{equation*}
Let $v$ be the global-in-time generalized mild solution of the system \eqref{eq15} with the initial datum $v_0 = a$, obtained by Theorem \ref{thm19}. Assume that $\| \nabla w \|_{L^2} < \infty$.

We attack $(\mathrm{ii})$. We assume that $a \in L^{3,\infty}_{0 , \sigma} (\Omega )$. We shall prove that
\begin{align}
& \lim_{t \to 0 + 0} \| v ( t ) -a \|_{L^{3,\infty}} = 0,\label{qq1}\\
& \lim_{t \to \infty} \| v (t) \|_{L^{3,\infty}} = 0\label{qq2}.
\end{align}
Fix $\varepsilon >0$. By the definition of $L^{3,\infty}_{0 , \sigma} (\Omega )$, there is $b \in C_{0,\sigma}^\infty (\Omega)$ such that
\begin{equation*}
\| a - b \|_{L^{3,\infty}} < \varepsilon \text{ and }\| b \|_{L^{3,\infty}} < \min \{ \delta_0, \delta_1 \}.
\end{equation*}
From Lemma \ref{lem34}, we see that the system \eqref{eq15} admits a unique global-in-time generalized mild solution $u$ satisfying \eqref{r2}-\eqref{r5}. Applying \eqref{eq32} and \eqref{r3}, we see that for $t>0$
\begin{multline*} 
\| v(t) - a \|_{L^{3,\infty}} \leq \| v(t) - u(t) \|_{L^{3,\infty}} + \| u(t) - b \|_{L^{3,\infty}} + \| b - a \|_{L^{3,\infty}}\\
\leq C \| a - b \|_{L^{3,\infty}} + \| u(t) - b \|_{L^{3,\infty}} \leq C \varepsilon \text{ as }t \to 0+ 0.
\end{multline*}
Since $\varepsilon$ is arbitrary, we see \eqref{qq1}. Similarly, we use \eqref{eq32} and \eqref{r5} to observe that 
\begin{multline*} 
\| v(t) \|_{L^{3,\infty}} \leq \| v(t) - u(t) \|_{L^{3,\infty}} + \| u(t) \|_{L^{3,\infty}} \\
\leq C \| a - b \|_{L^{3,\infty}} + \| u(t) \|_{L^{3,\infty}} \leq C \varepsilon \text{ as }t \to \infty.
\end{multline*}
Note that $b \in C^\infty_{0 , \sigma} (\Omega) \subset L^{\alpha, \infty}_\sigma (\Omega)$. Since $\varepsilon$ is arbitrary, we see \eqref{qq2}. Therefore we finish prove the assertion $(\mathrm{ii})$. 

The proof of the assertion $(\mathrm{iii})$ of Theorem \ref{thm110} is left to the reader, since we can prove $(\mathrm{iii})$ by the same argument to show $(\mathrm{ii})$ with the assertions $(\mathrm{i})$ and $(\mathrm{ii})$ of Lemma \ref{lem33}. Therefore Theorem \ref{thm110} is proved. 
\end{proof}

 \section{Global-in-Time Mild Solution}\label{sect4}

In this section we show the existence of a unique global-in-time mild solution of the system \eqref{eq15} when $v_0 \in L^{3,\infty}_{0 , \sigma} ( \Omega)$ and both $\| v_0 \|_{L^{3,\infty}}$ and $\| w \|_{L^{3,\infty}}$ are sufficiently small.  

Throughout this section we assume that $w$ is as in Assumption \ref{assB} and that $|u_\infty | \leq \gamma$ for some $\gamma > 0$. {\color{red}{Let $6 < p < \infty$ such that ${3 \beta}/(2 \beta -3)< p$. Fix $p$. Let $p'$ be the positive number such that $1/p+1/{p'}=1$}}. The symbols $L$ and $L^*$ represent the two linear operators defined by Subsection \ref{subsec23}. The aim of this section is to prove Theorem \ref{thm111}. In order to prove Theorem \ref{thm111}, we first constructs a unique local-in-time mild solution of the system \eqref{eq15}. Secondly, we prove that the system \eqref{eq15} admits a unique global-in-time mild solution of the system \eqref{eq15} when $\| v_0 \|_{L^{3,\infty}}$ and $\| w \|_{L^{3,\infty}}$ are sufficiently small. Moreover, we derive $L^\infty$-decay rate with respect to time of the solution.

\subsection{Local-in-Time Mild $L^{3,\infty}$-Solutions}\label{subsec41} 
Let us construct a unique local-in-time mild solution of the system \eqref{eq15}. For $T \in (0,\infty]$, we define the Banach space $Z_T$ as follows: 
\begin{equation*}
Z_T := \{ f \in BC ((0,T) ; L_\sigma^{3,\infty} (\Omega)) \cap C((0, T ) ; L^{p,\infty}_\sigma (\Omega )) ; { \ }\| f \|_{Z_T} < + \infty \},
\end{equation*}
\begin{multline*}
\| f \|_{Z_T}:= \\
\sup_{0<t <T} \{ \| f (t) \|_{L^{3,\infty}} \} + \sup_{0<t <T} \{ t^\frac{p- 3}{2 p} \| f (t ) \|_{L^{p,\infty}} \} + \sup_{0<t <T}\{ t^\frac{1}{2}\| \nabla f (t) \|_{L^{3,\infty}} \} .
\end{multline*}
We now attack the following proposition.
\begin{proposition}\label{prop41}
There are $\delta_* = \delta_* (\gamma ,p )>0$ and $K_{*} = K_{*} (\gamma , p ) >0$ such that if
\begin{equation*}
a \in L^{3,\infty}_{0, \sigma} (\Omega ),{ \ }\| a \|_{L^{3,\infty}} < \delta_*, \text{ and } \| w \|_{L^{3,\infty}} < \delta_* ,
\end{equation*}
then the system \eqref{eq15} with the initial datum $v_0 = a$ admits a unique local-in-time mild solution $v$ in $Z_{T_*}$ for some $T_* >0$. Moreover,
\begin{align}
& \| v \|_{Z_{T_*}} \leq K_{*} \| a \|_{L^{3,\infty}}, \label{eq41}\\
& v \in BC ([0,T_*) ; L_{0, \sigma}^{3,\infty} (\Omega)) \label{eq42}.
\end{align}
Moreover, assume in addition that $w$ satisfies Assumption $\ref{assC}$. Then for each fixed $0 < t < T_*$,
\begin{equation}\label{eq43}
\nabla v (t ) \in L^{3 , \infty}_{0} ( \Omega ) .
\end{equation} 
Here $T_* $ depends only on $\gamma$, $\beta$, $p$, $\| w \|_{L^{  \frac{ 3 \beta }{ 3-\beta}, \infty }}$, and $\| \nabla w \|_{L^{\beta , \infty}}$. 
\end{proposition}
\begin{proof}[Proof of Proposition \ref{prop41}]
Let $T>0$. For $m \in \mathbb{N}$ and $t \in ( 0, T ) $,
\begin{align*}
& v_1 ( t ) := \mathrm{e}^{ - t L } a,\\
& v_{m + 1}(t) := \mathrm{e}^{-t L} a - \int_0^t \mathrm{e}^{- ( t - s) L } P \{ (v_m , \nabla ) v_m+ (v_m, \nabla ) w + (w , \nabla ) v_m \} d s .
\end{align*}
We now solve the above integral forms in $Z_T$ to construct a local-in-time mild solution of the system \eqref{eq15}. We first consider $\| v_1 \|_{Z_T}$. By \eqref{e1} and \eqref{eE3}, we see that
\begin{align*}
\| \mathrm{e}^{-t L} a \|_{L^{3,\infty}} \leq & C \| a \|_{L^{3,\infty}},\\
\| \mathrm{e}^{- t L} a \|_{L^{p,\infty}} \leq & C t^{- \frac{p - 3}{2 p}} \| a \|_{L^{3,\infty} } , \\ 
\| \nabla \mathrm{e}^{- t L} a \|_{L^{3,\infty} }\leq & C t^{- \frac{1}{2}} \| a \|_{L^{3,\infty}} .
\end{align*}
Therefore we obtain
\begin{equation}\label{eq44}
\|  v_1 \|_{Z_T} \leq C (\gamma , p) \| a \|_{L^{3,\infty}}.
\end{equation}
Next we consider $\| v_{m+1} \|_{Z_T}$. Using \eqref{eq211}, \eqref{eq23}, \eqref{eq24}, and \eqref{e8}, we check
\begin{multline*}
\left\| \int_0^t \mathrm{e}^{-(t-s) L} P ( v_m , \nabla ) v_m (s) d s \right\|_{L^{3,\infty}} \\
\leq C \sup_{\varphi \in C_{0,\sigma}^\infty, { \ }\| \varphi \|_{L^{\frac{3}{2} , 1}} \leq 1} \bigg| \int_0^t \dual{ v_m (s) \otimes v_m (s) , - \nabla \mathrm{e}^{- (t-s) L^*} \varphi } d s \bigg|\\
\leq C \sup_{\| \varphi \|_{L^{\frac{3}{2} , 1}} \leq 1}  \int_0^t \| v_m (s) \|_{L^{3 , \infty}}^2 \| \nabla \mathrm{e}^{- ( t- s)L^*} \varphi \|_{L^{3,1}} d s \leq C \| v_m \|_{Z_T}^2.
\end{multline*}
Similarly, we have
\begin{equation*}
\left\| \int_0^t \mathrm{e}^{-(t-s) L} P \{ ( v_m (s) , \nabla ) w + (w , \nabla ) v_m(s) \} d s \right\|_{L^{3,\infty}} \leq C \| w \|_{L^{3,\infty}} \| v_m \|_{Z_T} .
\end{equation*}
Applying \eqref{eq211}, \eqref{eq23}, \eqref{eq24}, \eqref{e4}, and \eqref{e8}, we see that
\begin{multline*}
\left\| \int_0^t \mathrm{e}^{- (t -s ) L} P (v_m(s) , \nabla ) v_m (s) d s \right\|_{L^{p,\infty}} \\
\leq C \sup_{ \varphi \in C_{0 , \sigma }^\infty, { \ }\| \varphi \|_{L^{p',1}} \leq 1 } \left| \int_0^t \dual{ v_m (s) \otimes v_m (s), - \nabla \mathrm{e}^{- (t-s )L^*} \varphi  }  d s \right| \\
\leq C \int_0^{\frac{t}{2} } \frac{ \| v_m ( s ) \|_{L^{3, \infty}} (s^\frac{p - 3}{2 p} \| v_m ( s ) \|_{L^{p,\infty}} )}{s^\frac{p -3}{2 p} ( t - s )  } d s\\
 + C \sup_{ \| \varphi \|_{L^{p',1}} \leq 1 } \int_{\frac{t}{2}}^t \frac{\| v_m (s) \|_{L^{3,\infty}} (s^\frac{p - 3}{2 p} \| v_m ( s ) \|_{L^{p , \infty }} )}{s^\frac{p - 3}{2 p}} \| \nabla \mathrm{e}^{- (t -s ) L^*} \varphi \|_{L^{\frac{3 p}{2 p -3},1}} d s \\
\leq C t^{-\frac{p- 3}{2 p}}\| v_m \|_{Z_T}^2 .
\end{multline*}
Similarly, we obtain
\begin{equation*}
\left\| \int_0^t \mathrm{e}^{- (t -s ) L} P \{ (v_m, \nabla ) w + (w , \nabla ) v_m \} d s \right\|_{L^{p,\infty}} \leq C t^{-\frac{p- 3}{2 p}} \| w \|_{L^{3,\infty}} \| v_m \|_{Z_T}.
\end{equation*}
Therefore we find that
\begin{multline}\label{eq45}
\sup_{0< t <T}\{ \| v_{m+1}(t) \|_{L^{3,\infty}} \} + \sup_{0 < t <T}\{ t^\frac{p-3}{2 p}\| v_{m+1} (t) \|_{L^{p,\infty}}\}\\
 \leq C(\gamma, p) ( \| v_1 \|_{Z_T} + \| w \|_{L^{3,\infty}}\| v_m \|_{Z_T} + \| v_m \|_{Z_T}^2 ).
\end{multline}
Now we consider $\| \nabla v_{m+1} (t) \|_{L^{3,\infty}}$. By Lemma \ref{lem25}, \eqref{eq23}, \eqref{eq24}, and \eqref{e16}, we check that for each $j =1,2,3$,
\begin{multline*}
\left\| \partial_j \int_0^t \mathrm{e}^{- (t -s ) L}  P ( v_m , \nabla ) v_m ( s ) d s   \right\|_{L^{3,\infty } } \\
\leq C \sup_{ \phi \in [ C_0^\infty ]^3, { \ }\| \phi \|_{L^{\frac{3}{2} , 1}} \leq 1} \left| \int_0^t \dual{ (v_m (s), \nabla ) v_m (s) , \mathrm{e}^{- (t - s )L^*} P( -\partial_j \phi ) }  d s \right| \\
\leq C \sup_{\| \phi \|_{L^{\frac{3}{2} , 1}} \leq 1}\int_0^t \frac{ (s^\frac{p - 3}{2 p} \| v_m (s) \|_{L^{p , \infty}}) (s^{\frac{1}{2}} \| \nabla v_m (s) \|_{L^{3,\infty}} )}{ s^{\frac{2 p - 3}{2 p}}} \\
\cdot \| \mathrm{e}^{- (t-s) L^*} P \partial_j \phi \|_{L^{\frac{3 p}{2 p - 3} , 1}} d s\\
\leq C \| v_m \|_{Z_T}^2 \int_0^t \frac{1}{ s^{\frac{2 p - 3}{2 p}} (t - s)^{\frac{p + 3}{2 p}}} d s \leq C t^{- \frac{1}{2}} \| v_m \|^2_{Z_T}.
\end{multline*}
Similarly, we observe that for each $j =1,2,3$,
\begin{multline*}
\left\| \partial_j \int_0^t \mathrm{e}^{- (t -s ) L}  P ( v_m , \nabla ) w d s   \right\|_{L^{3,\infty } } \\
\leq C \int_0^t \frac{ ( s^\frac{p - 3}{2 p} \| v_m (s) \|_{L^{p , \infty}} ) \| \nabla w \|_{L^{\beta , \infty}} }{ s^{\frac{p -3}{2 p}}(t - s )^{\frac{3 }{2 \beta} + \frac{3}{2 p}}  } d s \\
 \leq C t^{- \frac{3 - \beta}{2 \beta}} \| v_m \|_{Z_T} \| \nabla w \|_{L^{\beta , \infty}} .
\end{multline*}
Note that $3/(2 \beta) + 3/(2p) < 1$ since $(3 \beta )/(2 \beta - 3 ) <p $. We also check that for each $j=1,2,3$,
\begin{multline*}
\left\| \partial_j \int_0^t \mathrm{e}^{- (t -s ) L}  P ( w , \nabla ) v_m d s   \right\|_{L^{3,\infty } } \leq C \sup_{ \phi \in [ C_0^\infty]^3, { \ }\| \phi \|_{L^{\frac{3}{2},1}} \leq 1} \bigg\{ \\
\int_0^t \frac{ ( s^\frac{1}{2} \| \nabla v_m (s) \|_{L^{3, \infty}} ) \| w \|_{L^{ \frac{3 \beta}{3 - \beta} , \infty}} }{ s^{\frac{1}{2}}} \| \mathrm{e}^{- (t- s)L^*} P \partial_j \phi  \|_{L^{\frac{ \beta }{\beta - 1} ,1}} d s \bigg\} \\
 \leq C \| w \|_{L^{ \frac{3 \beta }{ 3 - \beta  } , \infty}} \| v_m \|_{Z_T} \int_0^t \frac{1}{ s^{\frac{1}{2}} (t - s)^{\frac{3 }{2 \beta } }} d s\\
\leq C t^{- \frac{3 - \beta }{2 \beta }} \| w \|_{L^{ \frac{3 \beta }{3 - \beta } , \infty}} \| v_m \|_{Z_T} .
\end{multline*}
Note that $3/2 < \beta /(\beta -1) <3$. As a result, we have
\begin{multline}\label{eq46}
\sup_{0<t <T}\{ t^\frac{1}{2} \| \nabla v_m (t) \|_{L^{3,\infty}} \} \\
\leq C (\gamma, p ) ( \| a \|_{L^{3,\infty}} + \| v_m \|_{Z_T}^2 + \| w \|_{L^{3,\infty}} \| v_m \|_{Z_T})\\
  + C (\gamma , \beta, p ) T^\frac{2\beta - 3}{2 \beta} ( \| w \|_{L^{ \frac{3 \beta }{3 - \beta }  , \infty}}  +  \| \nabla w \|_{L^{\beta , \infty}} )  \| v_m \|_{Z_T} .
\end{multline}
From \eqref{eq44}-\eqref{eq46}, we conclude that there are $C_* = C_* (\gamma, p) >0$ and $C_{**} = C_{**} (\gamma, \beta ,p) >0$ such that
\begin{equation*}
\| v_1 \|_{Z_T} \leq C_* \| a \|_{L^{3,\infty}},
\end{equation*}
\begin{multline*}
\| v_{m+1} \|_{Z_T} \leq C_* ( \| a \|_{L^{3,\infty}} + \| v_m \|_{Z_T}^2 +  \| w \|_{L^{3,\infty}} \| v_m \|_{Z_T}) \\
+ C_{**} T^\frac{2\beta - 3}{2 \beta} ( \| w \|_{L^{ \frac{3 \beta }{3 - \beta } , \infty}} +  \| \nabla w \|_{L^{\beta , \infty}} ) \| v_m \|_{Z_T}.
\end{multline*}
By the previous arguments, we find that
\begin{multline*}
\| v_{m+2} - v_{m+1} \|_{Z_T} \leq \{ C_* \| v_m \|_{Z_T} + C_* \| v_{m+1} \|_{Z_T} + C_* \| w \|_{L^{3,\infty}} \\
+  C_{**} T^\frac{2\beta - 3}{2 \beta} ( \| w \|_{L^{ \frac{3 \beta}{3 - \beta}, \infty}} + \| \nabla w \|_{L^{\beta , \infty}} ) \} \| v_{m+1} - v_m \|_{Z_T} .
\end{multline*}
From now on we assume that
\begin{align*}
& \| a \|_{L^{3,\infty}} \leq \frac{1}{16 C_*^2},\\
& \| w \|_{L^{3,\infty}} \leq \frac{1}{8 C_*},
\end{align*}
and we choose $T_* >0$ such that
\begin{equation*}
T_*^\frac{2 \beta - 3}{2 \beta} \leq \frac{1}{8 C_{**} ( \| w \|_{L^{\frac{3 \beta }{3 - \beta} , \infty}} + \| \nabla w \|_{L^{\beta,\infty}} )} .
\end{equation*}
By induction, we check that for each $m \in \mathbb{N}$
\begin{align*}
& \| v_m \|_{Z_{T_*}} \leq 2 C_* \| a \|_{L^{3,\infty}} ,\\
& \| v_{m+2} - v_{m+1} \|_{Z_{T_*}} \leq \frac{1}{2} \| v_{m+1} - v_m \|_{Z_{T_*}} .
\end{align*}
Since $Z_{T_*}$ is a Banach space, we apply a fixed point argument to obtain a unique function $v \in Z_{T_*}$ satisfying
\begin{align}
& \| v \|_{Z_{T_*}} \leq 2 C_* \| a \|_{L^{3,\infty}}, \label{eq47}\\
& \| v_m - v \|_{Z_{T_*}} \to 0 \text{ as } m \to \infty .\label{eq48}
\end{align}
Since
\begin{equation*}
v_{m + 1}(t) = \mathrm{e}^{-t L} a - \int_0^t \mathrm{e}^{- ( t - s) L } P \{ (v_m , \nabla ) v_m + (v_m, \nabla ) w + (w , \nabla ) v_m \} d s,
\end{equation*}
we apply \eqref{eq47} and \eqref{eq48} to check that $v$ satisfies for $0<t <T_*$
\begin{equation*}
v(t) = \mathrm{e}^{-t L} a - \int_0^t \mathrm{e}^{- ( t - s) L } P \{ (v(s) , \nabla ) v (s) + (v (s), \nabla ) w + (w , \nabla ) v ( s ) \} d s.
\end{equation*}
It is easy to check that for all $\varphi \in C_{0,\sigma}^\infty (\Omega)$ and $0<t <T_*$
\begin{equation*}
\dual{v (t) , \varphi} = \dual{ a , \mathrm{e}^{-t L^*} \varphi} + \int_0^t \dual{v(s) \otimes v(s) + v(s) \otimes w + w \otimes v(s) , \nabla \mathrm{e}^{-(t-s) L^*} \varphi} d s.
\end{equation*}
From Lemma \ref{lem31}, we find that $v$ is a unique if $\| a \|_{L^{3,\infty}} \leq \delta_{ \dagger}/(2C_*)$. Here $\delta_{\dagger }$ is the constant appearing in Lemma \ref{lem31}.

Next, we show that $v \in BC ([0,T_*); L^{3,\infty}_{0,\sigma} (\Omega))$. Since $\mathrm{e}^{-t L}a \in L^{3,\infty}_{0,\sigma}$ for $t \geq 0$, we only have to show that for fixed $0<t <T_*$
\begin{equation*}
\int_0^t \mathrm{e}^{- ( t - s) L } P \{ (v (s) , \nabla ) v(s) + (v (s), \nabla ) w + (w , \nabla ) v(s) \} d s \in L^{3,\infty}_{0,\sigma}( \Omega ).
\end{equation*}
To this end, we show that
\begin{equation*}
\left\| \int_0^t \mathrm{e}^{- ( t - s) L } P \{ (v (s), \nabla ) v (s)+ (v (s), \nabla ) w + (w , \nabla ) v(s) \} d s \right\|_{L^{\alpha_1, \infty}} < + \infty
\end{equation*}
for some $3p/(p+3) < \alpha_1 <3 $. Let $3p/(p+3) < \alpha_1 <3 $, and fix $\alpha_1$. Using \eqref{eq211}, \eqref{eq23}, \eqref{eq24}, and \eqref{e4}, we see that
\begin{multline*}
\left\| \int_0^t \mathrm{e}^{- (t -s ) L} P \{ (v(s) , \nabla ) v(s) + (v (s), \nabla ) w + (w , \nabla ) v(s) \} d s \right\|_{L^{\alpha_1,\infty}} \\
= \sup_{ \varphi \in C_{0 , \sigma }^\infty, { \ }\| \varphi \|_{L^{\alpha_1',1} } \leq 1 } \left| \int_0^t \dual{v \otimes v + v \otimes w + v \otimes w , - \nabla \mathrm{e}^{- (t-s )L^*} \varphi  }  d s \right|\\
\leq \sup_{ \| \varphi \|_{L^{\alpha_1',1}} \leq 1 } \int_0^t \| v \otimes v + v \otimes w + v \otimes w \|_{L^{\frac{3 p}{ p + 3}, \infty}} \| \nabla \mathrm{e}^{- (t-s )L^*} \varphi \|_{L^{\frac{3 p}{ 2 p - 3}, 1}} d s \\
\leq \int_0^t \frac{ ( \| v (s) \|_{L^{3,\infty} } + \| w \|_{L^{3, \infty}} ) (s^{\frac{p-3}{2 p}} \| v ( s ) \|_{L^{p,\infty}} )}{ s^{\frac{p - 3}{2 p}} (t -s )^{1 - \frac{3}{2}(\frac{1}{\alpha_1} - \frac{1}{p})} } d s\\
\leq C t^\frac{3 - \alpha_1}{2\alpha_1}(\| a \|_{L^{3,\infty}} + \| w \|_{L^{3,\infty}} )\| v \|_{Z_T} \leq CT_*^\frac{3 - \alpha_1}{2 \alpha_1} < + \infty.
\end{multline*}
Note that $ \alpha_1' < 3 p/(2p - 3) \leq 3$ since $3 p/(p + 3) < \alpha_1$. Since for fixed $0<t <T_*$
\begin{equation*}
\int_0^t \mathrm{e}^{- ( t - s) L } P \{ (v , \nabla ) v + (v, \nabla ) w + (w , \nabla ) v\} d s \in L_\sigma^{p,\infty}(\Omega) \cap L_\sigma^{\alpha_1, \infty} (\Omega),
\end{equation*}
it follows from Lemma \ref{lem21} that for fixed $0 < t <T_*$
\begin{equation*}
\int_0^t \mathrm{e}^{- ( t - s) L } P \{ (v , \nabla ) v + (v, \nabla ) w + (w , \nabla ) v\} d s \in L^{3,\infty}_{0,\sigma}( \Omega ).
\end{equation*}
Therefore we conclude that $v \in BC ([0,T_*); L^{3,\infty}_{0,\sigma} (\Omega))$.

Finally, we deduce \eqref{eq43}. Assume that $w$ satisfies Assumption \ref{assC}. Since $a \in L^{3,\infty}_{0 , \sigma} ( \Omega )$, it follows from Lemma \ref{lem212} to see that $\nabla v_1 (t) \in L^{3, \infty}_0 ( \Omega )$ for each $t>0$. From $v_1 (t) \in L_0^{3,\infty} ( \Omega ) \cap L^{p,\infty} ( \Omega )$, $\nabla v_1 (t) \in L^{3,\infty}_0 (\Omega )$, and Assumption \ref{assC}, we easily check that $\{ ( v_1 , \nabla ) v_1 + ( v_1 , \nabla ) w + ( w , \nabla ) v_1 \}(t) \in L^{ 3 , \infty }_0 ( \Omega ) $ for each $0 < t <T_*$. Therefore we see that $\nabla v_2 (t) \in L^{3,\infty}_0 ( \Omega )$. By induction, we see that for $m \in \mathbb{N}$ and $0 < t < T_*$, $\nabla v_m (t) \in L^{3 , \infty}_0 ( \Omega )$. This implies that for $0  < t < T_*$, $\nabla v (t) \in L^{3 , \infty}_0 ( \Omega )$. Therefore Proposition \ref{prop41} is proved.
\end{proof}

\subsection{Global-in-Time Mild Solution and $L^\infty$-Asymptotic Stability}\label{subsec42} 
We now construct a unique global-in-time mild solution of the system \eqref{eq15}, satisfying $L^\infty$-asymptotic stability. Let $3p/(p + 3) < \alpha <3$. Let $\delta_0$, $K_0$ be the constants appearing in Theorem \ref{thm19}, let $\delta_1$ be the constant appearing in Theorem \ref{thm110}, let $\delta_{\dagger}$ be the constant appearing in Lemma \ref{lem31}, let $\delta_{\dagger \dagger}$ be the constant appearing in Lemma \ref{lem33}, and let $\delta_*$, $K_*$ be the constants appearing in Proposition \ref{prop41}.
\begin{proof}[Proof of Theorem \ref{thm111}]
Assume that $a \in L^{3 , \infty}_{0 , \sigma} ( \Omega )$,
\begin{align*}
& \| a \|_{L^{3,\infty}} \leq \min \left\{ \delta_0, \delta_1 , \delta_*, \delta_{\dagger}/{K_0}, \delta_{\dagger}/{K_*}, \min \{ \delta_0 , \delta_1 , \delta_{\dagger} , \delta_* \}/(K_*)^2 \right\} ,\\
& \| w \|_{L^{3,\infty}} \leq \min \{ \delta_0 , \delta_1, \delta_*, {\delta_{\dagger }/K_0}, \delta_{\dagger}/{K_*} \} .
\end{align*}
Let $v_\sharp$ be the global-in-time generalized mild solution of the system \eqref{eq15} with the initial datum $v_0 = a$, obtained by Theorem \ref{thm19}, satisfying
\begin{align*}
& \sup_{0 < t < \infty} \{ \| v_\sharp (t) \|_{L^{3,\infty}} \} \leq K_0 \| a \|_{L^{3,\infty } } \leq \delta_{\dagger } ,\\
& \sup_{0 < t < \infty} \{ t^\frac{p -3}{2 p} \| v_\sharp (t) \|_{L^{3,\infty}} \} \leq K_0 \| a \|_{L^{3,\infty } } \leq \delta_{\dagger }.
\end{align*}
From Proposition \ref{prop41}, we see that there exists a unique local-in-time mild solution in $Z_{T_*}$ of the system \eqref{eq15} with the initial datum $v_0 = a$ for some $T_* >0$, satisfying
\begin{align*}
& v \in BC ([ 0, T_* ) ; L^{3,\infty}_{0,\sigma} (\Omega ) ), \\
& \sup_{0 \leq t <  T_*} \{ \| v (t) \|_{L^{3,\infty}} \} \leq K_* \| a \|_{L^{3,\infty}} \leq \delta_{\dagger} .
\end{align*}
From Theorem \ref{thm19} and Lemma \ref{lem31}, we find that $v_\sharp ( t ) = v (t)$ on $[0, T_*)$.

Now we let $\tau \in (0,T_*)$. Since $\| v_\sharp (\tau) \|_{L^{3,\infty}} \leq \min \{ \delta_0 , \delta_* \}$, we apply Proposition \ref{prop41} to see that there exists a unique local-in-time mild solution of the system \eqref{eq15} with the initial datum $v_\sharp (\tau)$ for some $T_{**} >0$, satisfying 
\begin{align*}
v \in BC ([ \tau, \tau + T_{**}) ; L^{3,\infty}_{0,\sigma} (\Omega ) ), \\
\sup_{\tau \leq t < \tau + T_{**}} \{ \| v (t) \|_{L^{3,\infty}} \} \leq K_* \| v_\sharp ( \tau ) \|_{L^{3,\infty}} \\
\leq (K_*)^2 \| a \|_{L^{3,\infty}} \leq \min \{ \delta_0 , \delta_1 , \delta_{ \dagger} , \delta_* \}. 
\end{align*}
Note that $T_{**}$ depends only on $(\gamma , \beta , p , \| w \|_{L^{\frac{3 \beta}{3 - \beta}, \infty}} , \| \nabla w \|_{L^{\beta , \infty}} , K_0 , \delta_0)$, that is, $T_{**}$ does not depend on $\tau$. From Lemma \ref{lem31}, we see that $v_\sharp ( t ) = v (t)$ on $(\tau , \tau + T_{**})$. Thus, we repeat the above argument to construct a unique global-in-time mild solution of the system \eqref{eq15} with the initial datum $v_0 = a$.

Next we show that for each fixed $T >0$
\begin{align}
& \sup_{0  < t <T} \{ t^{ \frac{1}{2} } \| \nabla v (t) \|_{L^{3,\infty}} \} \leq C (T) < + \infty , \label{eq49}\\
& \sup_{0 < t < T} \{ t^\frac{1}{2}   \| v ( t ) \|_{L^\infty} \} \leq C (T) < \infty \label{eq410}.
\end{align}
Now we derive \eqref{eq49}. Fix $T >0$. By an argument similar to that in the proof of Proposition \ref{prop41}, we see that there is $C >0$ independent of both $t$ and $T$ such that for $0<t <T$ 
\begin{multline*}
t^\frac{1}{2} \| \nabla v (t) \|_{L^{3, \infty } } \leq C \| a \|_{L^{3, \infty}} + C T^\frac{2 \beta -3}{2 \beta } K_0 \| a \|_{L^{3,\infty}} \| \nabla w \|_{L^{\beta , \infty}} \\
+ \int_0^t C t^{ \frac{1}{2} } \cdot \frac{(s^\frac{2 p -3}{2 p} \| v (s) \|_{L^{p,\infty}} )(s^\frac{1}{2} \| \nabla v (s) \|_{L^{3,\infty}} )}{s^\frac{2 p -3}{2 p} ( t - s )^\frac{p + 3}{2 p} } d s\\
+ \int_0^t C t^{ \frac{1}{2} } \cdot \frac{\| w \|_{L^{ \frac{3 \beta}{3 - \beta} , \infty}} (s^{ \frac{1}{2}} \| \nabla v (s) \|_{L^{3, \infty}}) }{s^\frac{1}{2}( t - s)^\frac{3}{2 \beta} } d s .
\end{multline*}
Applying the Gronwall inequality, we check that
\begin{equation*}
\sup_{0  < t <T} \{ t^{ \frac{1}{2} } \| \nabla v (t) \|_{L^{3,\infty}} \} \leq C (T) < + \infty .
\end{equation*}
Here we used the facts that
\begin{align*}
& \int_0^t \frac{t^\frac{1}{2} }{s^\frac{2 p -3}{2 p} (t - s)^\frac{p + 3}{2 p} } d s \leq C,\\
& \int_0^t \frac{t^\frac{1}{2}}{ s^\frac{1}{2} ( t - s )^{\frac{3}{2 \beta}} } d s \leq C T^\frac{2 \beta - 3}{2 \beta } .
\end{align*}
Next we derive \eqref{eq410}. Let $t \in (0,T]$. We see at once that
\begin{multline*}
\| v (t) \|_{L^{\infty } } \leq \|  \mathrm{e}^{-t L} a \|_{L^{\infty} } + \left\|  \int_0^t \mathrm{e}^{- ( t - s) L} P (v, \nabla ) v (s) d s \right\|_{L^{\infty} } \\
 + \left\|  \int_0^t \mathrm{e}^{- ( t - s) L} P (v(s), \nabla ) w d s \right\|_{L^{\infty} } + \left\|  \int_0^t \mathrm{e}^{- ( t - s) L} P ( w, \nabla ) v (s) d s \right\|_{L^{\infty} }.
\end{multline*}
From \eqref{p5}, we have
\begin{equation*}
\| \mathrm{e}^{- t L} a \|_{L^\infty } \leq C t^{-\frac{1}{2} } \| a \|_{L^{3,\infty}} .
\end{equation*}
Applying \eqref{eq26}, \eqref{eq23}, \eqref{eq24}, and \eqref{p4}, we observe that
\begin{multline*}
\left\| \int_0^t \mathrm{e}^{- ( t - s ) L} P (v, \nabla ) v (s) d s \right\|_{L^\infty} \\
\leq C \sup_{\phi \in [C_0^\infty ]^3, { \ }\| \phi \|_{L^1} \leq 1} \left| \int_0^t \dual{ (v ( s )  , \nabla ) v (s) , \mathrm{e}^{-(t -s) L^*} P \phi } d s \right| \\
\leq C \int_0^t \frac{ ( s^\frac{p - 3}{2 p} \| v (s) \|_{L^{p , \infty}}) (s^\frac{1}{2} \| \nabla v ( s ) \|_{L^{3 , \infty}}) }{ ( t - s)^\frac{p + 3}{2 p}s^{\frac{2 p - 3}{2 p}} } d s \leq C (T) t^{- \frac{1}{2} } .
\end{multline*}
In the same manner, we check that
\begin{equation*}
\left\| \int_0^t \mathrm{e}^{- ( t - s ) L} P (v, \nabla ) w d s \right\|_{L^\infty} + \left\| \int_0^t \mathrm{e}^{- ( t - s ) L} P (w, \nabla ) v d s \right\|_{L^\infty} \leq C(T) t^{ - \frac{3 - \beta}{2 \beta} }  .
\end{equation*}
Therefore we see that
\begin{equation}\label{eq4r}
\sup_{0  < t <T} \{ t^{ \frac{1}{2} } \| v (t) \|_{L^{\infty}} \} \leq C (T) < + \infty .
\end{equation}
Now we derive \eqref{eq111}. Using \eqref{eq26}, we have 
\begin{multline}\label{eq411}
\| v (  t ) \|_{L^\infty} \leq \| \mathrm{e}^{- t L} a \|_{L^\infty} \\
+ C \sup_{\phi \in [ C_0^\infty ]^3, { \ }\| \phi \|_{L^1} \leq 1} \bigg| \int_0^t \dual{  \mathrm{e}^{- \frac{(t-s)}{2} L} P (v , \nabla ) v  , \mathrm{e}^{- \frac{(t-s)}{2} L^*} P \phi  } d s \bigg|\\
+ C\sup_{\phi \in [ C_0^\infty ]^3, { \ }\| \phi \|_{L^1} \leq 1} \bigg| \int_0^t \dual{  \mathrm{e}^{- \frac{(t-s)}{2} L} P \{ (v , \nabla ) w + ( w , \nabla ) v \} , \mathrm{e}^{- \frac{(t-s)}{2} L^*} P \phi  } d s \bigg| .
\end{multline}
By \eqref{eq23}, \eqref{eq24}, \eqref{p6}, and \eqref{e16}, we check that
\begin{multline}\label{eq412}
\sup_{\phi \in [ C_0^\infty ]^3, { \ }\| \phi \|_{L^1} \leq 1} \bigg| \int_0^t \dual{  \mathrm{e}^{- \frac{(t-s)}{2} L} P (v , \nabla ) v  , \mathrm{e}^{- \frac{(t-s)}{2} L^*} P \phi  } d s \bigg|\\
\leq C \sup_{ \| \phi \|_{L^1} \leq 1} \int_0^t \| \mathrm{e}^{- \frac{(t-s)}{2} L} P  (v , \nabla ) v \|_{L^{\ell , \infty}} \| \mathrm{e}^{- \frac{(t-s)}{2} L^*} P \phi  \|_{L^{\ell' ,1}} d s\\
\leq C \int_0^t  \frac{ \| v \otimes v \|_{L^{ \frac{p }{ 2} , \infty} }}{(t - s)^{\frac{1}{2} + \frac{3}{2} \left( \frac{2}{ p} - \frac{1}{\ell} \right) + \frac{3}{2} \left( 1  - \frac{1}{\ell'} \right) }} d s\\
\leq C \int_0^t \frac{(s^\frac{p -3}{2 p} \| v ( s ) \|_{L^{p,\infty}} )^2}{s^{\color{red}{\frac{p -3}{p}}} ( t - s )^\frac{p + 6}{2 p} } d s \leq C t^{- \frac{1}{2} } .
\end{multline}
{\color{red}{Note that $(p+6)/{2p} < 1$ since $6 < p$.}}

\noindent Fix $\ell > 3 \beta p/(3 p + 3 \beta - \beta p)$. Note that $3 \beta p/(3 p + 3 \beta - \beta p) >3/2$. From \eqref{eq23}, \eqref{p6}, and \eqref{e16}, we see that
\begin{multline}\label{eq413}
\sup_{\phi \in [ C_0^\infty ]^3, { \ }\| \phi \|_{L^1} \leq 1} \bigg| \int_0^t \dual{  \mathrm{e}^{- \frac{(t-s)}{2} L} P \{ (v , \nabla ) w + (w , \nabla ) v \} , \mathrm{e}^{- \frac{(t-s)}{2} L^*} P \phi  } d s \bigg|\\
\leq C \sup_{ \| \phi \|_{L^1} \leq 1} \int_0^t \| \mathrm{e}^{- \frac{(t-s)}{2} L} P \{ (v , \nabla ) w + (w , \nabla ) v \} \|_{L^{\ell , \infty}} \| \mathrm{e}^{- \frac{(t-s)}{2} L^*} P \phi  \|_{L^{\ell' ,1}} d s\\
\leq C \int_0^t  \frac{ \| v \otimes w + w \otimes v \|_{L^{\frac{3 p \beta }{3 p + 3 \beta - \beta p} , \infty}} }{(t - s)^{\frac{1}{2} + \frac{3}{2} \left( \frac{3 p + 3 \beta - \beta p}{3 \beta p} - \frac{1}{\ell} \right) + \frac{3}{2} \left( 1  - \frac{1}{\ell'} \right) }} d s\\
\leq C \int_0^t  \frac{ ( s^\frac{p -3}{2 p}\| v (s) \|_{L^{p,\infty}} ) \| w \|_{L^{\frac{3 \beta }{3 - \beta} , \infty}}}{ s^\frac{p - 3}{2 p} ( t - s)^{\frac{3}{2} \left( \frac{1}{p} + \frac{1}{\beta} \right)} }    d s \leq C t^{- \frac{3 - \beta }{ 2 \beta} } .
\end{multline}
Combining \eqref{eq411}-\eqref{eq413}, we see that
\begin{equation}
\| v (t) \|_{L^\infty} \leq C t^{ - \frac{1}{2} } + C t^{- \frac{3 - \beta }{2 \beta }} .\label{eq414}
\end{equation}
Therefore we obtain \eqref{eq111}.

Finally, we prove the assertions $(\mathrm{i})$-$(\mathrm{ii})$ of Theorem \ref{thm111}. From Proposition \ref{prop41}, we deduce $(\mathrm{ii})$. Now we attack $(\mathrm{i})$. Assume that $a \in L^{\alpha , \infty}_\sigma ( \Omega )$ for some $3p/(p +3) < \alpha <3$. Let $\delta_1$ be the constant appearing in Theorem \ref{thm110}. Assume that $\| a \|_{L^{3,\infty}} < \delta_1$ and that $\| w \|_{L^{3,\infty}} < \delta_1$. From Theorem \ref{thm110}, we see that
\begin{equation*}
\sup_{t>0}\{ t^{\frac{3}{2}\left( \frac{1}{\alpha} - \frac{1}{p} \right) } \| v (t) \|_{L^{p , \infty} } \} \leq \text{ Const.} < + \infty .
\end{equation*}
Using the same argument to derive \eqref{eq414} with $\| \mathrm{e}^{- t L} a \|_{L^\infty} \leq C t^{-\frac{3}{2 \alpha}} \| a \|_{L^{ \alpha , \infty } }$, we have \eqref{eq00}. Therefore Theorem \ref{thm111} is proved.
\end{proof}

\section{Global-in-Time Strong Solution}\label{sect5}

In this section we show the existence of a unique global-in-time strong solution of the system \eqref{eq15} when $v_0 \in L^{3,\infty}_{0 , \sigma} ( \Omega)$ and both $\| v_0 \|_{L^{3,\infty}}$ and $\| w \|_{L^{3,\infty}}$ are sufficiently small.  

Throughout this section we assume that $w$ is as in Assumption \ref{assC} and that $|u_\infty | \leq \gamma$ for some $\gamma > 0$. Let $3 \beta /(2 \beta - 3)< p < \infty$ and $1 < p' < 3/2$ such that $1/p+1/{p'}=1$. Fix $p$ and $p'$. The symbols $L$ and $L^*$ represent the two linear operators defined by Subsection \ref{subsec23}. Let us prove Theorem \ref{thm112}.
\begin{proof}[Proof of Theorem \ref{thm112}]
Let $\delta_2$ be the constant appearing in Theorem \ref{thm111}. Assume that $a \in L^{3,\infty}_{ 0 , \sigma} (\Omega)$, $\| w \|_{L^{ 3 , \infty } } < \delta_2$, and $\| a \|_{L^{3 ,\infty}} < \delta_2$. Let $v$ be the global-in-time mild solution of the system \eqref{eq15} with the initial datum $v_0 = a$, obtained by Theorem \ref{thm111}. Now we show that $v$ is a strong solution of \eqref{eq15}. From the assertion $(\mathrm{ii})$ of Theorem \ref{thm111}, we see that $\nabla v (t) \in L^{3 , \infty}_0 ( \Omega )$ for each fixed $t>0$. Since $w$ satisfies Assumption \ref{assC} and $v ( t ) \in L^{3,\infty}_{0 , \sigma} ( \Omega ) \cap L_\sigma^{p , \infty} ( \Omega )$, if follows from Lemma \ref{lem28} that for each fixed $t >0$
\begin{equation*}
P \{ (v , \nabla ) v + (v , \nabla ) w + (w , \nabla ) v \} ( t ) \in L^{3, \infty}_{0, \sigma} ( \Omega  ).
\end{equation*}
Note that $f g \in L^{3,\infty}_0( \Omega )$ if $f \in L^{3,\infty}_0 (\Omega )$ and $g \in L^\infty ( \Omega ) $. Since $\mathrm{e}^{- t L}$ is an analytic $C_0$-semigroup on $L^{3 , \infty}_{0 , \sigma } ( \Omega )$, we shall prove that
$P\{ (v,\nabla) v + (v,\nabla) w + (w , \nabla) v \}$ is locally H\"{o}lder continuous with respect to time in the $L^{3,\infty}$-norm. Let $\varepsilon, T > 0$ such that $\varepsilon <T$. Fix $\varepsilon$ and $T$. We shall show that there is $0 < h_0 < 1$ such that
\begin{equation*}
P\{ (v,\nabla) v + (v,\nabla) w + (w , \nabla) v \} \in C^{h_0} ([\varepsilon , T] ; L^{3,\infty}_{ 0 , \sigma } (\Omega ) ) .
\end{equation*}
From \eqref{eq1233} and \eqref{eq111}, we see that
\begin{align*}
\sup_{0 <  t \leq T}\{  t^\frac{1}{2} \| \nabla v (t) \|_{L^{3,\infty}} \} + \sup_{0 <  t \leq T}\{  t^\frac{1}{2} \| v (t) \|_{L^{\infty}} \} \leq C (T) < + \infty,\\
\sup_{\varepsilon \leq t \leq T}\{  \| \nabla v (t) \|_{L^{3,\infty}} \} + \sup_{\varepsilon \leq t \leq T}\{ \| v (t) \|_{L^{\infty}} \} \leq C ( \varepsilon, T) < + \infty
\end{align*}

Let $t_1,t_2>0$ such that $\varepsilon \leq t_1 \leq t_2 \leq T$. It is easy to check that
\begin{multline*}
\| P\{ (v(t_2),\nabla) v(t_2) + (v(t_2),\nabla) w + (w , \nabla) v(t_2) \} \\
- P\{ (v(t_1),\nabla) v(t_1) + (v(t_1),\nabla) w + (w , \nabla) v(t_1) \} \|_{L^{3,\infty}} \\
\leq C \| v (t_2) - v (t_1) \|_{L^\infty} \| \nabla v (t_2) \|_{L^{3 ,\infty}} + C \| v( t_1 ) \|_{L^\infty} \| \nabla ( v (t_2) - v(t_1) ) \|_{L^{3,\infty}} \\
+ C \| v( t_2 ) - v ( t_1 ) \|_{L^\infty} \| \nabla w \|_{L^{3,\infty}} + C \| w \|_{L^\infty} \| \nabla (v(t_2) - v (t_1)) \|_{L^{3,\infty}}\\
\leq C (\varepsilon , T) \{ \| v(t_2) - v(t_1) \|_{L^\infty} + \| \nabla ( v (t_2) - v (t_1) )\|_{L^{3,\infty}} \} .
\end{multline*}
Therefore we only have to show that there are $C_{h} = C_h(\varepsilon ,T) >0$ and $0 < h_1 < 0$ such that
\begin{align}
& \| v (t_2) - v (t_1) \|_{L^\infty} \leq C_h (t_2 - t_1 )^{h_1},\label{eq51}\\
& \| \nabla (v (t_2) - v (t_1) ) \|_{L^{3,\infty}} \leq C_h (t_2 - t_1)^{h_1} \label{eq52} .
\end{align}

We first attack \eqref{eq51}. Since
\begin{multline}\label{mf}
v (t_2) - v (t_1) = ( \mathrm{e}^{- (t_2 - t_1) L} - 1 ) \mathrm{e}^{- t_1 L} a\\
 + \int_0^{t_1} ( \mathrm{e}^{- (t_2-t_1) L} - 1 ) \mathrm{e}^{(t_1 - s)L} P\{ (v,\nabla) v + (v,\nabla) w + (w , \nabla) v \} (s) d s\\
 + \int_{t_1}^{t_2} \mathrm{e}^{- (t_2 - t_1) L} P\{ (v,\nabla) v + (v,\nabla) w + (w , \nabla) v \} ( s ) d s ,
\end{multline}
we use \eqref{eq26} and a duality argument to see that
\begin{equation*}
\|  v (t_2) - v (t_1 ) \|_{L^\infty} \leq C \{ H_1 + H_2 + H_3 + H_4 + H_5 \} (t_2 , t_1) .
\end{equation*}
Here
\begin{align*}
& H_1 (t_2 ,t_1):= \sup_{\phi \in [C_0^\infty ]^3, { \ }\| \phi \|_{L^1} \leq 1} | \dual{ \mathrm{e}^{- \frac{t_1}{2} L}  a, (\mathrm{e}^{- ( t_2 - t_1 )L^* } -1 ) \mathrm{e}^{- \frac{t_1}{2} L^*} P \phi } | ,\\
& H_2 (t_2,t_1) :=  \sup_{\| \phi \|_{L^1} \leq 1} \left| \int_0^{t_1} \dual{ (v , \nabla ) v , ( \mathrm{e}^{- (t_2 -t_1) L^*} - 1 ) \mathrm{e}^{ - (t_1 - s) L^*} P \phi} d s \right| , \\
& H_3 :=  \sup_{\| \phi \|_{L^1} \leq 1} \left| \int_0^{t_1} \dual{ (v , \nabla ) w +( w , \nabla ) v  , ( \mathrm{e}^{- (t_2 -t_1) L^*} - 1 ) \mathrm{e}^{ - (t_1 - s) L^*} P \phi} d s \right| , \\
& H_4 (t_2,t_1) :=  \sup_{\phi \in [C_0^\infty ]^3, { \ }\| \phi \|_{L^1} \leq 1} \left| \int_{t_1}^{t_2} \dual{ (v , \nabla ) v , \mathrm{e}^{ - (t_2 - s) L^*} P \phi} d s \right| ,\\
& H_5 (t_2,t_1) :=  \sup_{\phi \in [C_0^\infty ]^3, { \ }\| \phi \|_{L^1} \leq 1} \left| \int_{t_1}^{t_2} \dual{ (v , \nabla ) w +( w , \nabla ) v  , \mathrm{e}^{ - (t_2 - s) L^*} P \phi} d s \right| .
\end{align*}
Using \eqref{eq23}, \eqref{eq24}, \eqref{e1}, and \eqref{6f}, we see that
\begin{align*}
H_1 (t_2 ,t_1) & \leq C \| \mathrm{e}^{- \frac{t_1}{2} L} a \|_{L^{3,\infty}} \sup_{\| \phi \|_{L^1} \leq 1} \| ( \mathrm{e}^{- (t_2 - t_1) L^*} - 1) \mathrm{e}^{- \frac{t_1}{2} L^*} P \phi \|_{L^{\frac{3}{2} , 1}}\\
& \leq C \| a \|_{L^{3 , \infty}} (t_2 - t_1 )^\frac{1}{4} \mathrm{e}^{\frac{t_1}{2} + (t_2 -t_1) } t_1^{ -\frac{3}{4} }\\
& \leq C ( \varepsilon , T) (t_2 - t_1 )^\frac{1}{4} .
\end{align*}
Applying \eqref{eq23}, \eqref{eq24}, \eqref{e1}, and \eqref{6f}, we observe that
\begin{multline*}
H_2 (t_2 ,t_1) \leq C \sup_{\| \phi \|_{L^1} \leq 1} \bigg\{ \int_0^{t_1} \frac{(s^\frac{p-3}{2 p} \| v (s) \|_{L^{p,\infty}} )(s^\frac{1}{2} \| \nabla v (s) \|_{L^{3,\infty}} )}{s^\frac{2 p -3}{2 p}}\\
\cdot \| (\mathrm{e}^{- (t_2 -t_1) L^*} -1) \mathrm{e}^{- (t_1 -s)L^*} P \phi \|_{L^{\frac{3 p }{2 p - 3} ,1 }  } d s \bigg\}\\
\leq C(T) \int_0^{t_1} \frac{1}{s^\frac{2 p -3}{2 p}} (t_2 - t_1)^\frac{1}{2 p} \mathrm{e}^{ \frac{t_1 -s}{2} + (t_2-t_1)} (t_2 - t_1)^{- \frac{2 p +4}{2 p} } d s\\
\leq C (\varepsilon , T) (t_2 - t_1)^\frac{1}{2 p} .
\end{multline*}
Similarly, we check that
\begin{multline*}
H_3 (t_2 ,t_1) \leq \\
C \sup_{\| \phi \|_{L^1} \leq 1} \bigg\{ \int_0^{t_1} \frac{ (s^\frac{1}{2} \| v (s) \|_{L^{\infty}}) \| \nabla w \|_{L^{3 ,\infty}} + \| w \|_{L^\infty} (s^\frac{1}{2} \| \nabla v (s) \|_{L^{3,\infty}} ) }{s^\frac{1}{2}} \\
\cdot \| (\mathrm{e}^{- (t_2 -t_1) L^*} -1) \mathrm{e}^{- (t_1 -s)L^*} P \phi \|_{L^{\frac{3}{2} ,1 }  } d s \bigg\}\\
\leq C(T) \int_0^{t_1} \frac{1}{s^\frac{1}{2}} (t_2 - t_1)^\frac{1}{4} \mathrm{e}^{ \frac{t_1 -s}{2} + (t_2-t_1)} (t_1 - s)^{-\frac{1}{4} - \frac{1}{2}} d s\\
\leq C (\varepsilon , T) (t_2 - t_1)^\frac{1}{4} .
\end{multline*}
Applying \eqref{eq23}, \eqref{eq24}, and \eqref{p4}. we observe that
\begin{multline*}
H_4 (t_2 ,t_1) \leq C \sup_{\| \phi \|_{L^1} \leq 1} \bigg\{ \\
\int_{t_1}^{t_2} \frac{(s^\frac{p-3}{2 p} \| v (s) \|_{L^{p,\infty}} )(s^\frac{1}{2} \| \nabla v (s) \|_{L^{3,\infty}} )}{s^\frac{2 p -3}{2 p}} \| \mathrm{e}^{- (t_2 -s)L^*} P \phi \|_{L^{\frac{3 p }{2 p -3} ,1 }  } d s \bigg\}\\
\leq C(T) t_1^{- \frac{2 p -3}{2 p}} \int_{t_1}^{t_2} \frac{1}{ (t_2 -s)^\frac{p + 3}{2 p}} d s \leq C (\varepsilon , T) (t_2 - t_1)^\frac{p - 3}{2 p}
\end{multline*}
and that
\begin{multline*}
H_5 \leq C \sup_{\| \phi \|_{L^1} \leq 1} \bigg\{ \int_{t_1}^{t_2} \frac{ (s^\frac{1}{2} \| v (s) \|_{L^{\infty}})\| \nabla w \|_{L^{3 ,\infty}} + \| w \|_{L^\infty} (s^\frac{1}{2} \| \nabla v (s) \|_{L^{3,\infty}} ) }{s^\frac{1}{2}} \\
\cdot \| \mathrm{e}^{- (t_2 -s)L^*} P \phi \|_{L^{\frac{3}{2} ,1 }  } d s \bigg\}\\
\leq C(T) t_1^{- \frac{1}{2} }\int_{t_1}^{t_2} \frac{1}{ (t_2 - s)^\frac{1}{2} } d s\leq C (\varepsilon , T) (t_2 - t_1)^\frac{1}{2} .
\end{multline*}
Therefore we see \eqref{eq51}.

Next we derive \eqref{eq52}. Fix $j = 1,2,3$. Using \eqref{mf} and Lemma \ref{lem25}, we have
\begin{equation*}
\|  \partial_j v (t_2) - \partial_j v (t_1 ) \|_{L^{3 , \infty}} \leq C\{ H_6 + H_7 + H_8 + H_9 + H_{10}\} (t_2,t_1) .
\end{equation*}
Here
\begin{align*}
& H_6 (t_2 ,t_1):= \sup_{\phi \in [C_0^\infty ]^3, { \ }\| \phi \|_{L^{\frac{3}{2}, 1}} \leq 1} | \dual{ \mathrm{e}^{- \frac{t_1}{2} L}  a, (\mathrm{e}^{- ( t_2 - t_1 )L^* } -1 ) \mathrm{e}^{- \frac{t_1}{2} L^*} P  ( - \partial_j \phi ) } | ,\\
& H_7  :=  \sup_{\| \phi \|_{L^{\frac{3}{2}, 1}} \leq 1} \left| \int_0^{t_1} \dual{ (v , \nabla ) v , ( \mathrm{e}^{- (t_2 -t_1) L^*} - 1 ) \mathrm{e}^{ - (t_1 - s) L^*} P ( - \partial_j \phi ) } d s \right| , \\
& H_9 (t_2,t_1) :=  \sup_{\phi \in [C_0^\infty ]^3, { \ }\| \phi \|_{L^{\frac{3}{2} ,1}} \leq 1} \left| \int_{t_1}^{t_2} \dual{ (v , \nabla ) v , \mathrm{e}^{ - (t_2 - s) L^*} P  ( - \partial_j \phi )} d s \right| ,\\
& H_{10} (t_2,t_1) :=  \sup_{\| \phi \|_{L^{\frac{3}{2} ,1}} \leq 1} \left| \int_{t_1}^{t_2} \dual{ (v , \nabla ) w +( w , \nabla ) v  , \mathrm{e}^{ - (t_2 - s) L^*} P ( - \partial_j \phi ) } d s \right| ,
\end{align*}
\begin{multline*}
H_8 (t_2,t_1) :=  \sup_{\phi \in [C_0^\infty ]^3, { \ }\| \phi \|_{L^{\frac{3}{2} ,1}} \leq 1} \bigg| \\
\int_0^{t_1} \dual{ (v , \nabla ) w +( w , \nabla ) v  , ( \mathrm{e}^{- (t_2 -t_1) L^*} - 1 ) \mathrm{e}^{ - (t_1 - s) L^*} P ( - \partial_j \phi ) } d s \bigg| .
\end{multline*}
Using \eqref{eq23}, \eqref{eq24}, \eqref{e1}, and \eqref{6d}, we see that
\begin{align*}
H_6 (t_2 ,t_1) & \leq C \| \mathrm{e}^{- \frac{t_1}{2} L} a \|_{L^{3,\infty}} \sup_{\| \phi \|_{L^{ \frac{3}{2} , 1}} \leq 1} \| ( \mathrm{e}^{- (t_2 - t_1) L^*} - 1) \mathrm{e}^{- \frac{t_1}{2} L^*} P \partial_j \phi \|_{L^{\frac{3}{2} , 1}}\\
& \leq C \| a \|_{L^{3 , \infty}} (t_2 - t_1 )^\frac{1}{4} \mathrm{e}^{\frac{t_1}{2} + (t_2 -t_1) } t_1^{- \frac{3}{4} }\\
& \leq C (\varepsilon , T) (t_2 - t_1 )^\frac{1}{4} .
\end{align*}
Applying \eqref{eq23}, \eqref{eq24}, \eqref{e1}, and \eqref{6d}, we check that
\begin{multline*}
H_7 (t_2 ,t_1) \leq C \sup_{\| \phi \|_{L^{\frac{3}{2},1}} \leq 1} \bigg\{ \int_0^{t_1} \frac{(s^\frac{p-3}{2 p} \| v (s) \|_{L^{p,\infty}} )(s^\frac{1}{2} \| \nabla v (s) \|_{L^{3,\infty}} )}{s^\frac{2 p -3}{2 p}} \\
\cdot \| (\mathrm{e}^{- (t_2 -t_1) L^*} -1) \mathrm{e}^{- (t_1 -s)L^*} P \partial_j \phi \|_{L^{\frac{3 p }{2 p - 3} ,1 }  } d s \bigg\}\\
\leq C(T) \int_0^{t_1} \frac{1}{s^\frac{p -3}{2 p}} (t_2 - t_1)^\frac{1}{2 p} \mathrm{e}^{ \frac{t_1 -s}{2} + (t_2-t_1)} (t_2 - t_1)^{- \frac{p + 4}{2 p}} d s\\
\leq C (\varepsilon , T) (t_2 - t_1)^\frac{1}{2 p} .
\end{multline*}
Similarly, we see that
\begin{multline*}
H_8 (t_2 ,t_1) \leq \\
C \sup_{\| \phi \|_{L^{\frac{3}{2} , 1}} \leq 1} \bigg\{ 
\int_0^{t_1} \frac{ (s^\frac{1}{2} \| v (s) \|_{L^{\infty}})\| \nabla w \|_{L^{3 ,\infty}} + \| w \|_{L^\infty} (s^\frac{1}{2} \| \nabla v (s) \|_{L^{3,\infty}} ) }{s^\frac{1}{2}} \\
\cdot \| (\mathrm{e}^{- (t_2 -t_1) L^*} -1) \mathrm{e}^{- (t_1 -s)L^*} P \phi \|_{L^{\frac{3}{2} ,1 }  } d s \bigg\}\\
\leq C(T) \int_0^{t_1} \frac{1}{s^\frac{1}{2}} (t_2 - t_1)^\frac{1}{4} \mathrm{e}^{ \frac{t_1 -s}{2} + (t_2-t_1)} ( t_1 - s)^{-\frac{3}{4}} d s\\
\leq C (\varepsilon , T) (t_2 - t_1)^\frac{1}{4} .
\end{multline*}
Using \eqref{eq23}, \eqref{eq24}, and, \eqref{e16}, we check that
\begin{multline*}
H_9 (t_2 ,t_1) \leq C \sup_{\| \phi \|_{L^{\frac{3}{2} , 1}} \leq 1} \bigg\{ \int_{t_1}^{t_2} \frac{(s^\frac{p-3}{2 p} \| v (s) \|_{L^{p,\infty}} )(s^\frac{1}{2} \| \nabla v (s) \|_{L^{3,\infty}} )}{s^\frac{2 p -3}{2 p}}\\
\cdot \| \mathrm{e}^{- (t_2 -s)L^*} P \partial_j \phi \|_{L^{\frac{3 p }{2 p -3} ,1 }  } d s \bigg\}\\
\leq C(T) t_1^{- \frac{2 p -3}{2 p}} \int_{t_1}^{t_2} \frac{1}{ (t_2 -s)^\frac{p + 3}{2 p}} d s\leq C (\varepsilon , T) (t_2 - t_1)^\frac{p - 3}{2 p} 
\end{multline*}
and that
\begin{multline*}
H_{10}(t_2,t_1) \leq \\
C \sup_{\| \phi \|_{L^{\frac{3}{2} , 1}} \leq 1} \bigg\{ \int_{t_1}^{t_2} \frac{ (s^\frac{1}{2} \| v (s) \|_{L^{\infty}})\| \nabla w \|_{L^{3 ,\infty}} + \| w \|_{L^\infty} (s^\frac{1}{2} \| \nabla v (s) \|_{L^{3,\infty}} ) }{s^\frac{1}{2}} \\
\cdot \| \mathrm{e}^{- (t_2 -s)L^*} P \partial_j \phi \|_{L^{\frac{3}{2} ,1 }  } d s \bigg\}\\
\leq C(T) t_1^{- \frac{1}{2} }\int_{t_1}^{t_2} \frac{1}{ (t_2 - s)^\frac{1}{2} } d s \leq C (\varepsilon , T) (t_2 - t_1)^\frac{1}{2} .
\end{multline*}
Therefore we see \eqref{eq52}.

From \eqref{eq51} and \eqref{eq52}, we find that $P\{ (v, \nabla ) v + ( v , \nabla ) w + (w , \nabla )v  \}$ is H\"{o}lder continuous on $[\varepsilon , T]$ in the $L^{3,\infty}$-norm. Since $\varepsilon$ and $T$ are arbitrary, it follows from the analyticity of the semigroup $\mathrm{e}^{- t L}$ to see that
\begin{equation*}
v \in C ( (0 . \infty ) ; L^{3,\infty}_{ 0 , \sigma } (\Omega ) ) \cap C ( (0 , \infty) ; D (L) ) \cap C^1 ( (0, \infty ) ; L^{3,\infty}_{ 0 , \sigma } ( \Omega ) )
\end{equation*}
and that $v$ is a strong solution of the system \eqref{eq15} with the initial datum $v_0 = a$. Therefore Theorem \ref{thm112} is proved. Note that $D (L) = D (\mathscr{L}_{3,\infty})$ (see Lemma \ref{lem212}).
\end{proof}

\section{Appendix: Characterization of the Lorentz Norm}\label{sect6}
In the Appendix, we characterize the Lorentz norm. Let $n, m \in \mathbb{N}$, and let $U \subset \mathbb{R}^n$ be a domain. Let $1 < q ,q' < \infty$ and $1 \leq r ,r' \leq \infty$ such that $1/q +1/{q'}=1$ and $1/r + 1/{r'}=1$. We define the Lorentz space $L^{q,r}(U)$ as follows:
\begin{align*}
& L^{q,r}  ( U ) := \{ f \in L^1 (U) + L^\infty ( U) ; { \ }\| f \|_{L^{q,r} (U)} < +\infty  \} , \\
& \| f \|_{L^{q,r}(U)} := 
\begin{cases}
\left( \int_0^\infty (t^{1/q} f^{**} (t) )^r \frac{d t}{t} \right)^{1/r} & \text{ if } 1 \leq r < + \infty, \\
\sup_{t >0} \{ t^{1/q} f^{**} ( t ) \}  & \text{ if } r = \infty ,
\end{cases}\\
&f^* (t) := \inf \{ \sigma >0  ; { \ }  \mu \{ x \in U; { \ } | f (x) | > \sigma \} \leq t \}, { \ } t \geq 0 ,\\
&f^{**}(t) := t^{-1}\int_0^t f^*(s) d s, { \ }t \geq 0.
\end{align*}
Here $\mu \{ \cdot \}$ denotes the $n$-dimensional Lebesgue measure. \\For $f = (f^1, \cdots,f^m) \in [L^{q,r}( U )]^m$ and $g \in [L^{q',r'}(U)]^m$,
\begin{align*}
& \dual{f , g} := \int_U f (x )\cdot g (x) d x,\\
& \| f \|_{L^{q,r} ( U )} \equiv \| f \|_{X^{q,r}( U )} : = \| \sqrt{(f^1)^2 + (f^2)^2 + \cdots + (f^m)^2}  \|_{L^{q,r}( U )},\\
& \| f \|_{Y^{q,r}( U )} := \max_{1 \leq j \leq m} \{ \| f^j  \|_{L^{q,r}( U )} \},\\
& \| f \|_{Z^{q,r}( U )} := \sup_{\phi \in [ L^{q',r'}( U ) ]^m,{ \ } \| \phi \|_{L^{q',r'}( U )} \leq 1 } | \dual{f , \phi }|.
\end{align*}
Here
\begin{equation*}
\dual{f , \phi} := \int_U f (x ) \cdot \phi (x) d x .
\end{equation*}
The purpose of the Appendix is to prove the following proposition.
\begin{proposition}\label{prop61}
There is $C_\vartriangle = C_\vartriangle (m, q,r) >0$ such that \\for $f = (f^1,\cdots,f^m) \in [L^{q,r}(U)]^m$
\begin{align}
& \| f \|_{X^{q,r}( U )} \leq C_\vartriangle \| f \|_{Y^{q,r}( U )},\label{eq61}\\
& \| f \|_{Y^{q,r}( U )} \leq \| f \|_{Z^{q,r} ( U )},\label{eq62}\\
& \| f \|_{Z^{q,r}( U )} \leq \| f \|_{X^{q,r}( U )}\label{eq63}.
\end{align}
Moreover, assume in addition that $[C_0^\infty (U)]^m$ is dense in $[ L^{q',r'}(U) ]^m$. Then
\begin{equation*}
\| f \|_{Z^{q,r}(U)} \leq C_\vartriangle \sup_{\phi \in [C_0^\infty (U) ]^m, { \ }\| \phi \|_{L^{q',r'}} \leq 1} | \dual{f , \phi }| .
\end{equation*}
\end{proposition}
To prove Proposition \ref{prop61}, we prepare the three lemmas. From \cite[Chapter 1]{BL76} and \cite[IV Lemma 4.5]{BS88}, we have
\begin{lemma}\label{lem62}
$(\mathrm{i})$ Assume that $f \in L^{q,r} (U)$. Then for $t>0$
\begin{align*}
f^* ( t )  \leq & f^{**} (t),\\
\left( \int_0^\infty (t^{1/q} f^{**} (t) )^r \frac{d t}{t} \right)^{1/r} \leq & \frac{q}{q-1}\left( \int_0^\infty (t^{1/q} f^{*} (t) )^r \frac{d t}{t} \right)^{1/r} , \text{ if }1 \leq r <  \infty,\\
\sup_{t >0} \{ t^{1/q} f^{**} ( t ) \} \leq & \frac{q}{ q -1} \sup_{t >0} \{ t^{1/q} f^{*} ( t ) \}, \text{ if } r = \infty .
\end{align*}
\noindent $( \mathrm{ii})$ For all $f \in [ L^{q,r} ( U ) ]^m$ and $g \in [ L^{q' , r'} ( U ) ]^m$,
\begin{equation*}
| \dual{ f , g } | \leq \| f \|_{L^{q,r} (U)} \| g \|_{L^{q',r'} (U) } .  
\end{equation*}
\end{lemma}
Since the $n$-dimensional Lebesgue measure is nonatomic (\cite[Chapter 21]{Fre03}), we apply \cite[II Theorem 2.7 and IV Theorem 4.7]{BS88} to obtain 
\begin{lemma}\label{lem63}
For every $f \in L^{q,r}(\mathbb{R}^n )$
\begin{equation*}
\| f \|_{L^{q,r}(\mathbb{R}^n)} \leq \sup_{\Phi \in L^1_{loc} (\mathbb{R}^n), { \ }\| \Phi \|_{L^{q',r'}( \mathbb{R}^n )} \leq 1 } | \dual{f , \Phi }| \leq \| f \|_{L^{q,r} (\mathbb{R}^n)} .
\end{equation*}
\end{lemma}
Next we prove
\begin{lemma}\label{lem64}
For every $f \in L^{q,r}(U )$ 
\begin{equation*}
\| f \|_{L^{q,r}(U )} \leq \sup_{\phi \in L^{q',r'}(U ), { \ }\| \phi \|_{L^{q',r'}( U )} \leq 1 } | \dual{f , \phi }| \leq \| f \|_{L^{q,r} (U )} .
\end{equation*}
\end{lemma}
\begin{proof}[Proof of Lemma \ref{lem64}]
For $U \subset \mathbb{R}^n$,
\begin{equation*}
1_U := 
\begin{cases}
1, & x \in U\\
0, & x \in \mathbb{R}^n \setminus U .
\end{cases}
\end{equation*}
Fix $f \in L^{q,r} ( U )$. By the definition of $\| \cdot \|_{L^{q,r}}$ and Lemma \ref{lem63}, we see that
\begin{align*}
\| f \|_{L^{q,r } (U )} \leq & \| f 1_U \|_{L^{q,r} (\mathbb{R}^n)}\\
= & \sup_{\Phi \in L^1_{loc} (\mathbb{R}^n), { \ }\| \Phi \|_{L^{q',r'}( \mathbb{R}^n )} \leq 1 } | \dual{f 1_U , \Phi }| = \text{ (RHS)} .
\end{align*}
It is clear that
\begin{multline*}
\text{ (RHS) }= \sup_{\Phi \in L^1_{loc} (\mathbb{R}^n), { \ }\| \Phi \|_{L^{q',r'}( \mathbb{R}^n )} \leq 1 } \left|  \int_{\mathbb{R}^n} f 1_U (x) \Phi (x) d x \right| \\
\sup_{\Phi_1, \Phi_2 \in L^1_{loc} (\mathbb{R}^n), { \ } \| \Phi_1 \|_{L^{q',r'}( U )} + \| \Phi_2 \|_{L^{q',r'}(\mathbb{R}^n \setminus U )} \leq 1 } \left|  \int_{\mathbb{R}^n} f 1_U (x) (\Phi_1 (x) + \Phi_2 (x) ) d x \right| \\
= \sup_{\Phi_1, \Phi_2 \in L^1_{loc} (\mathbb{R}^n), { \ } \| \Phi_1 \|_{L^{q',r'}( U )} + \| \Phi_2 \|_{L^{q',r'}(\mathbb{R}^n \setminus U )} \leq 1 } \left|  \int_{U} f (x) \Phi_1 (x) d x \right| \\
= \sup_{\Phi_1 \in L^1_{loc} (\mathbb{R}^n), { \ } \| \Phi_1 \|_{L^{q',r'}( U )} \leq 1 } \left|  \int_{U} f (x) \Phi_1 (x) d x \right| .
\end{multline*}
Therefore we have
\begin{equation*}
\| f \|_{L^{q,r}(U )} \leq \sup_{\phi \in L^{q',r'}(U ),{ \ } \| \phi \|_{L^{q',r'}( U )} \leq 1 } | \dual{f , \phi }|  .
\end{equation*}
By Lemma \ref{lem62}, we see that
\begin{equation*}
\sup_{\phi \in L^{q',r'}(U ), { \ }\| \phi \|_{L^{q',r'}( U )} \leq 1 } | \dual{f , \phi }| \leq \| f \|_{L^{q,r} (U )} .
\end{equation*}
\end{proof}
Let us attack Proposition \ref{prop61}.
\begin{proof}[Proof of Proposition \ref{prop61}]
We first show \eqref{eq61}. Since
\begin{multline*}
\mu \{ x \in U; \sqrt{ (f^1 (x))^2 + \cdots ( f^m(x))^2 } > t \} \\
\leq \mu \{ x \in U; |f (x)| > t/\sqrt{m} \}+ \cdots + \mu \{ x \in U; |f^m(x)| > t/\sqrt{m} \},
\end{multline*}
we apply Lemma \ref{lem62} and the definition of $\| \cdot \|_{L^{q,r}}$ to check that
\begin{align*}
\| f \|_{X^{q,r}( U )} \leq & C (m,q,r) (\| f^1 \|_{L^{q,r}} + \cdots + \| f^m \|_{L^{q,r}})\\
\leq & C(m,q,r)\| f \|_{Y^{q,r}}.
\end{align*}
Therefore we see \eqref{eq61}.

Secondly, we derive \eqref{eq62}. Let $f = (f^1,f^2,\cdots,f^m) \in [L^{q,r}(U)]^m$. For each $1 \leq j \leq m$, we check that
\begin{multline*}
\| f^j \|_{L^{q,r}(U )} \leq \sup_{\phi^j \in L^{q',r'}(U ),{ \ } \| \phi^j \|_{L^{q',r'}( U )} \leq 1 } | \dual{f^j , \phi^j }|  \\
= \sup_{\phi = (\phi^1, \cdots, \phi^j, \cdots, \phi^m ) \in [L^{q',r'}(U )]^m,{ \ } \| \phi \|_{L^{q',r'}( U )} \leq 1 } | \dual{f^j , \phi^j }| \\
 \leq \sup_{\phi = (\phi^1, \cdots, \phi^j, \cdots, \phi^m ) \in [L^{q',r'}(U )]^m, { \ }\| \phi \|_{L^{q',r'}( U )} \leq 1 } | \dual{f , \phi }|  .
\end{multline*}
Note: Since $\dual{f , \phi} = \dual{f^j , \phi^j }$ when $\phi = (0, \cdots, 0, \phi^j, 0 ,\cdots, 0)$, it follows that
\begin{multline*}
\sup_{\phi = (\phi^1, \cdots, \phi^j, \cdots, \phi^m ) \in [L^{q',r'}(U )]^m,{ \ } \| \phi \|_{L^{q',r'}( U )} \leq 1 } | \dual{f^j , \phi^j }| \\
 \leq \sup_{\phi = (\phi^1, \cdots, \phi^j, \cdots, \phi^m ) \in [L^{q',r'}(U )]^m,{ \ } \| \phi \|_{L^{q',r'}( U )} \leq 1 } | \dual{f , \phi }|  .
\end{multline*}
Therefore we see that
\begin{equation*}
\max_{1 \leq j \leq m} \{ \| f^j \|_{L^{q,r}(U )} \} \leq \sup_{\phi \in [L^{q',r'}(U )]^m,{ \ } \| \phi \|_{L^{q',r'}( U )} \leq 1 } | \dual{f , \phi }|,
\end{equation*}
which is \eqref{eq62}. From Lemma \ref{lem62}, we have \eqref{eq63}. 
\end{proof}

\begin{center}
Acknowledgment
\end{center}
{\color{red}{The author gratefully acknowledges the precious comments of Professor Hideo Kozono. The author would like to thank Professor Toshiaki Hishida for valuable discussions about the paper \cite{HS09}. He is also grateful to the anonymous referees for their valuable comments to improve this paper.

A part of this work was supported by the Japanese-German Graduate Externship-Mathematical Fluid Dynamics. A part of this work was done when the author of the paper was an assistant professor at Waseda University (Japan).}}





\end{document}